\documentclass[a4paper,10pt]{article}
\usepackage{amscd}
\usepackage{textcomp}
\usepackage{float}
\usepackage{amsmath}
\usepackage{bbding}
\usepackage{amsfonts}
\usepackage{subfigure}
\usepackage{cite}
\usepackage{makecell}
\usepackage{mathrsfs}
\usepackage[driverfallback=dvipdfm,
           pdfstartview=FitH,
           colorlinks=true,
             pdfborder=001,
           citecolor=blue]{hyperref}
\usepackage{amssymb}
\usepackage{graphicx,latexsym,amsmath,amssymb,amsthm,enumerate}
\allowdisplaybreaks[4]
9.4in\textwidth 6.5in \hoffset -2.0cm \DeclareFontEncoding{OT2}{}{}
\DeclareTextSymbol{\cyrsftsn}{OT2}{126}
\DeclareTextSymbol{\textnumero}{OT2}{125}

\theoremstyle{definition}
\newtheorem{theorem}{Theorem}[section]
\newtheorem{lemma}{Lemma}[section]
\newtheorem{corollary}{Corollary}[section]
\newtheorem{proposition}{Proposition}[section]
\newtheorem{definition}{Definition}[section]
\newtheorem{remark}{Remark}[section]

\begin{document}

\title{{\Large\bf Finite-time and fixed-time consensus control of multi-agent systems driven by parabolic partial differential equations}\thanks{This work was supported by  the National Natural Science Foundation of China (Nos.11671282 and 12171339).}}
\author{{Xu-hui Wang, Xue-song Li and  Nan-jing Huang\footnote{Corresponding author.  E-mail address: nanjinghuang@hotmail.com; njhuang@scu.edu.cn}  } \\
{\small\it Department of Mathematics, Sichuan University, Chengdu,
	Sichuan 610064, P.R. China}}
\date{ }
\maketitle
\begin{flushleft}
\hrulefill\\
\end{flushleft}
 {\bf Abstract}.
 This paper focuses on the study of the finite-time consensus (FTC) and fixed-time consensus (FXC) issues of multi-agent systems (MASs) driven by parabolic partial differential equations (PDEs). Compared with the study in the existing literature, the topic of FTC and FXC control is first embodied in MASs driven by parabolic PDEs. Based on the Lyapunov theorems, the FTC and FXC controllers are devised to ensure that the MASs converge to a stable state with external disturbance. Furthermore, we simplify the controllers to guarantee the FTC and FXC of MASs without external disturbance. Finally, two illustrative examples are given to verify the feasibility of controllers.
\begin{flushleft} \hrulefill \end{flushleft}

\section{Introduction}
\noindent

In recent decades, multi-agent systems (MASs) have aroused the extensive interest of researchers in many areas, especially in computer science, neural network, and cooperative control. The research on MASs provides necessary theoretical support for the practical application in engineering, such as control problems of complex systems \cite{Qin16}, UAV formation \cite{Avram18}, and target tracking \cite{Ren09}. It is well known that consensus problem is an essential research subject in MASs research, and its purpose is to make multiple agents achieve agreement under the control protocol. Because of its wide applications in modern industry, this part of the research has a long history. Reynolds \cite{Re87} first proposed the Boid model according to the behavioral characteristics of groups (such as sheep, cattle and school), which was the earliest MASs model. Later, Vicsek et al. \cite{Vi95} took a different approach and first established the model of MASs from the statistical point of view, which provided a powerful tool to study the group movement behavior. In 2004, Olfatid-Saber and Murray \cite{Ol04} proposed a general approach to solve the asymptotical consensus issue of MASs based on first-order integrators. Recently, asymptotical consensus control for MASs has inspired extensive research interest such as for second-order MASs \cite{LiH17}, for high-order MASs \cite{Hamed17}, for fractional-order MASs \cite{Wang19}, and so on. More literature about asymptotical consensus issues, see \cite{Li05,Ol07,Re07,Re08,Yu09,Wu18} and the references therein.

However, most of the results in previous works of literature have a limitation, that is, the terminal time cannot be controlled within a limited range. In practical applications, we often need a faster convergence speed and hope MASs to reach consensus in finite time \cite{Yang13,Liu15,WangL19}. Therefore, it is vital to investigate the finite-time consensus (FTC) issue of MASs. In order to make the systems reach the equilibrium point in finite time, finite-time control theory was proposed \cite{Bh00}. Cortes \cite{Co06} extended the theory of nonsmooth stability analysis to study the FTC problem of MASs and proposed a discontinuous FTC controller. Wang and Xiao \cite{Wa10} proposed a state feedback controller  to solve the FTC issue of first-order MASs. Li et al. \cite{Li11} studied second-order MASs with external disturbance and designed a non-smooth FTC control protocol using the backstepping method. Guan et al. \cite{Guan12} used the homogeneity theory to study the FTC problem of second-order MASs with a leader in fixed and switched topologies. For high-order systems, Fu et al. \cite{Fu16} designed an observer-based FTC protocol. Recently, Duan et al. \cite{Duan19} solved the bipartite FTC problem for heterogeneous MASs using finite-time observer based on event triggering mechanism. However, the initial value may not be known and used in advance, which makes the FTC technology unworkable. In order to make the terminal time independent of the initial condition, the fixed-time control theory was also proposed \cite{Po12}. Parsegov et al. \cite{Pa13} applied the fixed-time control theory to first-order linear MASs and designed a distributed controller to achieve the fixed-time consensus (FXC). In addition, Defoort et al. \cite{De15} considered the FXC problem of first-order nonlinear leader-following MASs. When the agent's dynamic behavior is nonlinear to some extent, the leader's state can be tracked by the followers' state via designing corresponding fixed-time control protocols. Zuo \cite{Zuo15} studied the FXC problem of second-order leader-following MASs and designed a FXC protocol using the sliding mode control method to make the systems achieve consensus in fixed-time. Zuo et al. \cite{Zuo16} analyzed MASs with interference and designed a robust FXC protocol to enhance the systems' robustness. Zuo et al. \cite{Zuo18} reviewed the theoretical research results of fixed-time cooperative control in recent years. Compared with the traditional asymptotical consensus control theory, the FTC and FXC control methods can achieve higher control accuracy, faster convergence speed, and better control effect. Recently,  various theoretical results and applications have been studied extensively for the FTC and FXC problems in the literature; for instance we refer the reader to \cite{Zh12, Zh14, Xi14, Li17, Tong18, Sa19, Liu19, LiuJ19, Hu19, Li22} and references therein.

It is worth to be noted that most of the current research results on the FTC and FXC are obtained under the assumption that dynamical systems of agents are driven by ordinary differential equations (ODEs).  However, in many cases, the dynamics system of agents is not only related to time but also affected by spatial location, such as the spread of infectious diseases, the vehicle formation, the traffic management of large-scale UAV systems, the urban traffic networks, the noise control, the temperature control in industrial production, and so on, see \cite{Ghods12,Ras19,Cas19,Wa11, Wa12, Ra81}. Thus, it is more suitable to use partial differential equations (PDEs) to model these dynamics systems than to use ODEs. Recently, Pilloni et al. \cite{Pi16} investigated the synchronization problem of MASs described by boundary-controlled heat equations and designed a consensus controller through a nonlinear sliding-mode method eliminating the influence of boundary disturbance. Yang et al. \cite{Ya17} took into account the synchronization problem of nonlinear spatio-temporal MASs based on boundary controllers and gave sufficient conditions to ensure systems' consistency. And Yang et al. \cite{Ya18} investigated the output synchronization problem for parabolic partial differential systems and designed the feedback controller according to output information to ensure the synchronization. Very recently, Wang et al.\cite{Wang21} considered a special kind of first-order hyperbolic systems and designed the boundary protocol to reach the FTC.  Nevertheless, in some practical cases, it is necessary to consider MASs based on parabolic PDEs with external disturbance having the following form:
\begin{eqnarray}\label{e1}
\begin{array}{llll}
y_{i,t}(x,t)=ky_{i,xx}(x,t)+u_i(x,t)+d_i(x,t),\\
\end{array}
\end{eqnarray}
where  $y_i(x,t) \in \mathbb{R}$ is the state variable of the $i$th agent, $u_i(x,t)\in \mathbb{R}$ is the protocol of the $i$th agent, $d_i(x,t)\in \mathbb{R}$ is unknown external disturbance, and $k>0$ is a constant.  Let us illustrate this fact with the following food web network example introduced by Wang and Wu \cite{Wu13}.

A food web can be described by a complex network model in which nodes represent each species.

In limited habitats, the distribution of species is usually uneven, with different spatial densities of species leading to corresponding different migration. The spatial distribution of species can be described as state variables of nodes in parabolic partial differential systems, which depends on both spatial and temporal variables.  In particular,  a simple food web network model can be obtained from \cite{Wu13} as follows:
\begin{eqnarray}\label{e2}
\begin{array}{llll}
w_{i,t}(x,t)=kw_{i,xx}(x,t)+u_i(x,t),\\
\end{array}
\end{eqnarray}
where  $w_i(x,t) \in \mathbb{R}$ is the spatial density of the $i$th species,  $u_i(x,t)\in \mathbb{R}$ is the control input of the $i$th species, and $k>0$ is the diffusion coefficient. Clearly,  MASs \eqref{e2} is a special case of MASs \eqref{e1} with $d_i(x,t)\equiv 0$. It should be noted that the state of nodes in MASs may not be accurately measured in general since they are often destroyed by various disturbances and noise. Inspired by this fact and taking into account the external disturbance, we consider MASs \eqref{e1} based on parabolic PDEs with external disturbance in this work.

The abovementioned observations show that it is very important for agents to complete all planned tasks in a limited time in engineering applications.  On the other hand, there is a lack of unified results addressing the FTC and FXC issues for MASs described by PDEs. Thus, how to design the controllers to solve the FTC and FXC issues for PDE systems is still an open problem, which motivates this study.

The overarching goal of this paper is to investigate the FTC and FXC issues for MASs described by parabolic PDEs with external disturbance. Compared with the study in the existing literature, the main contributions of this paper can be listed as follows.

 (1) A new model of  MASs driven by parabolic PDEs with external disturbance is proposed, which can be used to describe the heat process utilized in science, engineering and practical experience in chemical and biochemical engineering. It is the first time that the FTC and FXC distributed controllers are designed to satisfy different control effects in MASs with external disturbance.

 (2) Compared with other works on FTC and FXC control for MASs with ODEs, the main difficulties in studying MASs with PDEs are to construct the Lyapunov function, to estimate the state variables, and to design the controllers. In order to overcome such difficulties, we first construct the Lyapunov function in the integral form to avoid the estimation of multivariate Lyapunov functions. Because the state variables are correlated with $x$ and $t$, it is inevitable to estimate the state variables' first and second partial derivatives with respect to $x$. We then use the integration by parts method and Wirtinger's inequality to transform the estimations of the first and second partial derivatives of the state variables into the ones of the state variables, which can be estimated via some known inequalities in \cite{Ol04,Wa10,Shi20,Wu14}. Finally, we design the controllers using the Lyapunov method for PDEs \cite{Bh00,Po12}.

 (3) Different from the research on consensus control for MASs driven by PDEs, the control objective considered in this paper is to achieve consensus in finite time or fixed time, which is more difficult to use the information of neighboring agents to improve the convergence speed.  A systematic approach is first proposed to accelerate the convergence speed and improve the control effect by adding a new control coefficient to each controller.

The remainder of this work is structured as follows. The next section recalls some critical lemmas, and formulates the FTC and FXC problem. After that in Section III, we analyze the FTC and FXC of MASs with external disturbance via the controllers. Before summarizing the results in Section V,  we provide two simulation examples to demonstrate the feasibility of controllers in Section IV.

{\it Notations}. Let $\mathbb{R}=(-\infty,+\infty)$ and $\mathbb{R}^+=(0,+\infty)$. Assume that $\mathbb{R}^n$ and $\mathbb{R}^{n\times n}$  denote respectively the $n$ dimensional real vector space and the $n \times n$ dimensional real matrix space, $|\cdot|$ is the absolute value, $\langle\cdot,\cdot\rangle$ means the inner product on $\mathbb{R}^n$, and $sign(\cdot)$ represents the sign function. For a given vector or matrix $M$, $M^T$ denotes its transpose. For a given function $f(x,t)$, $f_{xx}(x,t)$ and $f_{t}(x,t)$ are the second-order partial derivative of $f(x,t)$ with the spatial variable $x$ and the first-order partial derivative of $f(x,t)$ with the temporal variable $t$, respectively. Let $\mathcal{L}^{\infty}([0,L], \mathbb{R}^n)=\left\{f(x):[0,L]\rightarrow\mathbb{R}^n|\right.$  there is a subset $E_0\subset [0,L]$  with measure zero such that $f(x)$ is bounded on  $\left. [0,L] \backslash E_0\right\}$.

\section{Preliminaries}

\subsection{Graph theory}
\noindent

In MASs, $G=\{\mathfrak{V},\mathfrak{E},A\}$ represent a fixed communication topology graph, where $\mathfrak{V}=\{v_{1},v_{2},\ldots,v_{n}\}$ means the vertex set of the graph $G$ with $n$ nodes, $\mathfrak{E}\in \mathfrak{V}\times \mathfrak{V}$ is an edge set, and $A=(a_{ij})$ is the adjacency matrix of the graph $G$. If the node $v_{i}$ can receive the data from $v_{j}$, then $a_{ij}>0$, else $a_{ij}=0$. Thus, we can define $N_i=\{v_j\in \mathfrak{V}: (v_j,v_i)\in \mathfrak{E}\}$. The Laplace's matrix $\mathcal{L}$ of the graph $G$ is

\begin{eqnarray*}
\mathcal{L}(A)=\left(l_{ij}\right)=\left\{
\begin{array}{lll}
-a_{ij},& v_{j}\in N_i;\\
\sum_{v_j\in N_i}a_{ij},&j=i;\\
0,&other.
\end{array}
\right.
\end{eqnarray*}

\begin{definition}
  An undirected graph is said to be connected if there exists a sequence $S =\{a_{ii_1}, a_{i_1i_2}, \cdots, a_{i_mj}\}$ such that each element of $S$ is positive for any pair of $v_{i}$ and $v_{j}$.
\end{definition}

\begin{definition}
  A directed graph is said to be
\begin{itemize}
\item[(i)] strongly connected (s-con) if there exists a sequence $S =\{a_{ii_1}, a_{i_1i_2}, \cdots, a_{i_mj}\}$ such that each element of $S$ is positive for any pair of $v_{i}$ and $v_{j}$;
\item[(ii)] detail-balanced (d-bal) if there exists a vector $\omega = [\omega _1, \omega _2, \cdots, \omega _N]^{T}$ such that $\omega_i >0$ and $\omega_ia_{ij}=\omega_ja_{ji}$ for $i, j = 1, 2, \cdots, N$.
\end{itemize}
\end{definition}

\begin{definition}
  If $A^T=A$, then $R(x;A)=\frac{\langle x, Ax \rangle}{\langle x, x \rangle}$ is said to be the Rayleigh quotient of $A$, where $x\neq 0$.
\end{definition}

\subsection{Critical Lemmas}
\noindent

In this subsection, we recall some important lemmas which will be used in our work.

\begin{lemma}\label{l1}\cite{Bh00}
 Assume that the differentiable function $\mathbb{V}(t):\{0\}\bigcup\mathbb{R}^+ \rightarrow \{0\}\bigcup\mathbb{R}^+$ is positive definite and there are $k>0$ and $0<\alpha<1$ such that
 \begin{eqnarray*}
 \frac{\mathrm{d}\mathbb{V}(t)}{\mathrm{d}t}\leq -k\mathbb{V}(t)^{\alpha}.
 \end{eqnarray*}
 Then, $\mathbb{V}(t)$ converges to $0$ in finite time $t^{*}$. Moreover, the terminal time $t^*$ satisfies $ t^* \leq \frac{V(0)^{1-\alpha}}{K(1-\alpha)}$.

\end{lemma}

\begin{lemma}\label{l5} \cite{Po12}
 Assume that the differentiable function $\mathbb{V}(t):\{0\}\bigcup\mathbb{R}^+ \rightarrow \{0\}\bigcup\mathbb{R}^+$ is positive definite and there are  $k, p, q>0$,  $0< \alpha k <1$ and $\beta k >1$ such that
 \begin{eqnarray*}
 \frac{\mathrm{d}\mathbb{V}(t)}{\mathrm{d}t}\leq -\left(p\mathbb{V}(t)^{\alpha}+q\mathbb{V}(t)^{\beta}\right)^k.
 \end{eqnarray*}
 Then, $\mathbb{V}(t)$ converges to $0$ in fixed time $T$. Moreover, the terminal time $T$ satisfies $ T \leq T_{\max}:=\frac{1}{p^k(1-\alpha k)}+\frac{1}{q^k(\beta k-1)}$.

\end{lemma}

\begin{lemma}\label{l2}\cite{Wa10}
  For any $0<\alpha<1$, $\beta >1 $  and $x_{i}\geq 0$ with $i=1,2,\cdots,n$, one has
 \begin{eqnarray*}
  \left(\sum_{i=1}^{n}x_{i}\right)^{\alpha} \leq \sum_{i=1}^{n}x_{i}^{\alpha}\leq n^{1-\alpha}\left(\sum_{i=1}^{n}x_{i}\right)^{\alpha},
 \end{eqnarray*}
   and
 \begin{eqnarray*}
 \sum_{i=1}^{n}x_{i}^{\beta} \leq \left(\sum_{i=1}^{n}x_{i}\right)^{\beta} \leq  n^{\beta-1}\sum_{i=1}^{n}x_{i}^{\beta}.
 \end{eqnarray*}
\end{lemma}

\begin{lemma}\label{l3} \cite{Ol04, Shi20}
  If $G=\{\mathfrak{V},\mathfrak{E},A\}$ is undirected and connected, then $\mathcal{L}(A)$ has the following properties:
  \begin{itemize}
    \item [(i)] $\zeta^{T}\mathcal{L}(A)\zeta=\frac{1}{2}\Sigma_{i,j=1}^{n}a_{ij}(\zeta_j-\zeta_i)^2$ for any $\zeta=[\zeta_1, \zeta_2, \cdots, \zeta_n]^{T}\in \mathbb{R}^n$;
    \item [(ii)] For any $\zeta=[\zeta_1, \zeta_2, \cdots, \zeta_n]^{T}\in \mathbb{R}^n$, if $\langle \zeta, \textbf{1}\rangle=0$ with $\textbf{1}=[1, 1, \cdots, 1]^T\in \mathbb{R}^n$, then $\zeta^{T}\mathcal{L}(A)\zeta\geq\lambda_2(\mathcal{L}(A))\zeta^{T}\zeta$, where $\lambda_2(\mathcal{L}(A))= \underset{\langle x, \textbf{1}\rangle=0}\min R(x, \mathcal{L}(A))$ is the second smallest eigenvalue of $\mathcal{L}(A)$;
    \item [(iii)] If $\lambda_{\omega}(\mathcal{L}(A))=\underset{\langle x, \omega\rangle=0}\min R(x, \mathcal{L}(A))$, then $\lambda_2(\mathcal{L}(A))\geq \lambda_{\omega}(\mathcal{L}(A)) \geq0$. Specially, $\omega$ is perpendicular to the left eigenvector of the zero eigenvalue of $\mathcal{L}(A)$ if and only if $\lambda_{\omega}(\mathcal{L}(A))=0$.
  \end{itemize}
\end{lemma}

\begin{lemma}\label{l4}\cite{Wu14}
  Let $\varphi(\cdot,t) \in \mathcal{H}^{1}(0,L) \triangleq \{f(x):[0,L]\rightarrow \mathbb{R}|\int_0^Lf^{T}(x)f(x)\mathrm{d}x<\infty, \int_0^Lf_{x}^{T}(x)f_{x}(x)\mathrm{d}x<\infty\}$  with $\varphi(0)=0$ (or $\varphi(L)=0$). Then for any real symmetric positive definite matrix $\mathcal{P}$, one has
 \begin{eqnarray*}
 \int_0^L\varphi^{T}(s)\mathcal{P}\varphi(s)\mathrm{d}s\leq 4L^2\pi^{-2}\int_{0}^{L}\varphi_{s}^{T}(s)\mathcal{P}\varphi_{s}(s)\mathrm{d}s.
 \end{eqnarray*}
\end{lemma}

\begin{lemma}\label{l6}\cite{Ev81}
  There exists a unique solution $\psi(x, t)$ of the system
\begin{eqnarray*}
\left\{
\begin{array}{llll}
\psi_{t}(x,t)-\psi_{xx}(x,t)=f(x,t),\\
\psi(0,t)=\psi(L,t)=0,\\
\psi(x,0)=\psi^0(x),
\end{array}
\right.
\end{eqnarray*}
where $x \in [0, L]$, $\psi^0(x)$ is continuous on $[0, L] $, and $f(x,t)$ is continuous on $[0, L] \times [0, T]$.
\end{lemma}

\begin{lemma}\label{l7} \cite{H087}
Consider the nonlinear parabolic PDE with the following form:
\begin{eqnarray}\label{e2.0.1}
\left\{
\begin{array}{llll}
\psi_{t}(x,t)-k\psi_{xx}(x,t)=f(\psi(x,t)),\\
\psi(0,t)=\psi(L,t)=0,\\
\psi(x,0)=\psi^0(x),
\end{array}
\right.
\end{eqnarray}
where $\psi(x, t)=[\psi_1(x, t), \psi_2(x, t),\cdots, \psi_n(x, t)]^T$, $\psi^0(\cdot)\in \mathcal{L}^{\infty}([0,L], \mathbb{R}^n)$, and $f(\cdot)$ is continuously differentiable. Then there exists $T_M>0$ ensuring the local existence and uniqueness of the classical solution of \eqref{e2.0.1} on $[0, L] \times [0, T_M)$. Moreover, if $T_M < \infty$, then
$$
\lim_{t\uparrow T_M}\sup_{x \in [0, L]}\{|\psi_i(x,t)|\}=\infty
$$
for  some $1\leq i\leq n$.
\end{lemma}

\subsection{Problem Formulation}
\noindent

In this work, we will consider the FTC and FXC problems of parabolic partial differential MASs with external disturbance. Assume that MASs including $N$ agents have the following form:
\begin{eqnarray}\label{e3.2.1}
\left\{
\begin{array}{llll}
y_{i,t}(x,t)=ky_{i,xx}(x,t)+u_i(x,t)+d_i(x,t),\\
y_{i}(0,t)= y_{i}(L,t)=0,\\
y_i(x,0)=y_i^0(x),\quad i=1, 2, \cdots, N,
\end{array}
\right.
\end{eqnarray}
where $x \in [0,L]$ is the location, $t \in \{0\}\bigcup\mathbb{R}^+$ is the elapsed time, $y_i(x,t) \in \mathbb{R}$ is the state variable of the $i$th agent, $y_i^0(\cdot)$ is the initial condition of the $i$th agent such that $y_i^0(\cdot)$ is bounded and continuous with respect to $x$, $u_i(x,t)\in \mathbb{R}$ is the protocol of the $i$th agent, $d_i(x,t)\in \mathbb{R}$ is the unknown external disturbance, and $k>0$ is a constant.

Assume that $d_i(x,t)\in \mathbb{R}$ satisfies the following conditions:
\begin{itemize}
    \item [(i)] $d_i(x,t)$ is continuous with respect to $x$ and $t$;
    \item [(ii)] $|d_i(x,t)|\leq d_{\max}$ for any given $x \in [0,L]$, $t \in \left[0,+\infty\right)$ and $i=1, 2, \cdots, N$, where the upper bound $d_{\max}$ is a known constant.
  \end{itemize}

\begin{definition}
MASs \eqref{e3.2.1} achieves the FTC if there exists a finite time $t^*$ for any given $x \in [0,L]$ and $i,j=1,2,\cdots,N$ satisfying
 \begin{eqnarray*}
  \lim_{t\rightarrow t^*}|y_i(x,t)-y_j(x,t)| = 0.
 \end{eqnarray*}
Moreover, MASs \eqref{e3.2.1}  achieves the FXC if the finite time $t^*<T_{\max}$, where $T_{\max}$ is a positive constant and independent of the initial condition.
\end{definition}

\begin{remark}
It is worth mentioning that the control objectives of MASs \eqref{e3.2.1} are to design the corresponding controllers to achieve the FTC or the FXC.
\end{remark}

For the sake of simplicity, let's note $\xi_{ij}(x,t)=y_i(x,t)-y_j(x,t)$.
To achieve the FTC, we design the controller $u_i(x,t)$ via state values by setting
\begin{eqnarray}\label{e3.2.2}
u_i(x,t)&=&\underbrace{-\sum_{j=1}^{N}a_{ij}\left(\int_0^L\xi_{ij}(x,t)^{2}\mathrm{d}x\right)^{\alpha-1}\xi_{ij}(x,t)}_{(I)}\nonumber\\
&&\underbrace{-\bar{d}\sum_{j=1}^{N}a_{ij}sign(\xi_{ij}(x,t))}_{(II)},
\end{eqnarray}
where $0<\alpha<1 $, $a_{ij}\ge 0$ and $\bar{d}>0$ are given constants.

\begin{remark}
In the controller \eqref{e3.2.2}, Part (I) ensure the FTC and the role of Part (II) is to counteract the influence of external disturbance.
\end{remark}

To reach the FXC, we design the controller $u_i(x,t)$ via state values by setting
\begin{eqnarray}\label{e13.2.2}
u_i(x,t)\nonumber&=&\underbrace{-\sum_{j=1}^{N}a_{ij}\left(\int_0^L\xi_{ij}(x,t)^{2}\mathrm{d}x\right)^{\alpha-1}\xi_{ij}(x,t)}_{(I)}\nonumber\\
&&\underbrace{-\sum_{j=1}^{N}a_{ij}\left(\int_0^L\xi_{ij}(x,t)^{2}\mathrm{d}x\right)^{\beta-1}\xi_{ij}(x,t)}_{(II)}\nonumber\\
&&\underbrace{-\bar{d}\sum_{j=1}^{N}a_{ij}sign(\xi_{ij}(x,t))}_{(III)},
\end{eqnarray}
where $0<\alpha<1 $, $\beta >1$, $a_{ij}\ge 0$ and $\bar{d}>0$ are given constants.

\begin{remark}
Compared with the controller \eqref{e3.2.2}, Part (II) of controller \eqref{e13.2.2} is to make the finite terminal time independent of the initial condition.
\end{remark}

\begin{remark}
Noting that the FXC problem can be considered as an evolution of the FTC problem. In fact, the FTC and the FXC controllers ensure the systems to achieve consensus in finite time, but they have the different convergence rates and the different estimations for the finite terminal time.
\end{remark}

Obviously, if $d_i(x,t)\equiv0$, then MASs \eqref{e3.2.1} is reduced to the following form:
\begin{eqnarray}\label{e3.1.1}
\left\{
\begin{array}{llll}
y_{i,t}(x,t)=ky_{i,xx}(x,t)+u_i(x,t),\\
y_{i}(0,t)= y_{i}(L,t)=0,\\
y_i(x,0)=y_i^0(x),\quad i=1,2,\cdots, N.
\end{array}
\right.
\end{eqnarray}
In other words,  MASs \eqref{e3.1.1} is  a special case of MASs \eqref{e3.2.1}.

 To achieve the FTC, we design the controller $u_i(x,t)$ via state values by setting
\begin{eqnarray}\label{e3.1.2}
u_i(x,t)=-\sum_{j=1}^{N}a_{ij}\left(\int_0^L\xi_{ij}(x,t)^{2}\mathrm{d}x\right)^{\alpha-1}\xi_{ij}(x,t),
\end{eqnarray}
where $0<\alpha<1$ and $a_{ij}\ge 0$ are constants.

Similarly, to achieve the FXC, we design the controller $u_i(x,t)$ via state values by setting
\begin{eqnarray}\label{e13.1.2}
u_i(x,t)
&=&-\sum_{j=1}^{N}a_{ij}\left(\int_0^L\xi_{ij}(x,t)^{2}\mathrm{d}x\right)^{\alpha-1}\xi_{ij}(x,t)\nonumber\\
&&-\sum_{j=1}^{N}a_{ij}\left(\int_0^L\xi_{ij}(x,t)^{2}\mathrm{d}x\right)^{\beta-1}\xi_{ij}(x,t),
\end{eqnarray}
where $0<\alpha<1$, $\beta>1$ and $a_{ij}\ge 0$ are constants.

\begin{remark}
  Since MASs \eqref{e3.1.1} is a special case of MASs \eqref{e3.2.1} with $d_i(x,t)\equiv0$, we can simplify the controller \eqref{e3.2.2} as the controller \eqref{e3.1.2} to reach FTC via omitting the Part (II). Similarly, we can obtain the fixed-time controller \eqref{e13.1.2} of MASs \eqref{e3.1.1}.
\end{remark}

\subsection{Unique Solvability of MASs \eqref{e3.2.1}}
\begin{proposition}\label{t0}
Suppose that graph $G$ is undirected and connected. Then MASs \eqref{e3.2.1} with the controller \eqref{e3.2.2} (or \eqref{e13.2.2}) has a unique solution.
\end{proposition}

{\bf Proof}\hspace{0.2cm} For convenience, the variables $(x,t)$ are omitted in the following proof. We first prove that the following system \eqref{e0.1} has a unique solution
\begin{eqnarray}\label{e0.1}
\left\{
\begin{array}{llll}
y_{t}=y_{xx}+d,\\
y(0,t)= y(L,t)=0,\\
y(x,0)=y^0(x),
\end{array}
\right.
\end{eqnarray}
where $y=[y_1, \cdots, y_N]^T$, $d=[d_1, \cdots, d_N]^T$, and $y^0(x)=[y_{1}^{0}(x), \cdots, y_
{N}^{0}(x)]^T$. In fact, since $y^0(x)$ and $d$ are continuous, it follows from Lemma \ref{l6} that the system \eqref{e0.1} has a unique solution.

Next we show that the following system \eqref{e0.2} has a unique solution
\begin{eqnarray}\label{e0.2}
\left\{
\begin{array}{llll}
y_{t}-ky_{xx}=u(y),\\
y(0,t)= y(L,t)=0,\\
y(x,0)=y^0(x),
\end{array}
\right.
\end{eqnarray}
where $y=[y_1, \cdots, y_N]^T$ and  $y^0(x)=[y_{1}^{0}(x), \cdots, y_
{N}^{0}(x)]^T$. Indeed, since $y^0(x)$ is bounded and continuous with respect to $x$, we know that $y^0(\cdot) \in \mathcal{L}^{\infty}([0,L], \mathbb{R}^n)$. Noting that the controllers \eqref{e3.2.2} and \eqref{e13.2.2} include sign functions, we can set the signal switching times to be $t_m$ with $0=t_0< t_1< \cdots < t_m <\cdots$. Thus, it follows that $u(\cdot)$ is continuously differentiable with respect to $y$ on $[t_m, t_{m+1})$. From Lemma \ref{l7}, there exists a constant $T_{M,m}>0$ such that the system \eqref{e0.2} has a unique solution on $[0,L] \times [t_m, \min\{t_{m+1},T_{M,m}\})$.  Construct a Lyapunov generalized energy function
\begin{eqnarray*}
\mathbb{V}(t)=\int_{0}^{L}\sum_{i=1}^{N}(y_i)^2\mathrm{d}x.
\end{eqnarray*}
Then it is easy to see that the derivative of $\mathbb{V}(t)$ with $t$ can be given by
\begin{eqnarray*}
\frac{\mathrm{d}\mathbb{V}(t)}{\mathrm{d}t}
&=&2\int_{0}^{L}\sum_{i=1}^{N}y_iy_{i,t}\mathrm{d}x\nonumber\\
&=&2k\sum_{i=1}^{N}\int_{0}^{L}y_iy_{i,xx}\mathrm{d}x+2\sum_{i=1}^{N}\int_{0}^{L}y_iu_i\mathrm{d}x\nonumber\\
&=& -2\sum_{i=1}^{N}\int_{0}^{L}y_i\sum_{j=1}^{N}a_{ij}\left(\int_0^L\xi_{ij}^{2}\mathrm{d}x\right)^{\alpha-1}\xi_{ij}\mathrm{d}x\nonumber\\
&& -2\sum_{i=1}^{N}\int_{0}^{L}y_i\bar{d}\sum_{j=1}^{N}a_{ij}sign(\xi_{ij})\mathrm{d}x-2k\sum_{i=1}^{N}\int_{0}^{L}(y_{i,x})^2\mathrm{d}x \nonumber\\
&\leq& -\sum_{i,j=1}^{N}\left(a_{ij}^{\frac{1}{\alpha}}\int_0^L\xi_{ij}^{2}\mathrm{d}x\right)^{\alpha}-\int_{0}^{L}\sum_{i,j=1}^{N}\bar{d}a_{ij}\left(\xi_{ij}^2\right)^{\frac{1}{2}}\mathrm{d}x\nonumber\\
&&-\frac{1}{2}kL^{-2}\pi^2\sum_{i=1}^{N}\int_{0}^{L}(y_i)^2\mathrm{d}x\nonumber\\
&\leq& 0.
\end{eqnarray*}
Thus, one has $\mathbb{V}(t) \leq \mathbb{V}(0) < \infty$. According to Lemma \ref{l7}, we have $T_{M,m}= \infty$ for $m=0,1,2,\cdots$ and so the system \eqref{e0.2} has a unique solution on $[0,L] \times [t_m, t_{m+1})$ for $m=0,1,2,\cdots$, which means that the system \eqref{e0.2} has a unique solution in the whole domain.

Thus, by the unique solvability of systems \eqref{e0.1} and \eqref{e0.2}, we know that  MASs \eqref{e3.2.1} with the controller \eqref{e3.2.2} (or \eqref{e13.2.2}) has a unique solution. This completes the proof.

\begin{remark}
Since MASs \eqref{e3.1.1} is  a special case of MASs \eqref{e3.2.1}, it follows from Proposition \ref{t0} that MASs \eqref{e3.1.1} with \eqref{e3.1.2} (or \eqref{e13.1.2}) has a unique solution.
\end{remark}

\section{Main results}
\noindent

Based on the controllers \eqref{e3.2.2} and \eqref{e13.2.2}, the sufficient conditions for the FTC and FXC of MASs \eqref{e3.2.1} are obtained. For convenience, the variables $(x,t)$ are omitted in the following proof.

\begin{theorem}\label{t5}
Suppose that graph $G$ is undirected and connected. Then the controller \eqref{e3.2.2} can ensure MASs \eqref{e3.2.1} to achieve the FTC providing $\bar{d} \geq\frac{2d_{\max}N^{\frac{1}{2}}}{\sqrt{2}\lambda^{\frac{1}{2}}_2(\mathcal{L}(E))}$ with $E =(a_{ij}^{2})\in \mathbb{R}^{N\times N}$ and the finite terminal time $t^*$ satisfies
 \begin{eqnarray*}
t^*\leq\frac{\left[\int_{0}^{L}\frac{1}{2N}\sum_{i,j=1}^{N}(y_i^0(x)-y_j^0(x))^2\mathrm{d}x\right]^{1-\alpha}}{C_1(1-\alpha)},
\end{eqnarray*}
where $C_1=\left(2\lambda_2(\mathcal{L}(B))\right)^{\alpha}$ with $B =(a_{ij}^{\frac{1}{\alpha}})\in \mathbb{R}^{N\times N}$.
\end{theorem}

{\bf Proof}\hspace{0.2cm} Design an auxiliary variable $y^*$ as follows
\begin{eqnarray*}
y^*=\frac{1}{N}\sum_{i=1}^{N}y_i.
\end{eqnarray*}
Obviously, the dynamic equation of $y^*(x,t)$ has the following form
\begin{eqnarray*}
\left\{
\begin{array}{llll}
 y^*_t=ky^*_{xx}+u^*+d^*,\\
y^*(0,t)= y^*(L,t)=0,\\
y^*(x,0)=y^{*,0}(x),
\end{array}
\right.
\end{eqnarray*}
where $u^*=\frac{1}{N}\sum_{i=1}^{N}u_i$ is the average protocol, $d^*=\frac{1}{N}\sum_{i=1}^{N}d_i$ is the average disturbance and $y^{*,0}(x)=\frac{1}{N}\sum_{i=1}^{N}y^0_i(x)$ is the average initial condition. Denote $\delta = [\delta _1, \delta _2, \cdots, \delta _N]^T$ with $\delta_i=y_i-y^*$. Then $\langle \delta, \textbf{1} \rangle =0$.

Now we need to prove that  $y_i$ achieve the consensus in finite time, which is equivalent to that $y_i$ converges to $y^*$ in finite time for $i=1,2,\cdots,N$. To this end, we construct a Lyapunov generalized energy function
\begin{eqnarray*}
\mathbb{V}(t)=\int_{0}^{L}\sum_{i=1}^{N}(y_i-y^*)^2\mathrm{d}x=\int_{0}^{L}\sum_{i=1}^{N}\delta_i^2\mathrm{d}x,
\end{eqnarray*}
which is equivalent to the following form
\begin{eqnarray*}
\mathbb{V}(t)=\int_{0}^{L}\frac{1}{2N}\sum_{i,j=1}^{N}\xi_{ij}^2\mathrm{d}x.
\end{eqnarray*}
Clearly, the derivative of $\mathbb{V}(t)$ with $t$ can be given by
\begin{eqnarray*}
\frac{\mathrm{d}\mathbb{V}(t)}{\mathrm{d}t}
&=&2\int_{0}^{L}\sum_{i=1}^{N}\delta_i(y_{i,t}-y^*_t)\mathrm{d}x\nonumber\\
&=&\underbrace{2k\sum_{i=1}^{N}\int_{0}^{L}\delta_i(y_{i,xx}-y^*_{xx})\mathrm{d}x}_{(I)}\nonumber\\
&&+\underbrace{2\sum_{i=1}^{N}\int_{0}^{L}\delta_i(d_i-d^*)\mathrm{d}x}_{(II)}+\underbrace{2\sum_{i=1}^{N}\int_{0}^{L}\delta_i(u_i-u^*)\mathrm{d}x}_{(III)}.
\end{eqnarray*}
To estimate the derivative of $\mathbb{V}(t)$, we need to estimate Part (I), Part (II) and Part (III), respectively.  Using Lemma \ref{l4} to estimate Part (I), one has
\begin{eqnarray}\label{e3.2.3}
&&2k\sum_{i=1}^{N}\int_{0}^{L}\delta_i(y_{i,xx}-y^*_{xx})\mathrm{d}x\nonumber\\
&=&-2k\sum_{i=1}^{N}\int_{0}^{L}(y_{i,x}-y^*_{x})^2\mathrm{d}x \nonumber\\
&\leq& -\frac{1}{2}kL^{-2}\pi^2\sum_{i=1}^{N}\int_{0}^{L}\delta_i^2\mathrm{d}x.
\end{eqnarray}
To estimate Part (II),  it is easy to have
\begin{eqnarray}\label{e3.2.4}
&&2\sum_{i=1}^{N}\int_{0}^{L}\delta_i(d_i-d^*)\mathrm{d}x\nonumber\\
&=&2\sum_{i=1}^{N}\int_{0}^{L}\delta_id_i\mathrm{d}x
\leq 2d_{\max}\int_{0}^{L}\sum_{i=1}^{N}|\delta_i|\mathrm{d}x\nonumber\\
&\leq& 2d_{\max}N^{\frac{1}{2}}\int_{0}^{L}\left(\sum_{i=1}^{N}\delta_i^2\right)^{\frac{1}{2}}\mathrm{d}x.
\end{eqnarray}
Since $G$ is undirected, we obtain that $a_{ij}=a_{ji}$. Employing \eqref{e3.2.2} and Lemma \ref{l2} to estimate Part (III), we obtain
\begin{eqnarray}\label{e3.2.5}
&&2\sum_{i=1}^{N}\int_{0}^{L}\delta_i(u_i-u^*)\mathrm{d}x\nonumber\\
&=&-2\sum_{i=1}^{N}\int_{0}^{L}y_i\sum_{j=1}^{N}a_{ij}\left(\int_0^L\xi_{ij}^{2}\mathrm{d}x\right)^{\alpha-1}\xi_{ij}\mathrm{d}x-2\sum_{i=1}^{N}\int_{0}^{L}y_i\bar{d}\sum_{j=1}^{N}a_{ij}sign(\xi_{ij})\mathrm{d}x\nonumber\\
&=&-\int_{0}^{L}\sum_{i,j=1}^{N}a_{ij}\xi_{ij}^2dx\left(\int_0^L\xi_{ij}^{2}\mathrm{d}x\right)^{\alpha-1}-\int_{0}^{L}\sum_{i,j=1}^{N}\bar{d}a_{ij}\xi_{ij}sign(\xi_{ij})\mathrm{d}x\nonumber\\
&=&-\sum_{i,j=1}^{N}a_{ij}\left(\int_0^L\xi_{ij}^{2}\mathrm{d}x\right)^{\alpha}-\int_{0}^{L}\sum_{i,j=1}^{N}\bar{d}a_{ij}|\xi_{ij}|\mathrm{d}x\nonumber\\
&=&-\sum_{i,j=1}^{N}\left(a_{ij}^{\frac{1}{\alpha}}\int_0^L\xi_{ij}^{2}\mathrm{d}x\right)^{\alpha}-\int_{0}^{L}\sum_{i,j=1}^{N}\bar{d}a_{ij}\left(\xi_{ij}^2\right)^{\frac{1}{2}}\mathrm{d}x\nonumber\\
&\leq&-\left(\int_0^L\sum_{i,j=1}^{N}a_{ij}^{\frac{1}{\alpha}}\xi_{ij}^{2}\mathrm{d}x\right)^{\alpha}-\bar{d}\int_{0}^{L}\left(\sum_{i,j=1}^{N}a^2_{ij}\xi_{ij}^2\right)^{\frac{1}{2}}\mathrm{d}x\nonumber\\
&\leq&-\left(\int_0^L\frac{\sum_{i,j=1}^{N}a_{ij}^{\frac{1}{\alpha}}\xi_{ij}^{2}}{\sum_{i=1}^{N}\delta_i^{2}}\sum_{i=1}^{N}\delta_i^{2}\mathrm{d}x\right)^{\alpha}-\bar{d}\int_{0}^{L}\left(\frac{\sum_{i,j=1}^{N}a^2_{ij}\xi_{ij}^2}{\sum_{i=1}^{N}\delta_i^{2}}\sum_{i=1}^{N}\delta_i^{2}\right)^{\frac{1}{2}}\mathrm{d}x.
\end{eqnarray}
Define $B =(a_{ij}^{\frac{1}{\alpha}})\in \mathbb{R}^{N\times N}$ and $E =(a_{ij}^{2})\in \mathbb{R}^{N\times N}$. According to Lemma \ref{l3}, we have
\begin{eqnarray}\label{e3.2.6}
\frac{\sum_{i,j=1}^{N}a_{ij}^{\frac{1}{\alpha}}\xi_{ij}^{2}}{\sum_{i=1}^{N}\delta_i^{2}}
&=&\frac{\sum_{i,j=1}^{N}a_{ij}^{\frac{1}{\alpha}}(\delta_j-\delta_i)^{2}}{\sum_{i=1}^{N}\delta_i^{2}}\nonumber\\
&=&\frac{2\langle \delta, \mathcal{L}(B)\delta \rangle}{\langle \delta, \delta \rangle} \geq 2 \lambda_2(\mathcal{L}(B))
\end{eqnarray}
and
\begin{eqnarray}\label{e3.2.7}
\frac{\sum_{i,j=1}^{N}a_{ij}^{2}\xi_{ij}^{2}}{\sum_{i=1}^{N}\delta_i^{2}}
&=&\frac{\sum_{i,j=1}^{N}a_{ij}^{2}(\delta_j-\delta_i)^{2}}{\sum_{i=1}^{N}\delta_i^{2}} \nonumber\\
&=&\frac{2\langle \delta, \mathcal{L}(E)\delta \rangle}{\langle \delta, \delta \rangle} \geq 2 \lambda_2(\mathcal{L}(E)).
\end{eqnarray}
Substituting \eqref{e3.2.6} and \eqref{e3.2.7} into \eqref{e3.2.5}, one has
\begin{eqnarray}\label{e3.2.8}
&&2\sum_{i=1}^{N}\int_{0}^{L}\delta_i(u_i-u^*)\mathrm{d}x\nonumber\\
&\leq&-\left(\int_0^L 2\lambda_2(\mathcal{L}(B))\sum_{i=1}^{N}\delta_i^{2}\mathrm{d}x\right)^{\alpha}-\bar{d}\int_{0}^{L}\left(2\lambda_2(\mathcal{L}(E))\sum_{i=1}^{N}\delta_i^{2}\right)^{\frac{1}{2}}\mathrm{d}x\nonumber\\
&\leq&-(2\lambda_2(\mathcal{L}(B)))^{\alpha}\left(\int_0^L \sum_{i=1}^{N}\delta_i^{2}\mathrm{d}x\right)^{\alpha}-\sqrt{2}\bar{d}\lambda^{\frac{1}{2}}_2(\mathcal{L}(E))\int_{0}^{L}\left(\sum_{i=1}^{N}\delta_i^{2}\right)^{\frac{1}{2}}\mathrm{d}x\nonumber\\
&\leq&-C_1(V(t))^{\alpha}-\sqrt{2}\bar{d}\lambda^{\frac{1}{2}}_2(\mathcal{L}(E))\int_{0}^{L}\left(\sum_{i=1}^{N}\delta_i^{2}\right)^{\frac{1}{2}}\mathrm{d}x,
\end{eqnarray}
where $C_1=\left(2\lambda_2(\mathcal{L}(B))\right)^{\alpha}$.  Since $\bar{d} \geq\frac{2d_{\max}N^{\frac{1}{2}}}{\sqrt{2}\lambda^{\frac{1}{2}}_2(\mathcal{L}(E))}$, it follows from \eqref{e3.2.3}, \eqref{e3.2.4} and \eqref{e3.2.8} that
\begin{eqnarray*}
\frac{d\mathbb{V}(t)}{dt}
&\leq&-C_{1}(\mathbb{V}(t))^{\alpha}-\frac{1}{2}kL^{-2}\pi^2\sum_{i=1}^{N}\int_{0}^{L}\delta_i^2\mathrm{d}x\\
&&+\left(2d_{\max}N^{\frac{1}{2}}-\sqrt{2}\bar{d}\lambda^{\frac{1}{2}}_2(\mathcal{L}(E))\right)\int_{0}^{L}\left(\sum_{i=1}^{N}\delta_i^{2}\right)^{\frac{1}{2}}\mathrm{d}x\\
&\leq&-C_1(\mathbb{V}(t))^{\alpha}.
\end{eqnarray*}
Since $\alpha \in (0,1)$ and the undirected graph $G$ is connected, we know that $C_1>0$. Thus, Lemma \ref{l1} shows that MASs \eqref{e3.2.1} achieves the FTC in $t^*$ with
\begin{eqnarray*}
t^*\leq\frac{\mathbb{V}(0)^{1-\alpha}}{C_1(1-\alpha)}=\frac{\left[\int_{0}^{L}\frac{1}{2N}\sum_{i,j=1}^{N}(y_i^0(x)-y_j^0(x))^2dx\right]^{1-\alpha}}{C_1(1-\alpha)}.
\end{eqnarray*}
The proof is completed.

\begin{theorem}\label{t6}
Suppose that graph $G$ is undirected and connected. Then the controller \eqref{e13.2.2} can ensure MASs \eqref{e3.2.1} to achieve the FXC providing $\bar{d} \geq\frac{2d_{\max}N^{\frac{1}{2}}}{\sqrt{2}\lambda^{\frac{1}{2}}_2(\mathcal{L}(E))}$ with $E =(a_{ij}^{2})\in \mathbb{R}^{N\times N}$ and the finite terminal time $t^*$ satisfies
\begin{eqnarray*}
t^*\leq T_{\max}:=\frac{1}{C_1(1-\alpha )}+\frac{1}{(N^{1-\beta})^2C_2(\beta -1)},
\end{eqnarray*}
where
\begin{eqnarray*}
C_1=\left(2\lambda_2(\mathcal{L}(B))\right)^{\alpha}, \quad C_2=\left(2\lambda_2(\mathcal{L}(P))\right)^{\beta}
\end{eqnarray*}
with $B =(a_{ij}^{\frac{1}{\alpha}})\in \mathbb{R}^{N\times N}$ and  $P =(a_{ij}^{\frac{1}{\beta}})\in \mathbb{R}^{N\times N}$.
\end{theorem}

{\bf Proof}\hspace{0.2cm} Define a Lyapunov generalized energy function for MASs \eqref{e3.2.1} by setting
\begin{eqnarray*}
\mathbb{V}(t)=\int_{0}^{L}\sum_{i=1}^{N}\delta_i^2\mathrm{d}x.
\end{eqnarray*}
Similar to the discussion in Theorem \ref{t5}, the derivative of $\mathbb{V}(t)$ with $t$ can be given by
\begin{eqnarray*}
\frac{\mathrm{d}\mathbb{V}(t)}{\mathrm{d}t}
&=&2\int_{0}^{L}\sum_{i=1}^{N}\delta_i(y_{i,t}-y^*_t)\mathrm{d}x\nonumber\\
&\leq&-C_1(\mathbb{V}(t))^{\alpha}-\sum_{i,j=1}^{N}\left(a_{ij}^{\frac{1}{\beta}}\int_0^L\xi_{ij}^{2}\mathrm{d}x\right)^{\beta}\nonumber\\
&\leq&-C_1(\mathbb{V}(t))^{\alpha}-(N^{1-\beta})^2\left(\int_0^L\sum_{i,j=1}^{N}a_{ij}^{\frac{1}{\beta}}\xi_{ij}^{2}\mathrm{d}x\right)^{\beta}\nonumber\\
&\leq&-C_1(\mathbb{V}(t))^{\alpha}-(N^{1-\beta})^2C_2(\mathbb{V}(t))^{\beta},
\end{eqnarray*}
where
\begin{eqnarray*}
C_1=\left(2\lambda_2(\mathcal{L}(B))\right)^{\alpha}, \quad C_2=\left(2\lambda_2(\mathcal{L}(P))\right)^{\beta}
\end{eqnarray*}
with $B =(a_{ij}^{\frac{1}{\alpha}})\in \mathbb{R}^{N\times N}$ and $P =(a_{ij}^{\frac{1}{\beta}})\in \mathbb{R}^{N\times N}$. Thus, Lemma \ref{l5} shows that MASs \eqref{e3.2.1} achieves the FXC in $t^*$ with
\begin{eqnarray*}
t^*\leq T_{\max}:=\frac{1}{C_1(1-\alpha )}+\frac{1}{(N^{1-\beta})^2C_2(\beta -1)}.
\end{eqnarray*}
This completes the proof.

Since MASs \eqref{e3.1.1} is a special case of MASs \eqref{e3.2.1}, we have the following corollaries which show the FTC and FXC of MASs \eqref{e3.1.1} via the controllers \eqref{e3.1.2} and \eqref{e13.1.2} under the undirected and connected graph.

\begin{corollary}\label{t1}Suppose that graph $G$ is undirected and connected. Then the controller \eqref{e3.1.2} can ensure MASs \eqref{e3.1.1} to achieve the FTC and the finite terminal time $t^*$ satisfies
\begin{eqnarray*}
t^*\leq \frac{\left[\int_{0}^{L}\frac{1}{2N}\sum_{i,j=1}^{N}(y_i^0(x)-y_j^0(x))^2\mathrm{d}x\right]^{1-\alpha}}{C_1(1-\alpha)},
\end{eqnarray*}
where $C_1=\left(2\lambda_2(\mathcal{L}(B))\right)^{\alpha}$ with $B =(a_{ij}^{\frac{1}{\alpha}})\in \mathbb{R}^{N\times N}$.
\end{corollary}

\begin{corollary}\label{t2}
Suppose that graph $G$ is undirected and connected. Then the controller \eqref{e13.1.2} can ensure MASs \eqref{e3.1.1} to achieve the FXC and the finite terminal time $t^*$ satisfies
\begin{eqnarray*}
t^*\leq T_{\max}:=\frac{1}{C_1(1-\alpha )}+\frac{1}{(N^{1-\beta})^2C_2(\beta -1)},
\end{eqnarray*}
where
\begin{eqnarray*}
C_1=\left(2\lambda_2(\mathcal{L}(B))\right)^{\alpha}, \quad C_2=\left(2\lambda_2(\mathcal{L}(P))\right)^{\beta}
\end{eqnarray*}
with $B =(a_{ij}^{\frac{1}{\alpha}})\in \mathbb{R}^{N\times N}$ and $P =(a_{ij}^{\frac{1}{\beta}})\in \mathbb{R}^{N\times N}$.
\end{corollary}

Next we focus on the FTC and FXC problem of MASs \eqref{e3.1.1} via the controllers \eqref{e3.1.2} and \eqref{e13.1.2} under the directed, s-con and d-bal graph.

\begin{theorem}\label{t3}
Suppose that graph $G$ is directed, s-con and d-bal. Then the controller \eqref{e3.1.2} can ensure MASs \eqref{e3.1.1} to achieve the FTC and the finite terminal time $t^*$ satisfies
\begin{eqnarray*}
t^*\leq\frac{\left[\int_{0}^{L}\sum_{i=1}^{N}\omega_i(y_i^0(x)-\hat{y}^0(x))^2\mathrm{d}x\right]^{1-\alpha}}{C_3(1-\alpha)},
\end{eqnarray*}
where
\begin{eqnarray*}
\hat{y}^{0}(x)=\frac{1}{\sum_{i=1}^{N}\omega_i}\sum_{i=1}^{N}\omega_iy_i^0(x)
\end{eqnarray*}
and
\begin{eqnarray*}
C_3=\left(\frac{2}{\underset{i=1,\cdots,N}\max\omega_i} \lambda_{\omega}(\mathcal{L}(D))\right)^{\alpha}
\end{eqnarray*}
with $D =((\omega_ia_{ij})^{\frac{1}{\alpha}})\in \mathbb{R}^{N\times N}$.
\end{theorem}

{\bf Proof}\hspace{0.2cm}   Since graph $G$ is d-bal, there exists a vector $\omega = [\omega _1, \omega _2, \cdots, \omega _N]^{T}$ with $\omega_i >0$ satisfying $\omega_ia_{ij}=\omega_ja_{ji}(i, j = 1, \cdots, N) $.  Design an auxiliary variable $\hat{y}$ as follows
\begin{eqnarray*}
\hat{y}=\frac{1}{\sum_{i=1}^{N}\omega_i}\sum_{i=1}^{N}\omega_iy_i.
\end{eqnarray*}
Obviously, the dynamic equation of $\hat{y}$ satisfies
\begin{eqnarray*}
\left\{
\begin{array}{llll}
 \hat{y}_t=k\hat{y}_{xx}+\hat{u},\\
\hat{y}(0,t)= \hat{y}(L,t)=0,\\
\hat{y}(x,0)=\hat{y}^{0}(x),
\end{array}
\right.
\end{eqnarray*}
where $\hat{u}=\frac{1}{\sum_{i=1}^{N}\omega_i}\sum_{i=1}^{N}\omega_iu_i$ is the weighting average protocol, and $\hat{y}^{0}(x)=\frac{1}{\sum_{i=1}^{N}\omega_i}\sum_{i=1}^{N}\omega_iy_i^{0}(x)$ is the weighting average initial condition. Let $\hat{\delta} = [\hat{\delta} _1, \hat{\delta} _2, \cdots, \hat{\delta} _N]^T$ with $\hat{\delta}_i=y_i-\hat{y}$. Then $\langle \hat{\delta}, \omega\rangle =0$.

Now we need to prove that  $y_i$ achieves the consensus in finite time, which is equivalent to that $y_i$ converges to $\hat{y}$ in finite time for $i=1,2,\cdots,N$. To this end, we construct a Lyapunov generalized energy function
\begin{eqnarray*}
\mathbb{V}(t)
=\int_{0}^{L}\sum_{i=1}^{N}\omega_i(y_i-\hat{y})^2\mathrm{d}x=\int_{0}^{L}\sum_{i=1}^{N}\omega_i\hat{\delta}_i^2\mathrm{d}x.
\end{eqnarray*}
It is easy to check that
\begin{eqnarray*}
\frac{\mathrm{d}\mathbb{V}(t)}{\mathrm{d}t}&=&2\int_{0}^{L}\sum_{i=1}^{N}\omega_i\hat{\delta}_i(y_{i,t}-\hat{y}_t)\mathrm{d}x\nonumber\\
&=&\underbrace{2k\sum_{i=1}^{N}\int_{0}^{L}\omega_i\hat{\delta}_i(y_{i,xx}-\hat{y}_{xx})\mathrm{d}x}_{(I)}+\underbrace{2\sum_{i=1}^{N}\int_{0}^{L}\omega_i\hat{\delta}_i(u_i-\hat{u})\mathrm{d}x}_{(II)}.\nonumber\\
\end{eqnarray*}
To estimate the derivative of $\mathbb{V}(t)$, we need to estimate Part (I) and Part (II), respectively. Using Lemma \ref{l4} to estimate Part (I), one has
\begin{eqnarray}\label{e3.1.7}
&&2k\sum_{i=1}^{N}\omega_i\int_{0}^{L}\hat{\delta}_i(y_{i,xx}-\hat{y}_{xx})\mathrm{d}x\nonumber\\
&=&-2k\sum_{i=1}^{N}\omega_i\int_{0}^{L}(y_{i,x}-\hat{y}_{x})^2\mathrm{d}x\nonumber\\
&\leq&-\frac{1}{2}kL^{-2}\pi^2\sum_{i=1}^{N}\omega_i\int_{0}^{L}\hat{\delta}_i^2\mathrm{d}x.
\end{eqnarray}
Since $\omega_ia_{ij}=\omega_ja_{ji}$ $(i, j = 1, \cdots, N)$, using \eqref{e3.1.2} and Lemma \ref{l2} to estimate Part (II), we have
\begin{eqnarray}\label{e3.1.8}
&&2\sum_{i=1}^{N}\omega_i\int_{0}^{L}\hat{\delta}_i(u_i-\hat{u})\mathrm{d}x\nonumber\\
&=&-2\sum_{i=1}^{N}\omega_i\int_{0}^{L}\hat{\delta}_i\sum_{j=1}^{N}a_{ij}\left(\int_0^L\xi_{ij}^{2}\mathrm{d}x\right)^{\alpha-1}\xi_{ij}\mathrm{d}x\nonumber\\
&=&-2\sum_{i=1}^{N}\int_{0}^{L}\omega_i y_i\sum_{j=1}^{N}a_{ij}\left(\int_0^L\xi_{ij}^{2}\mathrm{d}x\right)^{\alpha-1}\xi_{ij}\mathrm{d}x\nonumber\\
&=&-2\int_{0}^{L}\sum_{i,j=1}^{N}y_i\omega_ia_{ij}\left(\int_0^L\xi_{ij}^{2}\mathrm{d}x\right)^{\alpha-1}\xi_{ij}\mathrm{d}x\nonumber\\
&=&-\int_{0}^{L}\sum_{i,j=1}^{N}\omega_ia_{ij}\xi_{ij}^2dx\left(\int_0^L\xi_{ij}^{2}\mathrm{d}x\right)^{\alpha-1}\nonumber\\
&=&-\sum_{i,j=1}^{N}\omega_ia_{ij}\left(\int_0^L\xi_{ij}^{2}\mathrm{d}x\right)^{\alpha} \nonumber\\
&=&-\sum_{i,j=1}^{N}\left((\omega_ia_{ij})^{\frac{1}{\alpha}}\int_0^L\xi_{ij}^{2}\mathrm{d}x\right)^{\alpha}\nonumber\\
&\leq&-\left(\int_0^L\sum_{i,j=1}^{N}(\omega_ia_{ij})^{\frac{1}{\alpha}}\xi_{ij}^{2}\mathrm{d}x\right)^{\alpha}\nonumber\\
&\leq&-\left(\int_0^L\frac{\sum_{i,j=1}^{N}(\omega_ia_{ij})^{\frac{1}{\alpha}}\xi_{ij}^{2}}{\sum_{i=1}^{N}\omega_i\hat{\delta}_i^{2}}\sum_{i=1}^{N}\omega_i\hat{\delta}_i^{2}\mathrm{d}x\right)^{\alpha}.
\end{eqnarray}
Define $D =((\omega_ia_{ij})^{\frac{1}{\alpha}})\in \mathbb{R}^{N\times N}$. According to Lemma \ref{l3}, one has
\begin{eqnarray}\label{e3.1.9}
\frac{\sum_{i,j=1}^{N}(\omega_ia_{ij})^{\frac{1}{\alpha}}\xi_{ij}^{2}}{\sum_{i=1}^{N}\omega_i\hat{\delta}_i^{2}}
&=&\frac{\sum_{i,j=1}^{N}(\omega_ia_{ij})^{\frac{1}{\alpha}}(\hat{\delta}_j-\hat{\delta}_i)^{2}}{\sum_{i=1}^{N}\hat{\delta}_i^{2}}\nonumber\\
&\geq& \frac{2\langle \hat{\delta}, \mathcal{L}(D)\hat{\delta} \rangle}{\underset{i=1,\cdots,N}\max\omega_i\langle \hat{\delta}, \hat{\delta} \rangle} \geq \frac{2}{\underset{i=1,\cdots,N}\max\omega_i} \lambda_{\omega}(\mathcal{L}(D)).
\end{eqnarray}
It is easy to see that $\langle\omega, \textbf{1}\rangle \neq 0$ and so $\lambda_{\omega}(\mathcal{L}(D))>0$ from Lemma \ref{l3}.  Substituting \eqref{e3.1.9} into \eqref{e3.1.8}, we obtain
\begin{eqnarray}\label{e3.1.10}
&&2\sum_{i=1}^{N}\omega_i\int_{0}^{L}\hat{\delta}_i(u_i-\hat{u})\mathrm{d}x\nonumber\\
&\leq&-\left(\int_0^L \frac{2}{\underset{i=1,\cdots,N}\max\omega_i} \lambda_{\omega}(\mathcal{L}(D))\sum_{i=1}^{N}\hat{\delta}_i^{2}\mathrm{d}x\right)^{\alpha}\nonumber\\
&\leq&-(\frac{2}{\underset{i=1,\cdots,N}\max\omega_i} \lambda_{\omega}(\mathcal{L}(D)))^{\alpha}\left(\int_0^L \sum_{i=1}^{N}\omega_i\hat{\delta}_i^{2}\mathrm{d}x\right)^{\alpha}\nonumber\\
&\leq&-C_3(\mathbb{V}(t))^{\alpha},
\end{eqnarray}
where
\begin{eqnarray*}
C_3=\left(\frac{2}{\underset{i=1,\cdots,N}\max\omega_i} \lambda_{\omega}(\mathcal{L}(D))\right)^{\alpha}>0.
\end{eqnarray*}
 Now, it follows from \eqref{e3.1.7} and \eqref{e3.1.10} that
\begin{eqnarray*}
\frac{\mathrm{d}\mathbb{V}(t)}{\mathrm{d}t}&\leq&-C_{3}(\mathbb{V}(t))^{\alpha}-\frac{1}{2}kL^{-2}\pi^2\sum_{i=1}^{N}\omega_i\int_{0}^{L}\hat{\delta}_i^2dx\\
&\leq&-C_3(\mathbb{V}(t))^{\alpha}.
\end{eqnarray*}
Thus, Lemma \ref{l1} shows that MASs \eqref{e3.1.1} achieves the FTC in $t^*$ with
\begin{eqnarray*}
t^*\leq\frac{\mathbb{V}(0)^{1-\alpha}}{C_3(1-\alpha)}=\frac{\left[\int_{0}^{L}\sum_{i=1}^{N}\omega_i(y_i^0(x)-\hat{y}^0(x))^2\mathrm{d}x\right]^{1-\alpha}}{C_3(1-\alpha)}.
\end{eqnarray*}
The proof is completed.

\begin{theorem}\label{t4}
Suppose that graph $G$ is directed, s-con and d-bal. Then the controller \eqref{e13.1.2} can ensure MASs \eqref{e3.1.1} to achieve the FXC and the finite terminal time $t^*$ satisfies
\begin{eqnarray*}
t^*\leq T_{\max}:=\frac{1}{C_3(1-\alpha )}+\frac{1}{(N^{1-\beta})^2C_4(\beta -1)},
\end{eqnarray*}
where
\begin{eqnarray*}
C_3=\left(\frac{2}{\underset{i=1,\cdots,N}\max\omega_i} \lambda_{\omega}(\mathcal{L}(D))\right)^{\alpha}
\end{eqnarray*}
and
\begin{eqnarray*}
C_4=\left(\frac{2}{\underset{i=1,\cdots,N}\max\omega_i} \lambda_{\omega}(\mathcal{L}(Q))\right)^{\beta}
\end{eqnarray*}
with $D =((\omega_ia_{ij})^{\frac{1}{\alpha}})\in \mathbb{R}^{N\times N}$ and $Q =((\omega_ia_{ij})^{\frac{1}{\beta}})\in \mathbb{R}^{N\times N}$.
\end{theorem}

{\bf Proof}\hspace{0.2cm} Construct a Lyapunov generalized energy function by setting
\begin{eqnarray*}
\mathbb{V}(t)=\int_{0}^{L}\sum_{i=1}^{N}\omega_i\hat{\delta}_i^2\mathrm{d}x.
\end{eqnarray*}
Similar to the discussion in Theorem \ref{t3}, the derivative of $\mathbb{V}(t)$ with $t$ can be obtained as
\begin{eqnarray*}
\frac{\mathrm{d}\mathbb{V}(t)}{\mathrm{d}t}
&=&2\int_{0}^{L}\sum_{i=1}^{N}\omega_i\hat{\delta}_i(y_{i,t}-\hat{y}_t)\mathrm{d}x\nonumber\\
&\leq&-C_3(\mathbb{V}(t))^{\alpha}-\sum_{i,j=1}^{N}\left((\omega_ia_{ij})^{\frac{1}{\beta}}\int_0^L\xi_{ij}^{2}\mathrm{d}x\right)^{\beta}\nonumber\\
&\leq&-C_3(\mathbb{V}(t))^{\alpha}-(N^{1-\beta})^2\left(\int_0^L\sum_{i,j=1}^{N}(\omega_ia_{ij})^{\frac{1}{\beta}}\xi_{ij}^{2}\mathrm{d}x\right)^{\beta}\nonumber\\
&\leq&-C_3(\mathbb{V}(t))^{\alpha}-(N^{1-\beta})^2C_4(\mathbb{V}(t))^{\beta},
\end{eqnarray*}
where
\begin{eqnarray*}
C_3=\left(\frac{2}{\underset{i=1,\cdots,N}\max\omega_i} \lambda_{\omega}(\mathcal{L}(D))\right)^{\alpha}
\end{eqnarray*}
and
\begin{eqnarray*}
C_4=\left(\frac{2}{\underset{i=1,\cdots,N}\max\omega_i} \lambda_{\omega}(\mathcal{L}(Q))\right)^{\beta}
\end{eqnarray*}
with $D =((\omega_ia_{ij})^{\frac{1}{\alpha}})\in \mathbb{R}^{N\times N}$ and $Q =((\omega_ia_{ij})^{\frac{1}{\beta}})\in \mathbb{R}^{N\times N}$. Therefore, Lemma \ref{l5} shows that MASs \eqref{e3.1.1} achieves the FXC in $t^*$ with
\begin{eqnarray*}
t^*\leq T_{\max}:=\frac{1}{C_3(1-\alpha )}+\frac{1}{(N^{1-\beta})^2C_4(\beta -1)}.
\end{eqnarray*}
This completes the proof.

\section{Simulation results}
\noindent

We provide two numerical simulations to verify the feasibility of the controllers designed above in this section. It should be noted that we solve parabolic MASs via the conventional finite difference method.

\begin{figure}[H]
   \centering
   \includegraphics[width=5cm]{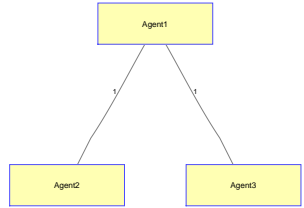}
   \caption{An undirected and connected graph with 3 agents.}\label{Fig11}
 \end{figure}

\textbf{Example 1} \quad Consider MASs with $3$ agents, the graph $G$ is undirected and connected given by Fig.\,\ref{Fig11}, where the number beside the edge means the adjacency element. Clearly, we can get the adjacency matrix
$A=\left(\begin{smallmatrix}
           0 & 1 & 1 \\
           1 & 0 & 0 \\
           1 & 0 & 0
         \end{smallmatrix}\right)$.

When $\alpha = 0.5 $, $\beta= 1.1$, $k= 10^{-3}$, $\bar{d}=8$ and $x \in [0,2]$, MASs can be obtained as
\begin{eqnarray}\label{e4.9}
\left\{
\begin{array}{llll}
y_{i,t}=10^{-3}\cdot y_{i,xx}+u_i+d_i,\\
y_{i}(0,t)= y_{i}(2,t)=0,\quad i=1, 2, 3,\\
\end{array}
\right.
\end{eqnarray}
the FTC controller can be formulated as
\begin{eqnarray}\label{e4.11}
u_i=-\sum_{j=1}^{3}a_{ij}\left(\int_0^2\xi_{ij}^{2}\mathrm{d}x\right)^{-0.5}\xi_{ij}-8\sum_{j=1}^{3}a_{ij}sign(\xi_{ij}),
\end{eqnarray}
and the FXC controller can be obtained as follows
\begin{eqnarray}\label{e4.12}
u_i&=&-\sum_{j=1}^{3}a_{ij}\left(\int_0^2\xi_{ij}^{2}\mathrm{d}x\right)^{-0.5}\xi_{ij}\nonumber\\
&&-\sum_{j=1}^{3}a_{ij}\left(\int_0^2\xi_{ij}^{2}\mathrm{d}x\right)^{0.1}\xi_{ij}-8\sum_{j=1}^{3}a_{ij}sign(\xi_{ij}).
\end{eqnarray}
The initial conditions are setting as
\begin{eqnarray}\label{e4.13}
\left\{
\begin{array}{lll}
y_1(x, 0)=y^0_1(x)=3\sin(\pi x),\\
y_2(x, 0)=y^0_2(x)=-2\cos(\pi x)+2,\\
y_3(x, 0)=y^0_3(x)=2\cos(\pi x)-2,
\end{array}
\right.
\end{eqnarray}
and the external disturbances are as follows
\begin{eqnarray}\label{e4.10}
\left\{
\begin{array}{lll}
d_1(x, t)=\sin(\pi  t)+x,\\
d_2(x, t)=\cos(\pi  t)+x,\\
d_3(x, t)=\sin(\pi  t  x).
\end{array}
\right.
\end{eqnarray}
Note that the distributed consensus controller considered in \cite{Fu18} was given by
\begin{eqnarray}\label{e4.c1}
u_i=-\sum_{j=1}^{3}a_{ij}\xi_{ij}
\end{eqnarray}
and the consensus boundary controller studied in \cite{Pi16} was designed as follows
\begin{eqnarray}\label{e4.c2}
y_{i,x}(2,t)=u_i(t)=-\sum_{j=1}^{3}a_{ij}\xi_{ij}(2,t).
\end{eqnarray}

\begin{figure}[H]
   \centering
   \includegraphics[width=4cm]{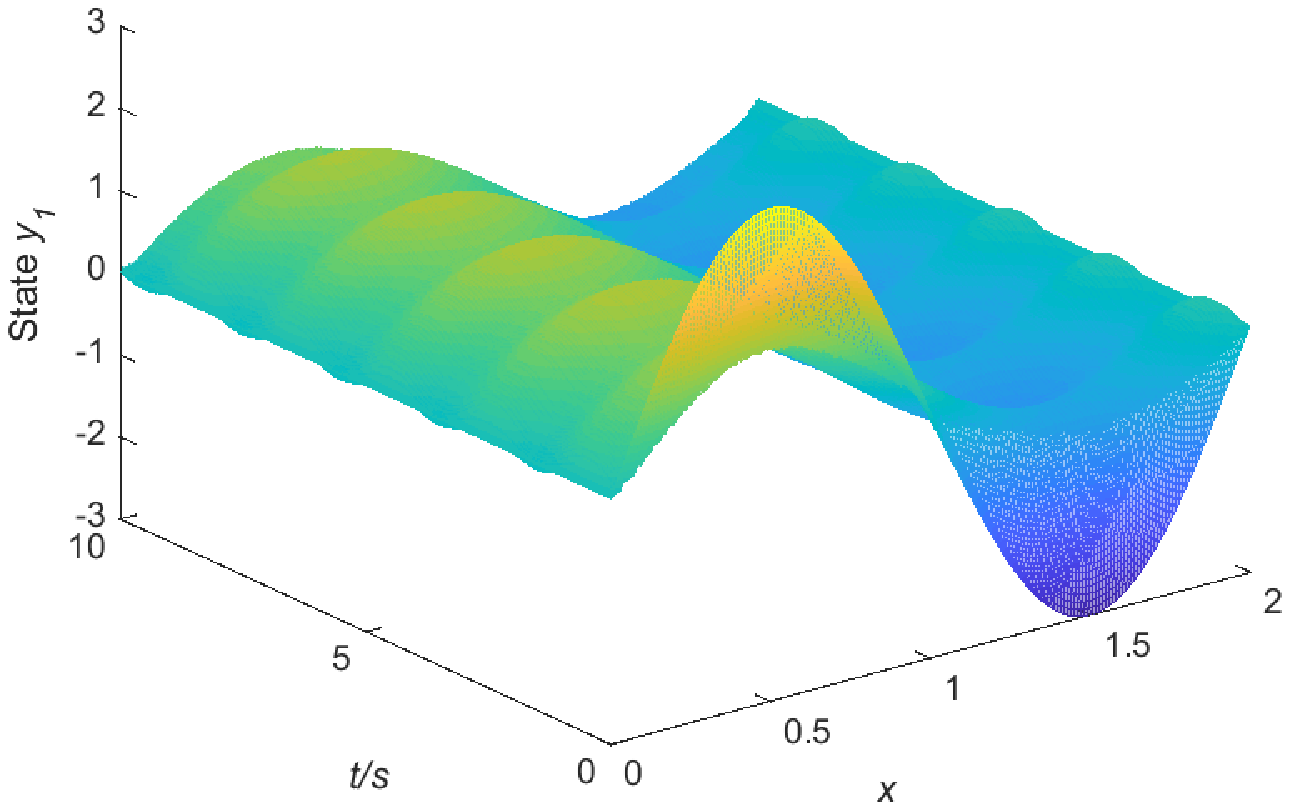}
   \includegraphics[width=4cm]{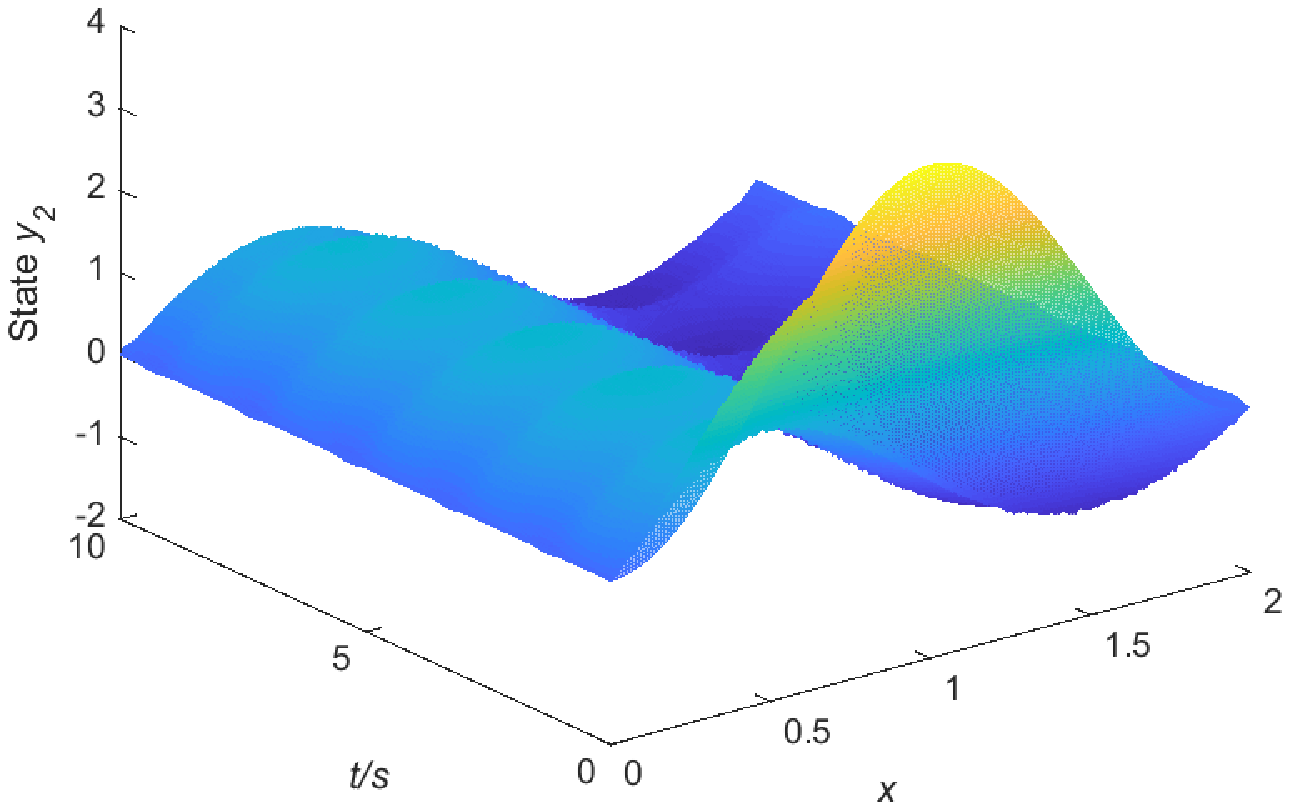}
   \includegraphics[width=4cm]{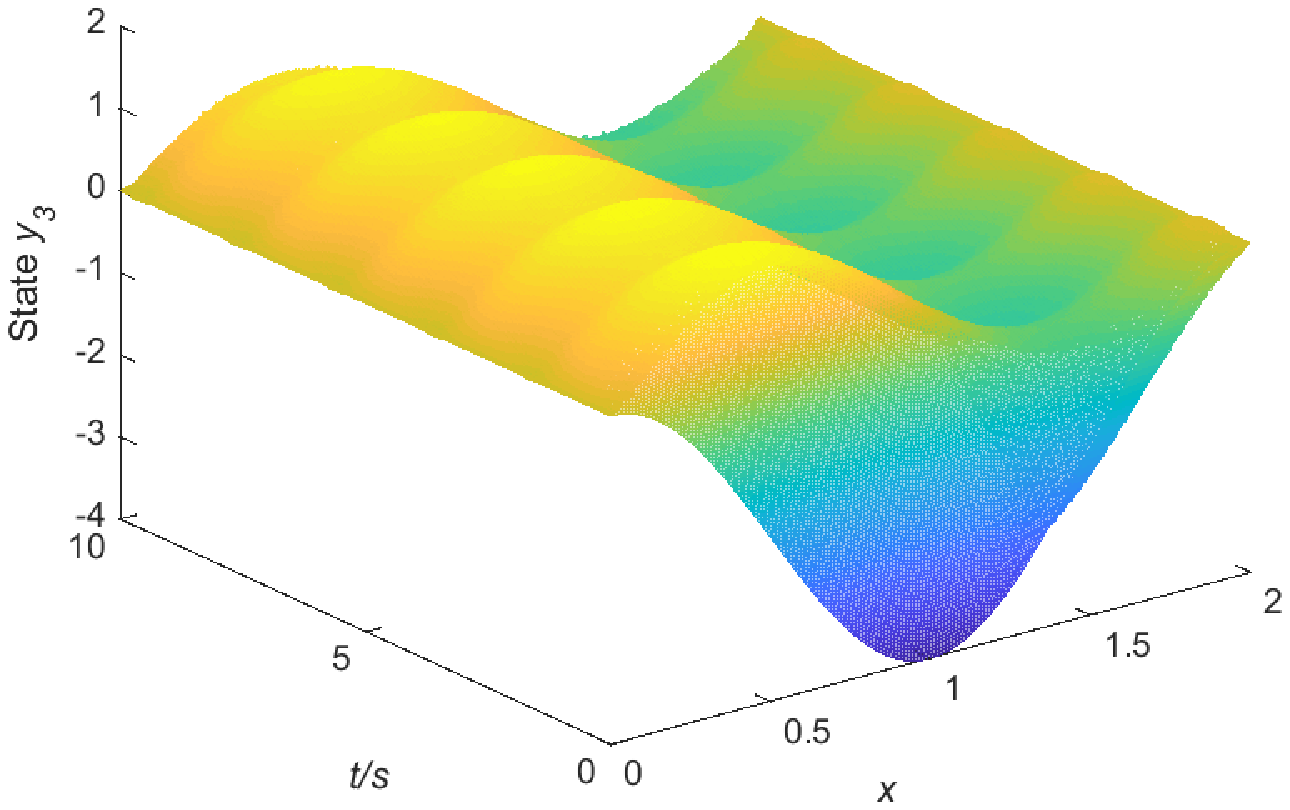}
   \caption{The states of MASs \eqref{e4.9} with \eqref{e4.13} and \eqref{e4.10} via the controller \eqref{e4.11}. }\label{Fig12}
 \end{figure}

\begin{figure}[H]
   \centering
   \includegraphics[width=4cm]{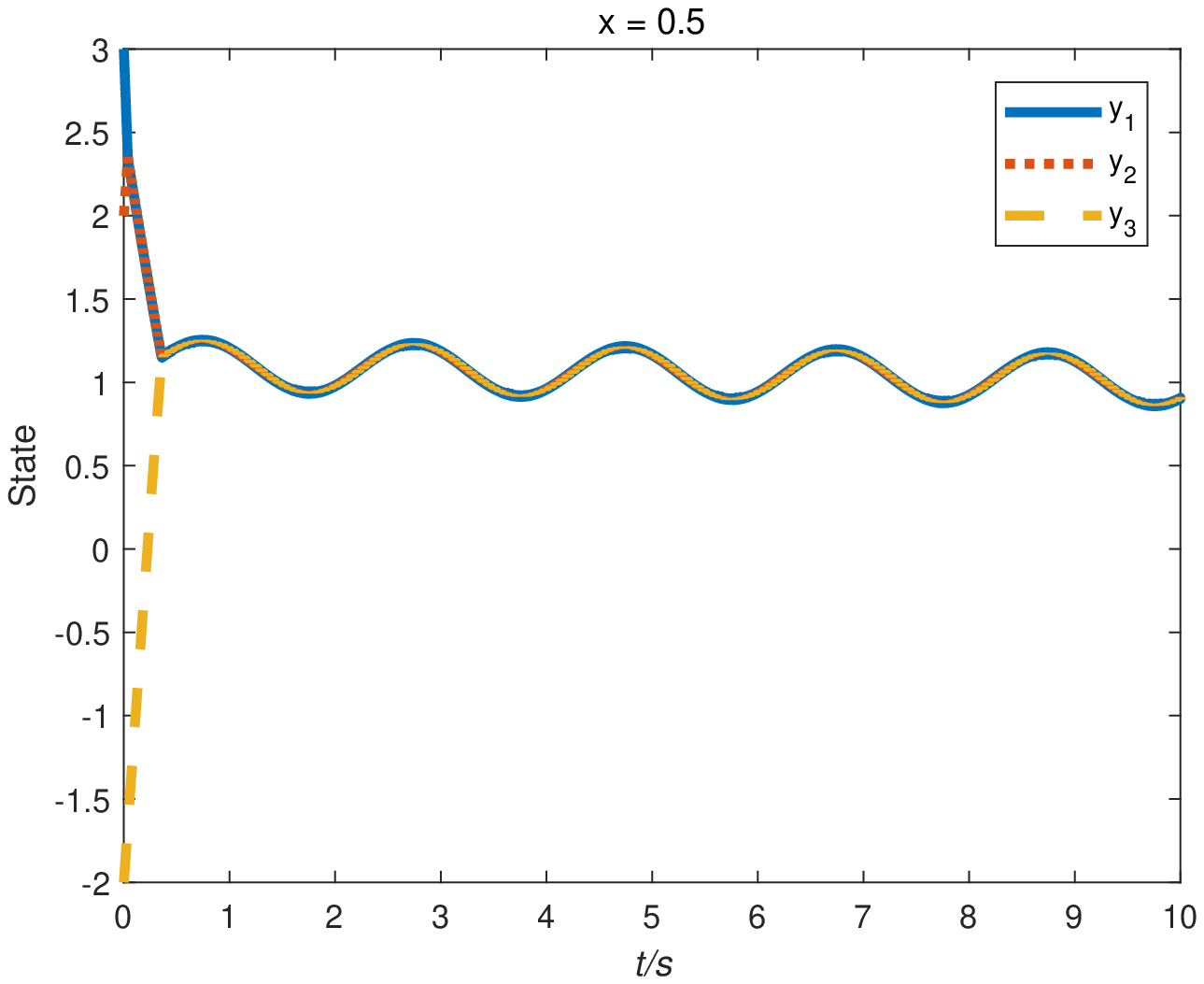}
   \includegraphics[width=4cm]{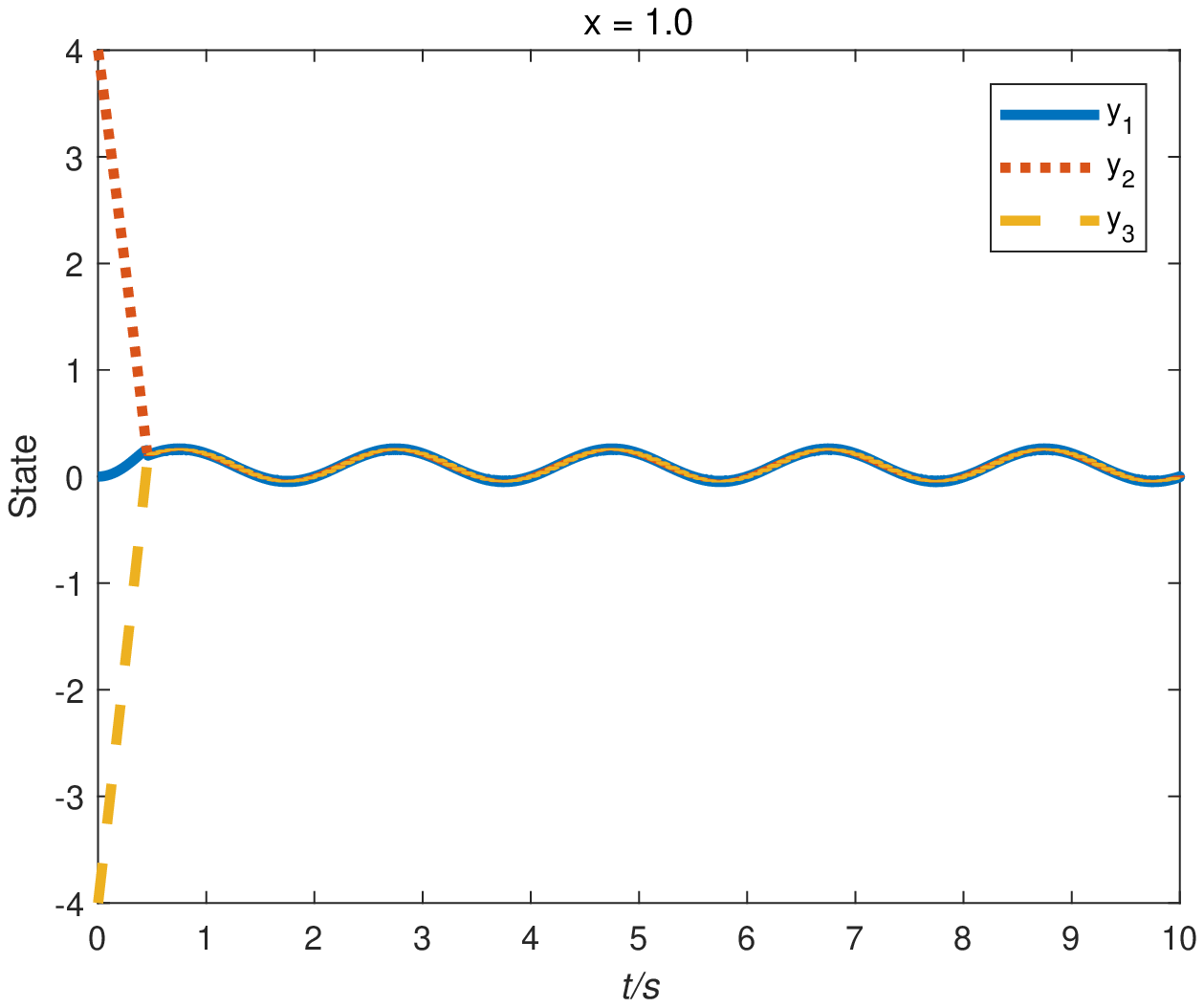}
   \includegraphics[width=4cm]{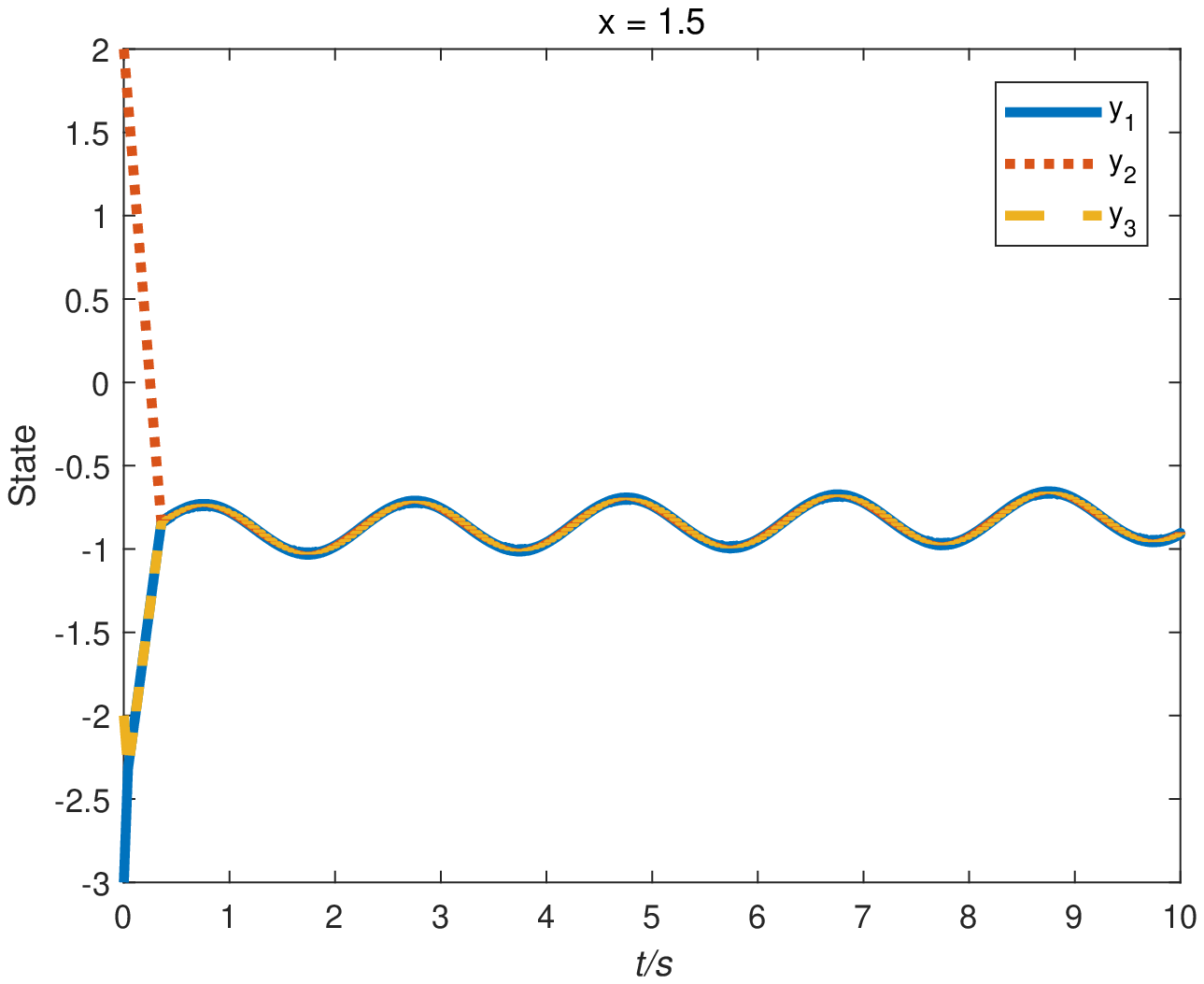}
   \caption{The sections of MASs \eqref{e4.9} with \eqref{e4.13} and \eqref{e4.10} via the controller \eqref{e4.11}.}\label{Fig13}
 \end{figure}

From Theorem \ref{t5}, it is easy to see that MASs \eqref{e4.9} with \eqref{e4.13} and \eqref{e4.10} via the controller \eqref{e4.11} can achieve the FTC at $t^*$ and $t^*\leq7.746$. Fig.\,\ref{Fig12} and Fig.\,\ref{Fig13} show the states and the sections of MASs \eqref{e4.9} with \eqref{e4.13} and \eqref{e4.10} via the controller \eqref{e4.11} at $x=0.5, 1.0, 1.5$, respectively. Fig.\,\ref{Fig16} displays the sections of MASs \eqref{e4.9} without control at $x=0.5, 1.0, 1.5$, respectively.  Thus, the controller \eqref{e4.11} can ensure the FTC of MASs \eqref{e4.9} with \eqref{e4.13} and \eqref{e4.10}.

\begin{figure}[H]
   \centering
   \includegraphics[width=4cm]{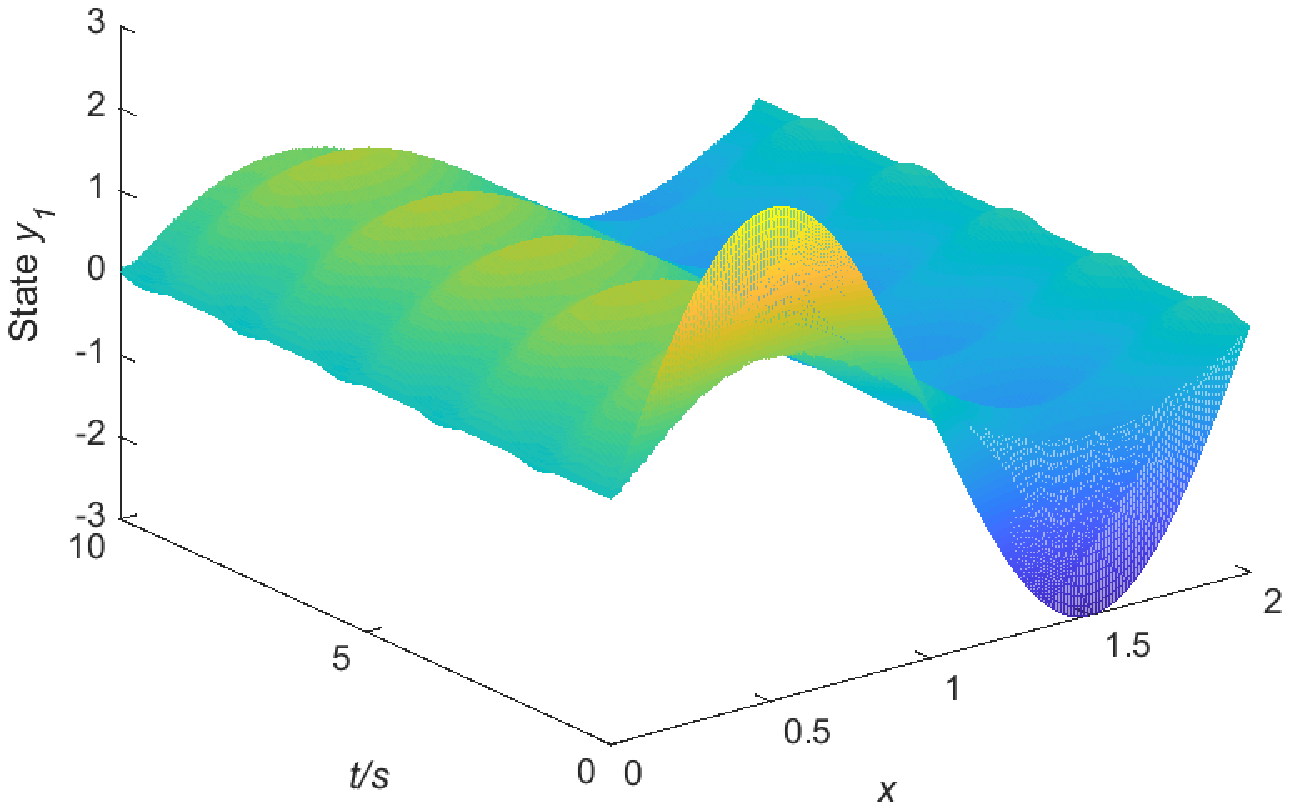}
   \includegraphics[width=4cm]{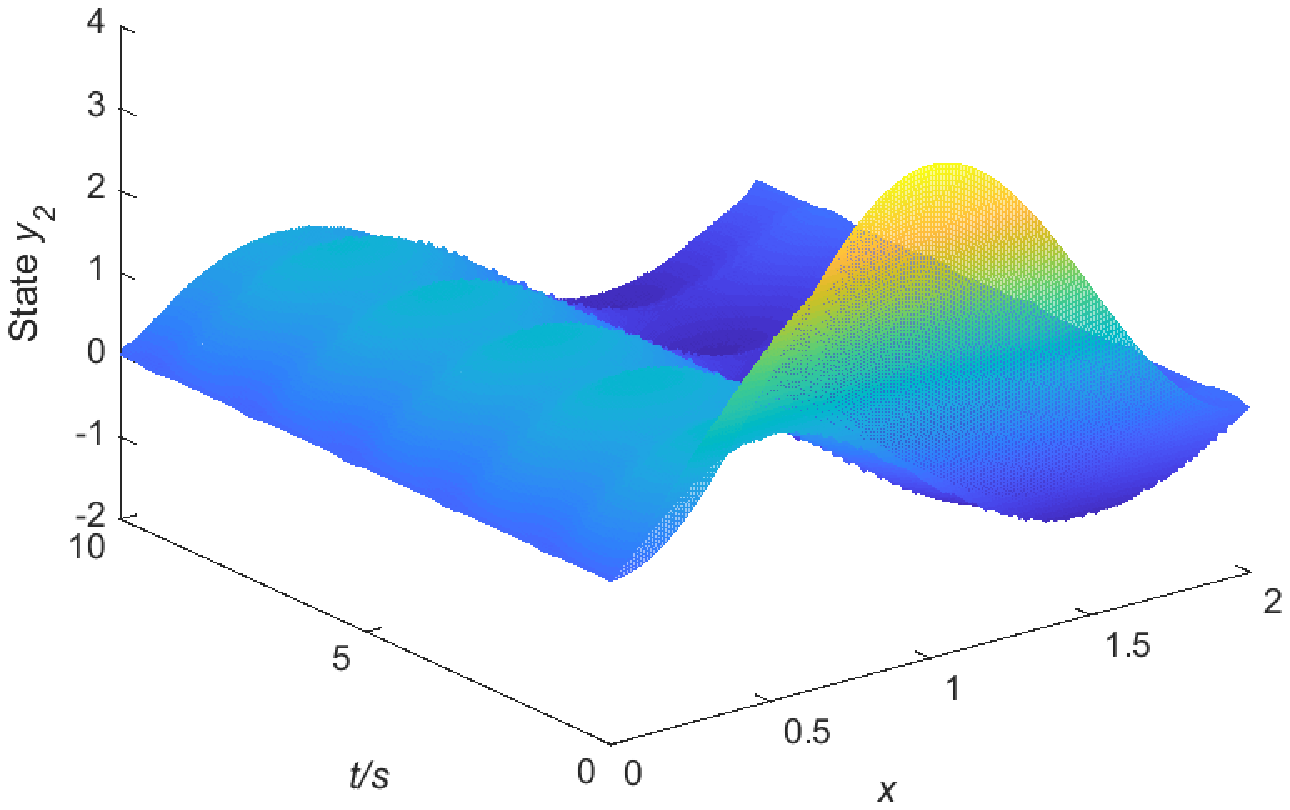}
   \includegraphics[width=4cm]{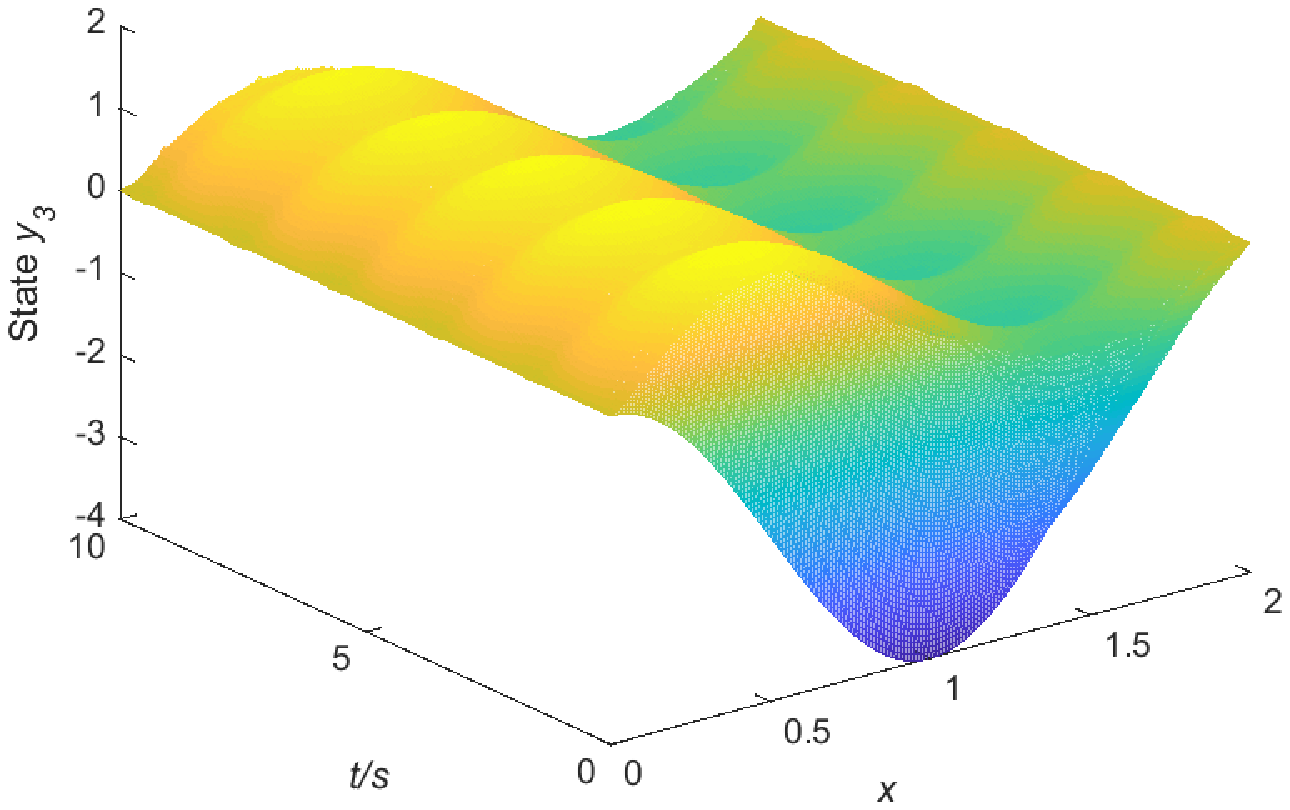}
   \caption{The states of MASs \eqref{e4.9} with \eqref{e4.13} and \eqref{e4.10} via the controller \eqref{e4.12}. }\label{Fig14}
 \end{figure}

\begin{figure}[H]
   \centering
   \includegraphics[width=4cm]{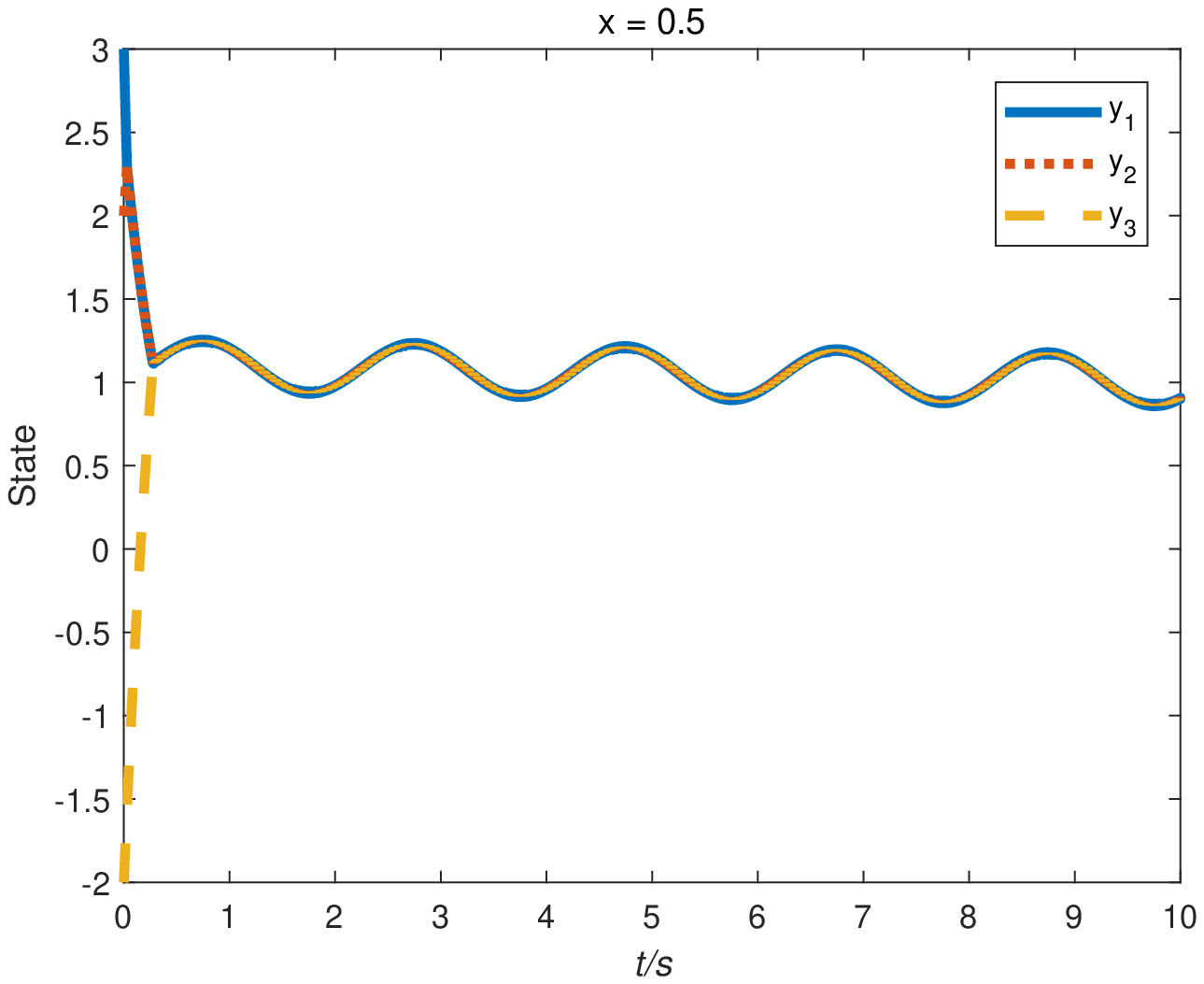}
   \includegraphics[width=4cm]{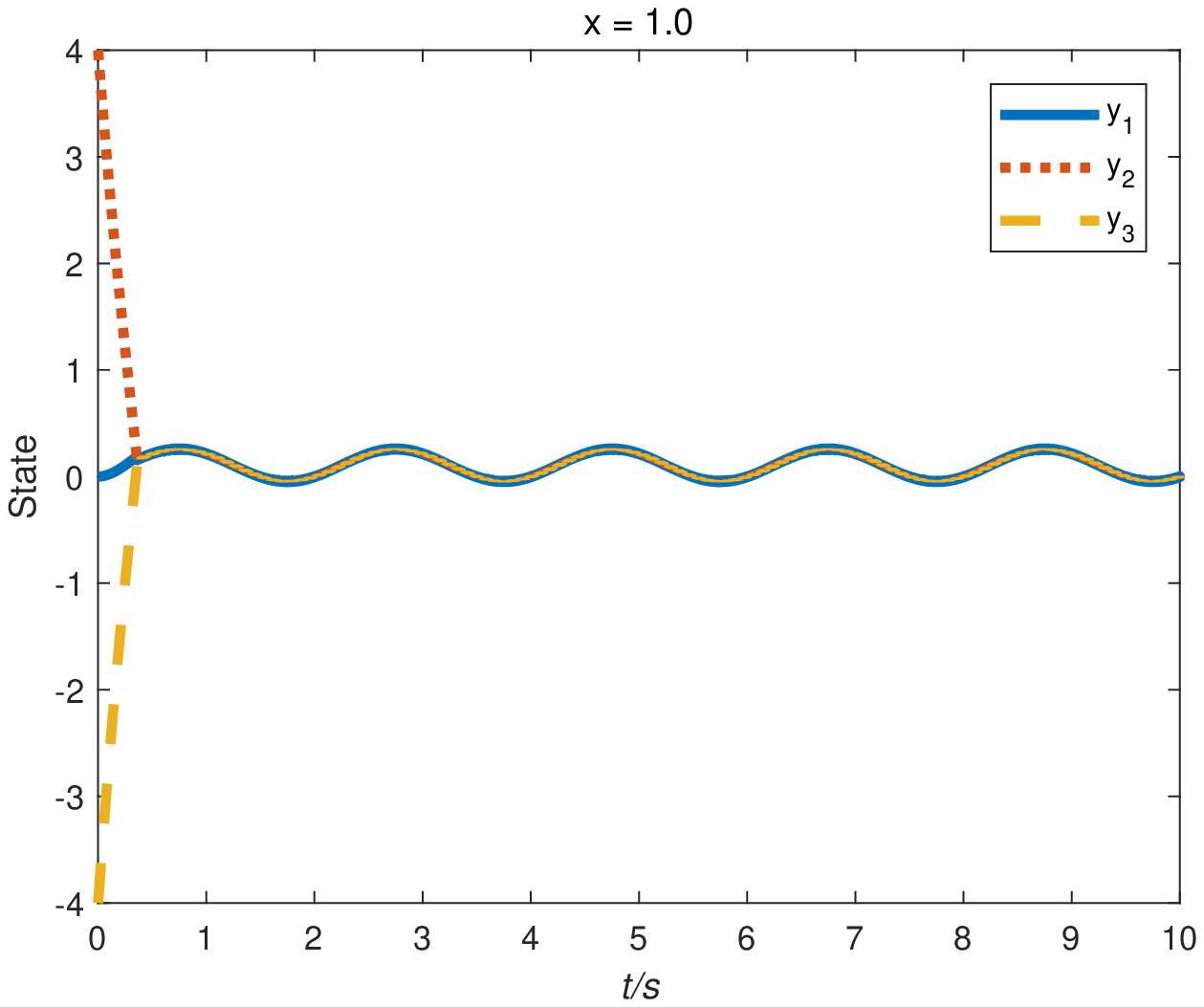}
   \includegraphics[width=4cm]{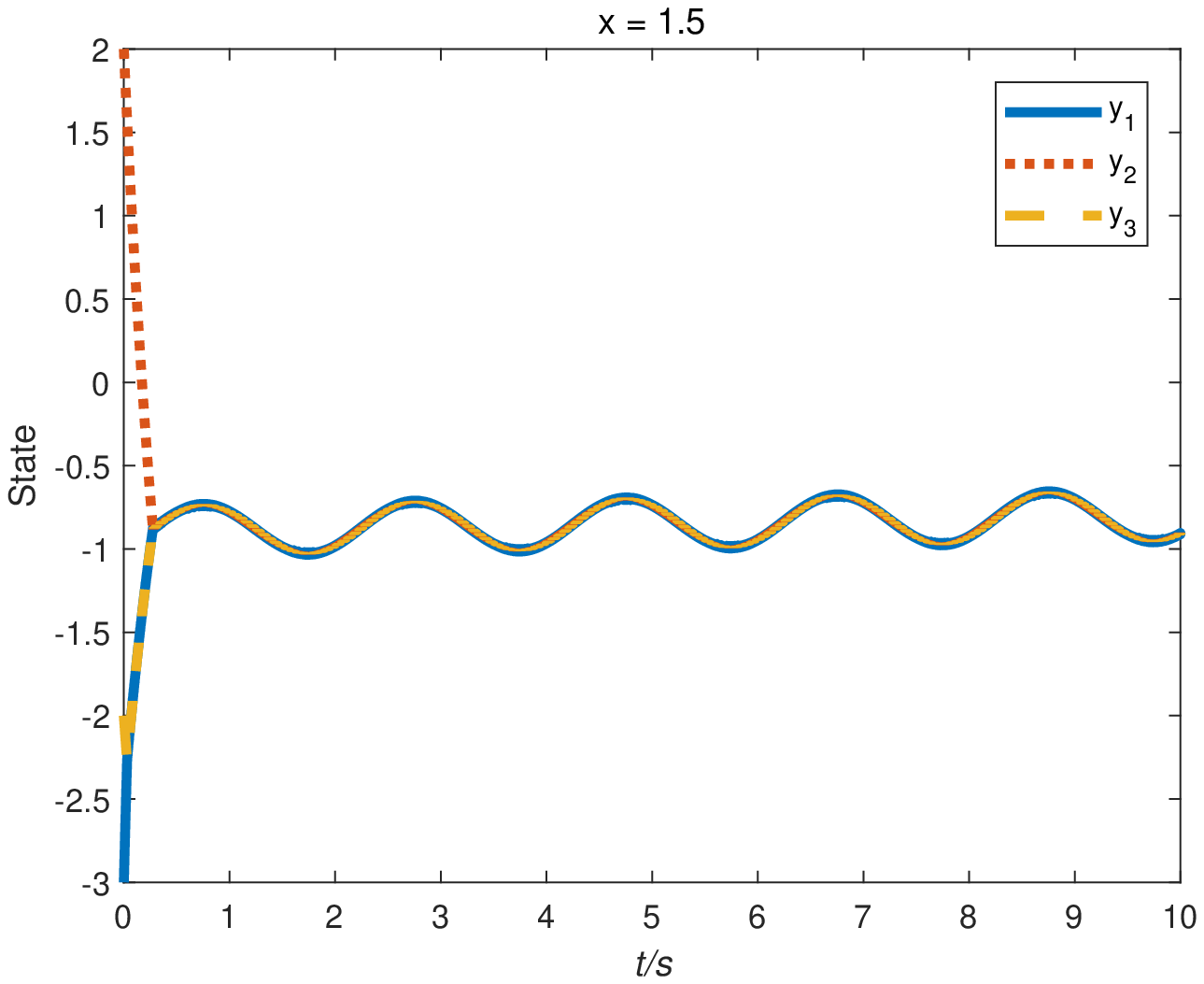}
   \caption{The sections of MASs \eqref{e4.9} with \eqref{e4.13} and \eqref{e4.10} via the controller \eqref{e4.12}.}\label{Fig15}
 \end{figure}

\begin{figure}[H]
   \centering
   \includegraphics[width=4cm]{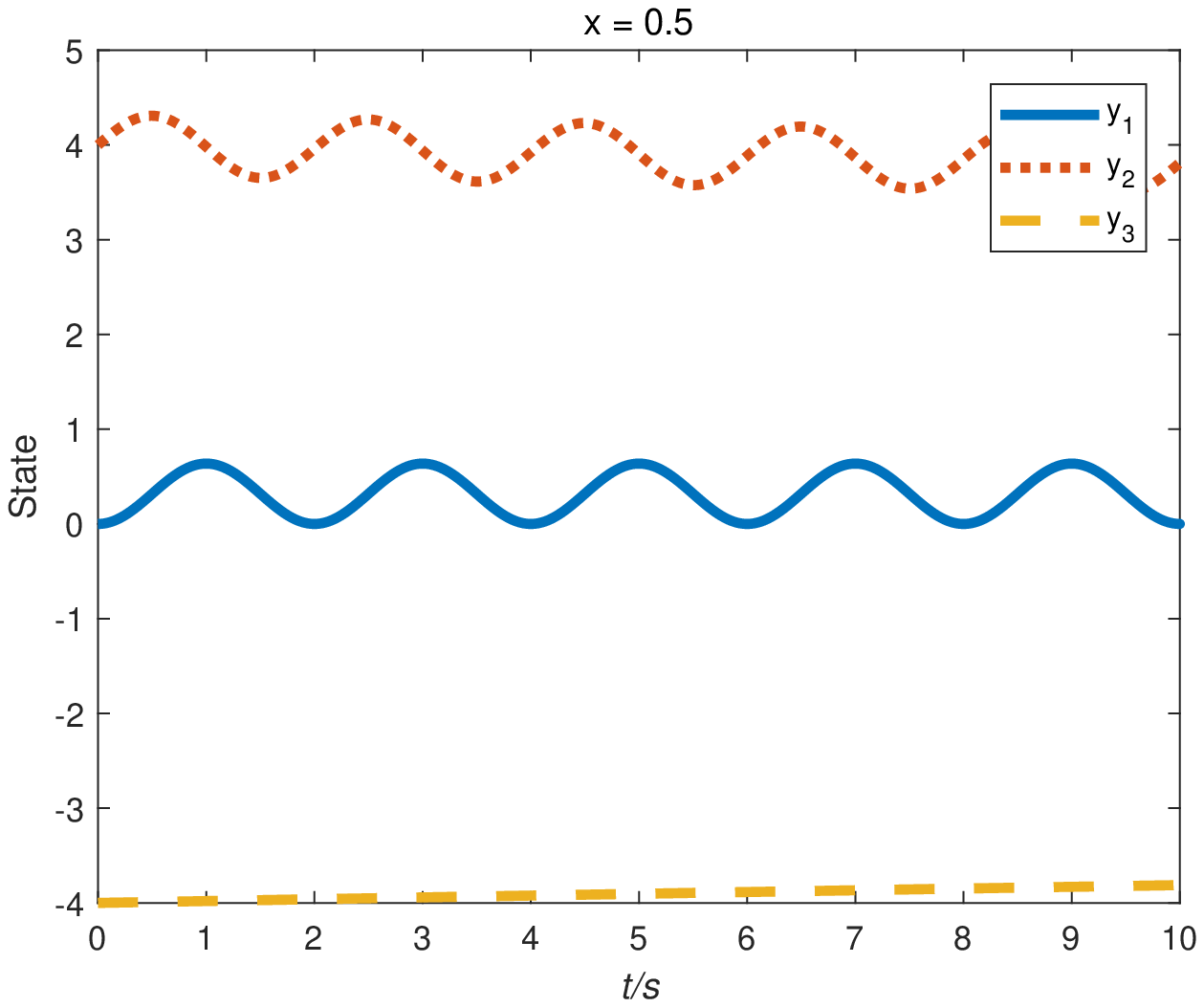}
   \includegraphics[width=4cm]{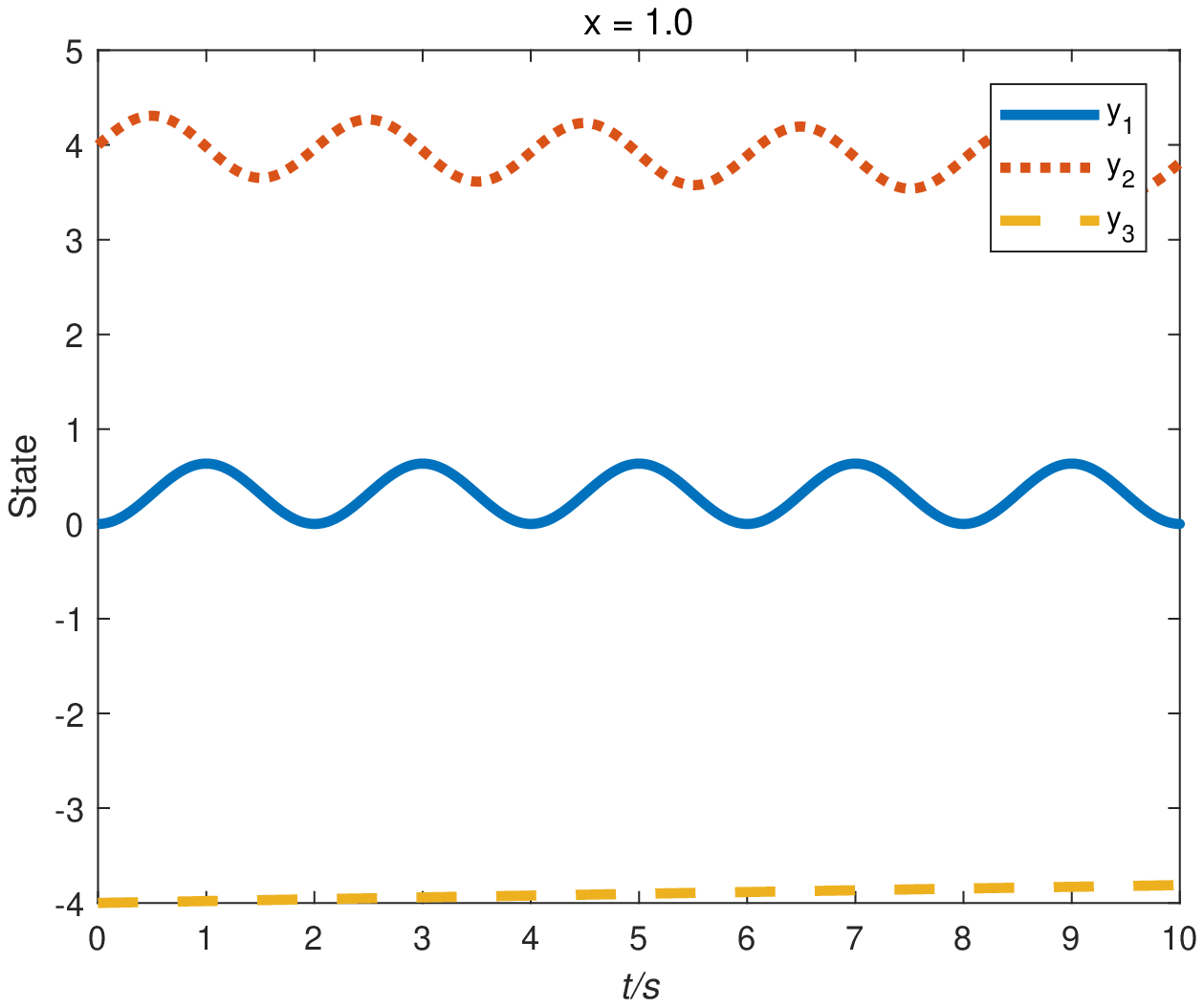}
   \includegraphics[width=4cm]{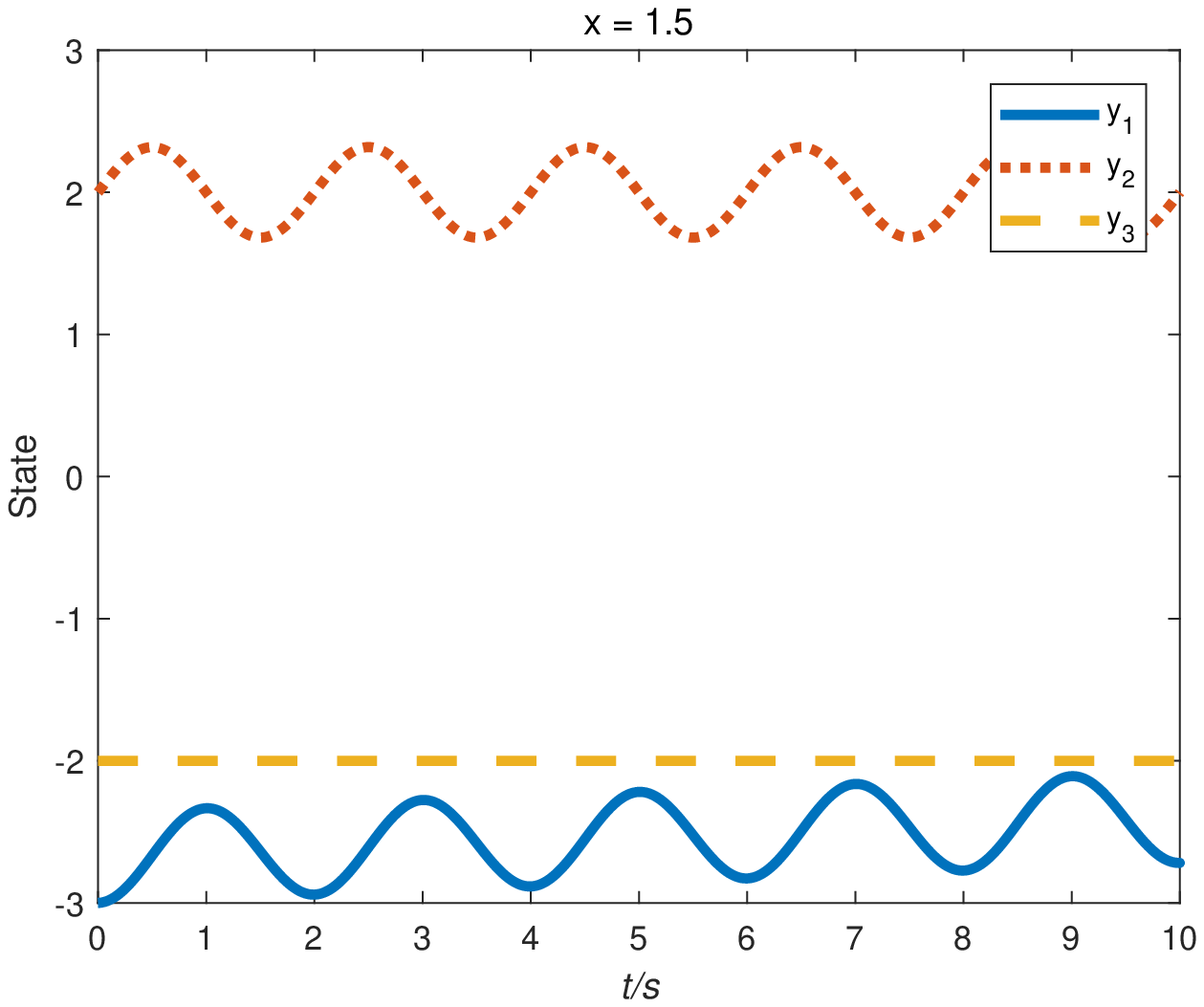}
   \caption{The sections of MASs \eqref{e4.9} with \eqref{e4.13} and \eqref{e4.10} without control.}\label{Fig16}
 \end{figure}

It follows from Theorem \ref{t6} that MASs \eqref{e4.9} with \eqref{e4.13} and \eqref{e4.10} via the controller \eqref{e4.12} can reach the FXC at $t^*$ and $t^*\leq T_{\max}= 7.226$. Fig.\,\ref{Fig14} and Fig.\,\ref{Fig15} indicate the states and the sections of MASs \eqref{e4.9} with \eqref{e4.13} and \eqref{e4.10} via the controller \eqref{e4.12} at $x=0.5, 1.0, 1.5$, respectively. Therefore, the controller \eqref{e4.12} can guarantee the FXC of MASs \eqref{e4.9} with \eqref{e4.13} and \eqref{e4.10}.

\begin{figure}[H]
   \centering
   \includegraphics[width=4cm]{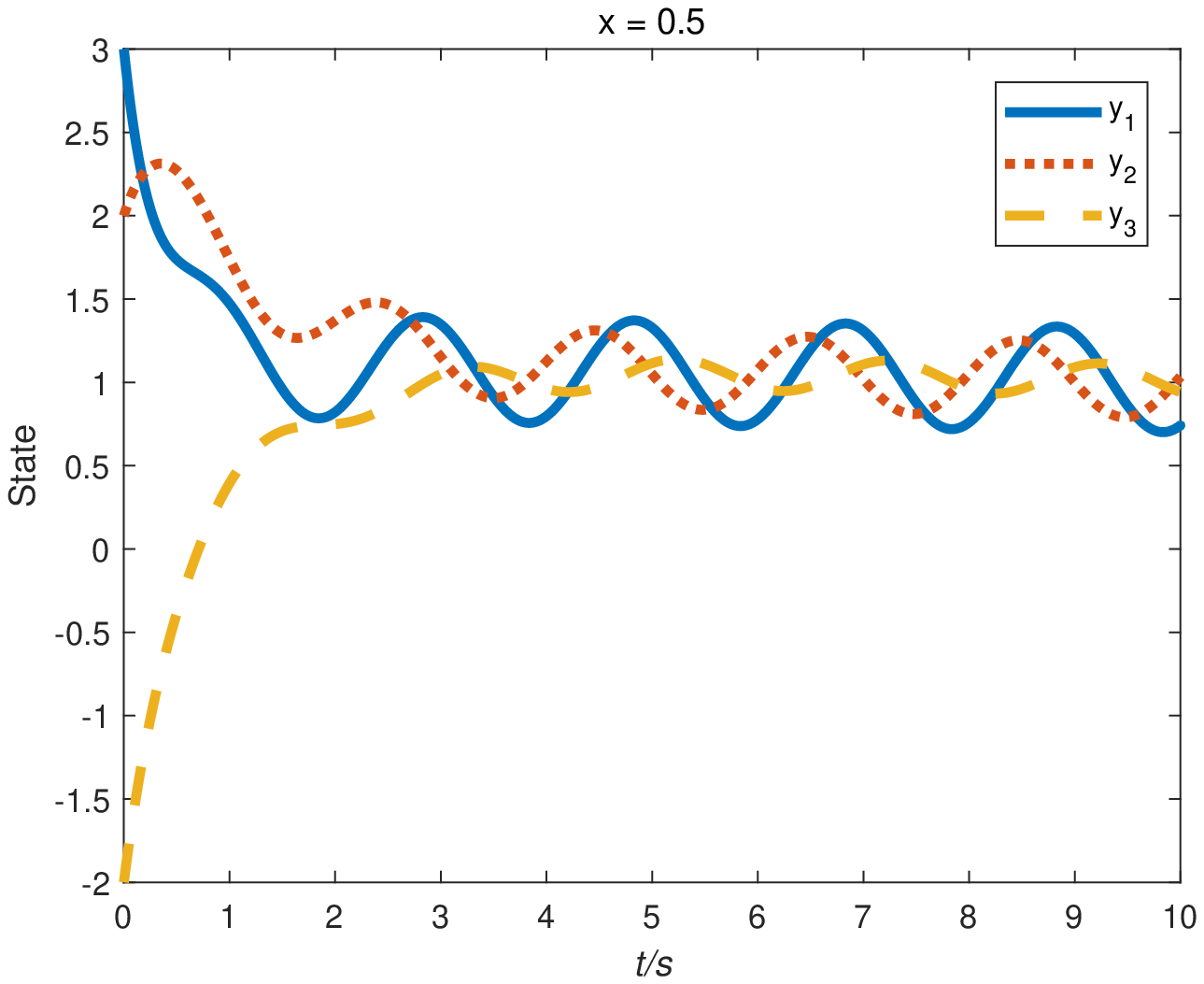}
   \includegraphics[width=4cm]{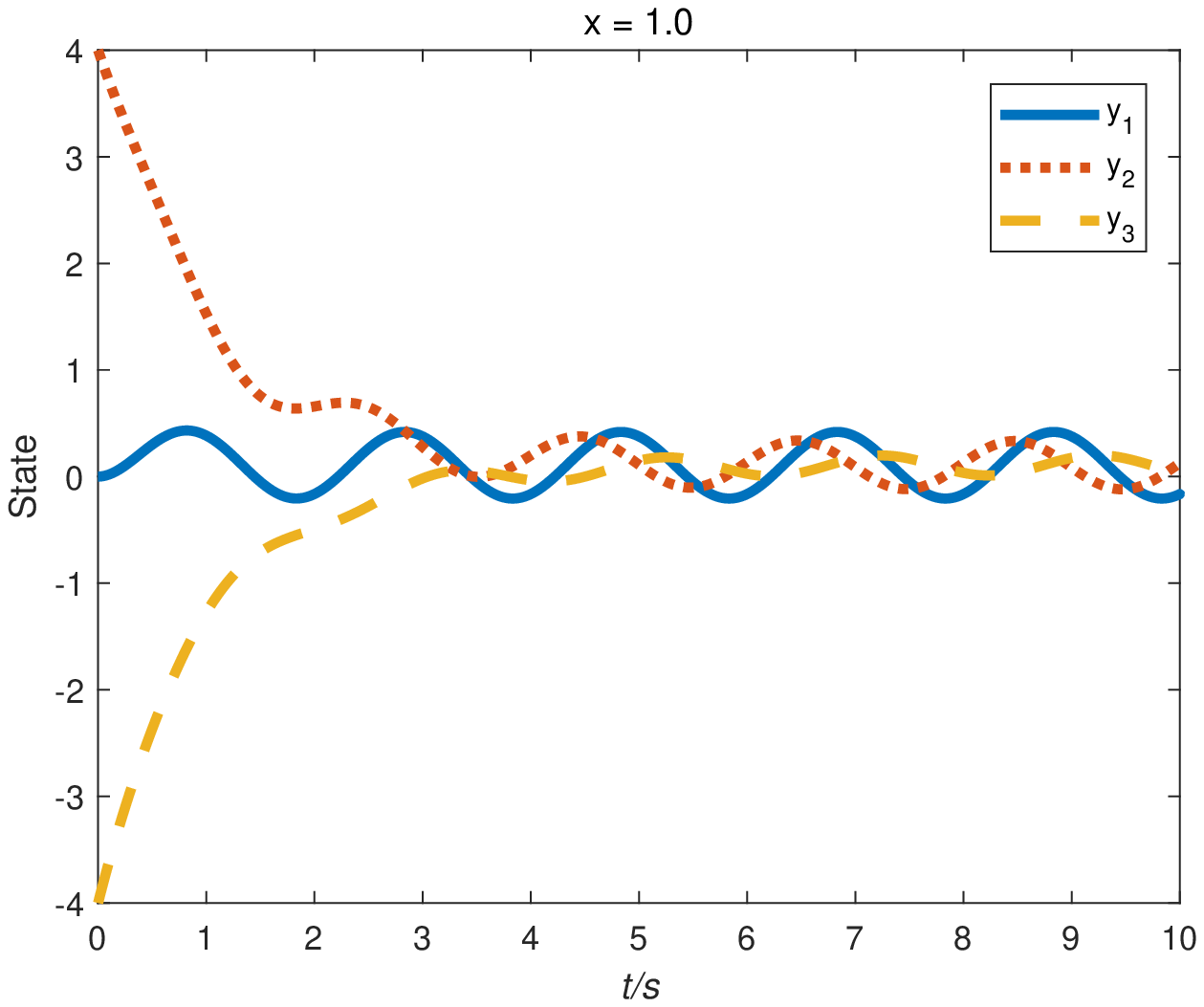}
   \includegraphics[width=4cm]{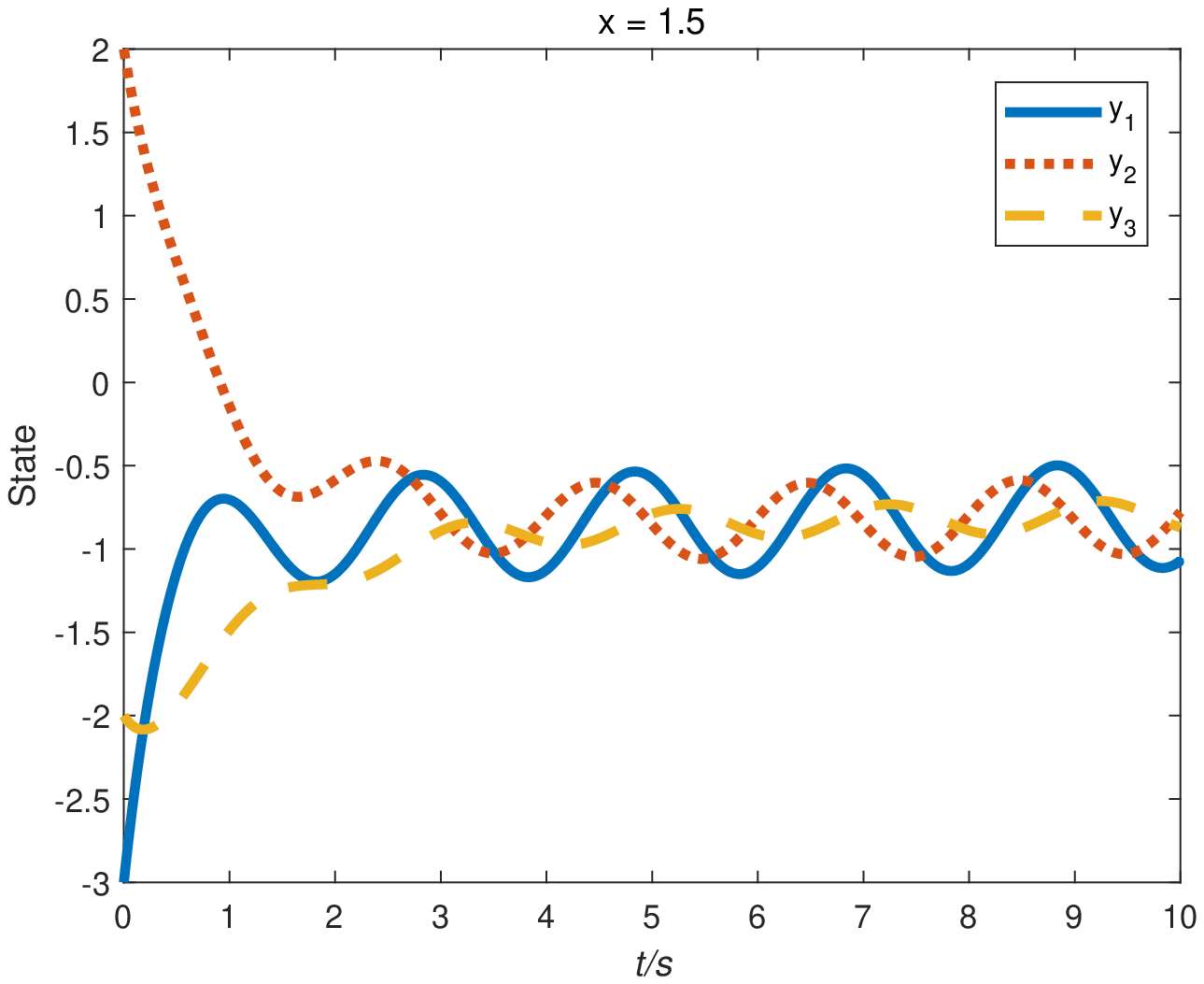}
   \caption{The sections of MASs \eqref{e4.9} with \eqref{e4.13} and \eqref{e4.10} via the controller \eqref{e4.c1}.}\label{Fig17}
 \end{figure}

\begin{figure}[H]
   \centering
   \includegraphics[width=4cm]{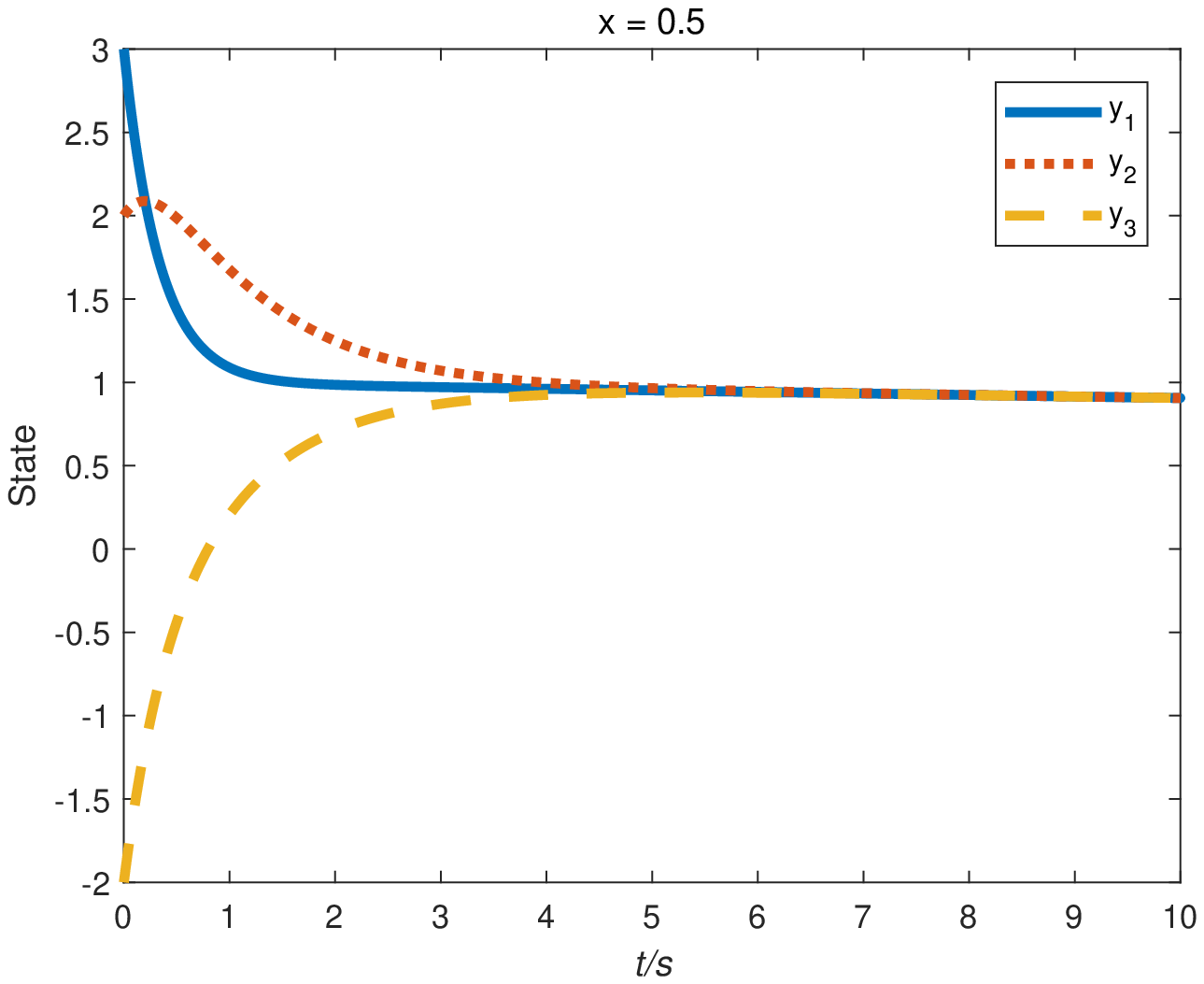}
   \includegraphics[width=4cm]{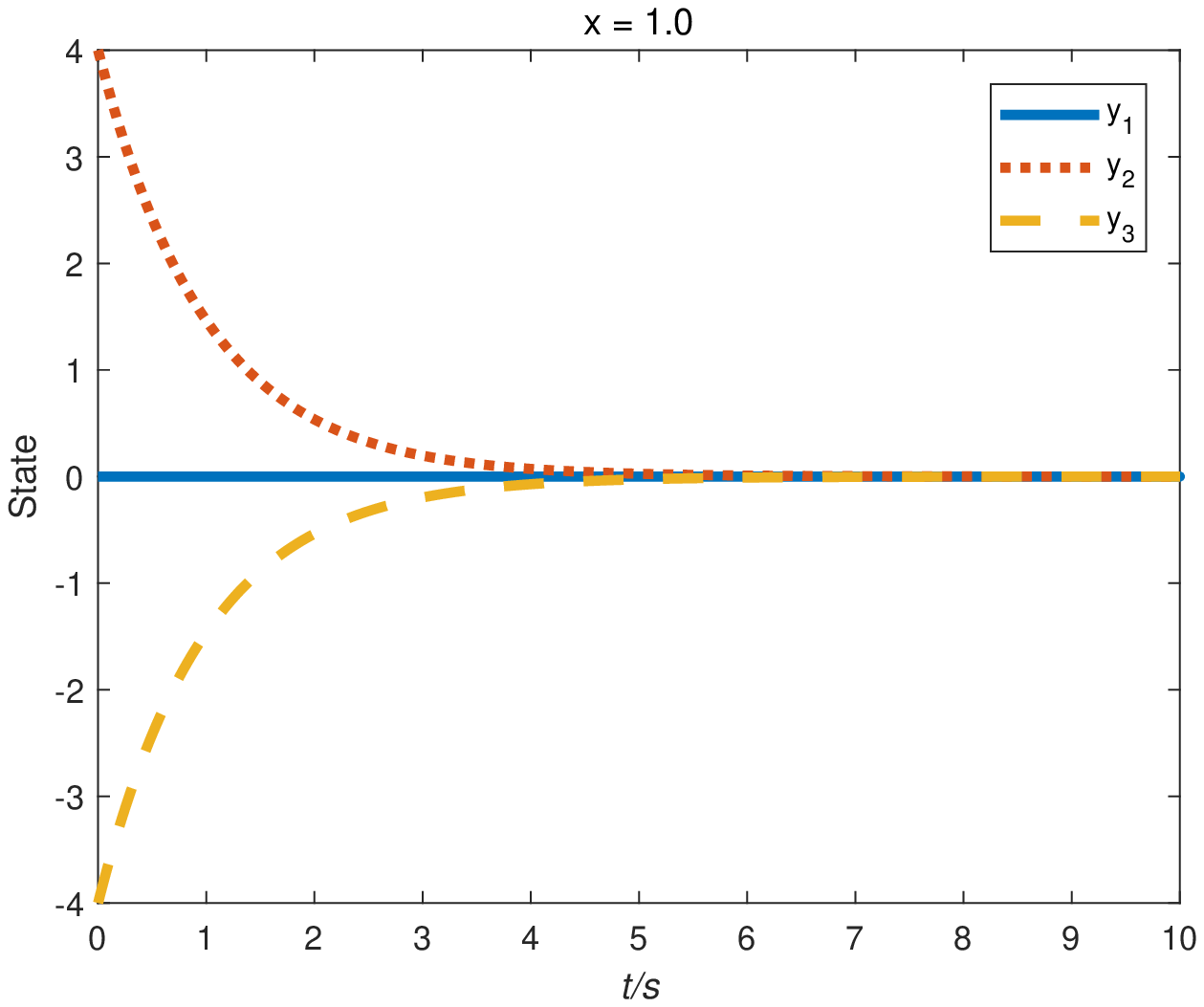}
   \includegraphics[width=4cm]{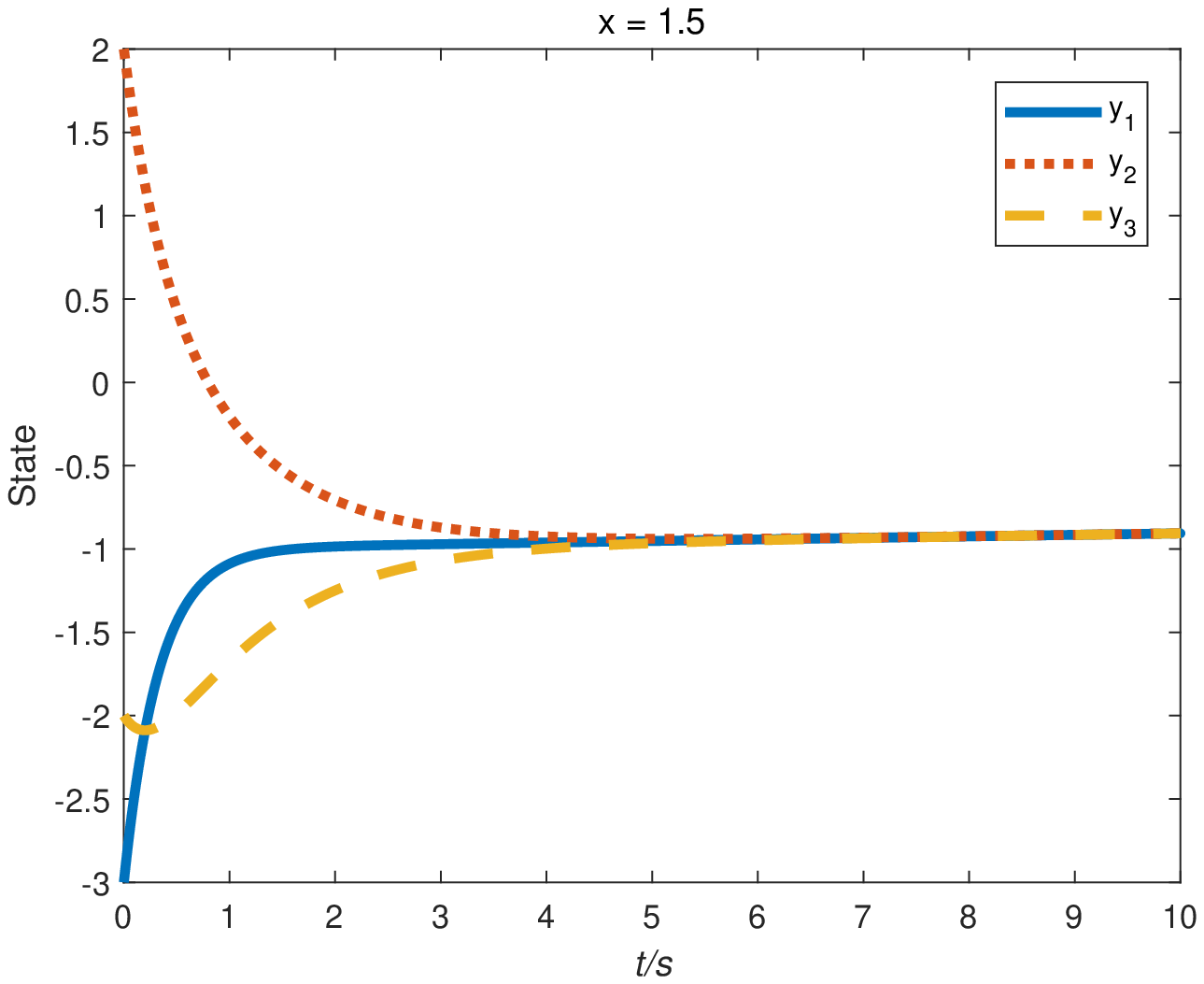}
   \caption{The sections of MASs \eqref{e4.9} with \eqref{e4.13} and without \eqref{e4.10} via the controller \eqref{e4.c1}.}\label{Fig18}
 \end{figure}

\begin{figure}[H]
   \centering
   \includegraphics[width=4cm]{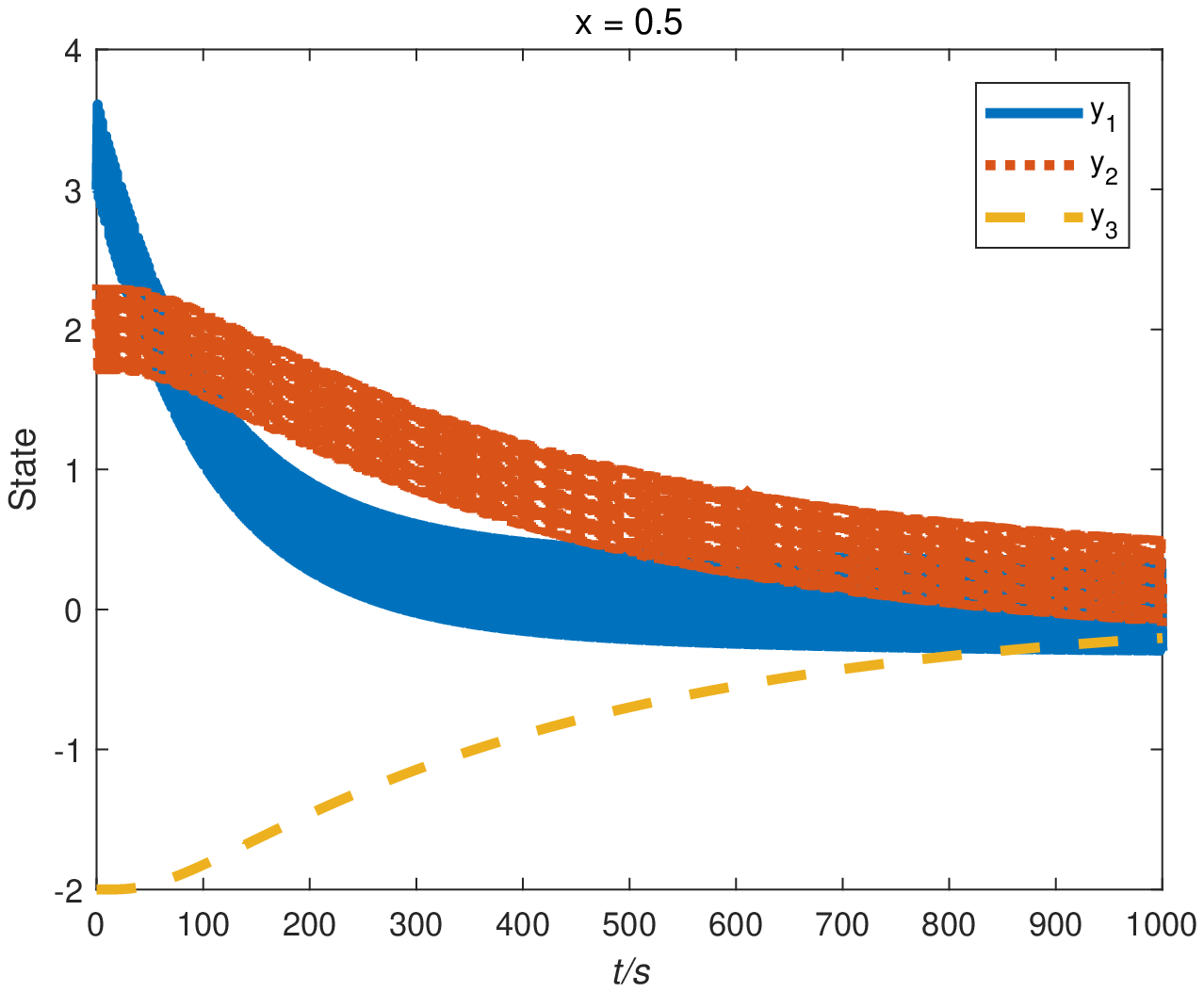}
   \includegraphics[width=4cm]{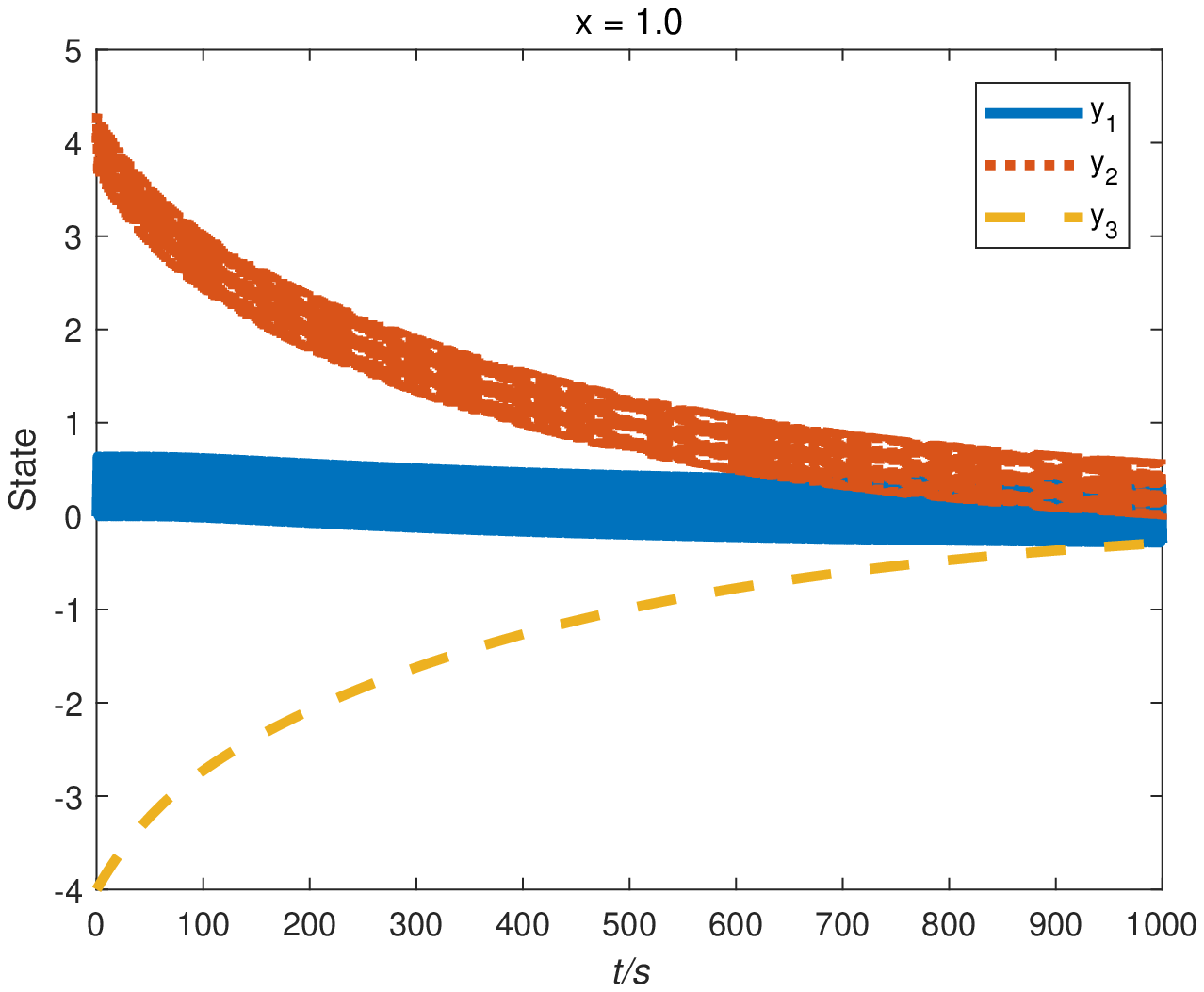}
   \includegraphics[width=4cm]{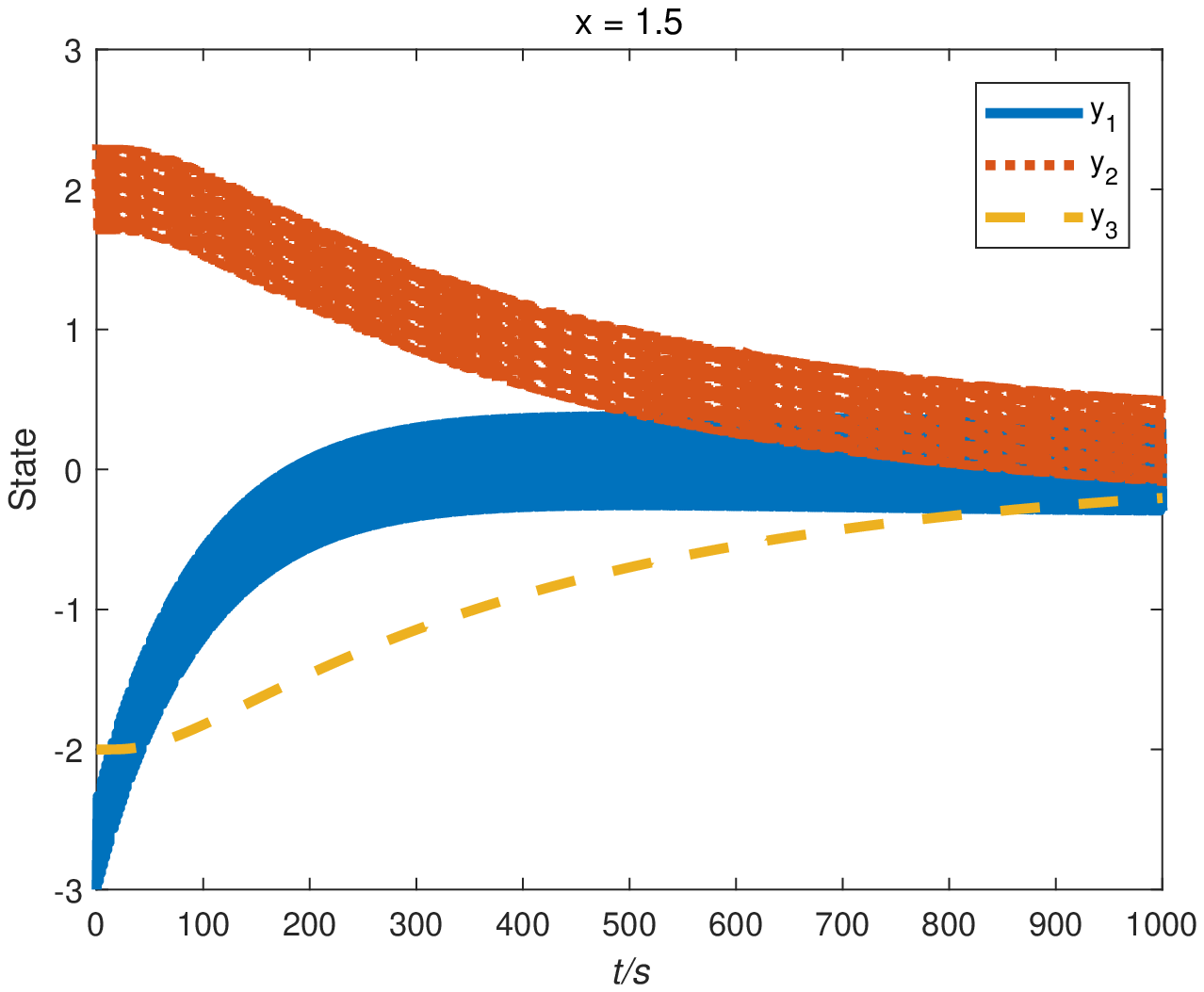}
   \caption{The sections of MASs \eqref{e4.9} with \eqref{e4.13} and \eqref{e4.10} via the controller \eqref{e4.c2}.}\label{Fig19}
 \end{figure}

\begin{figure}[H]
   \centering
   \includegraphics[width=4cm]{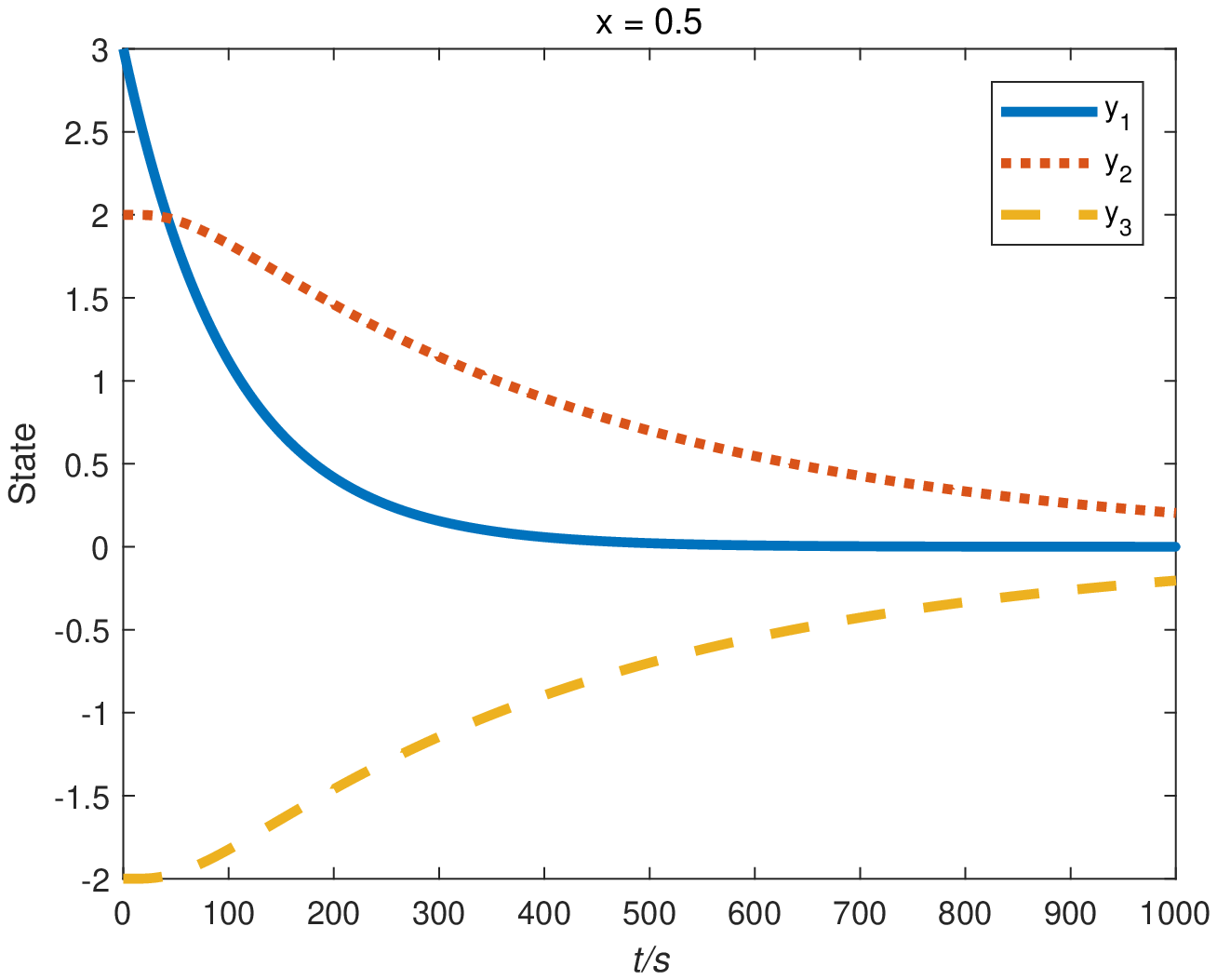}
   \includegraphics[width=4cm]{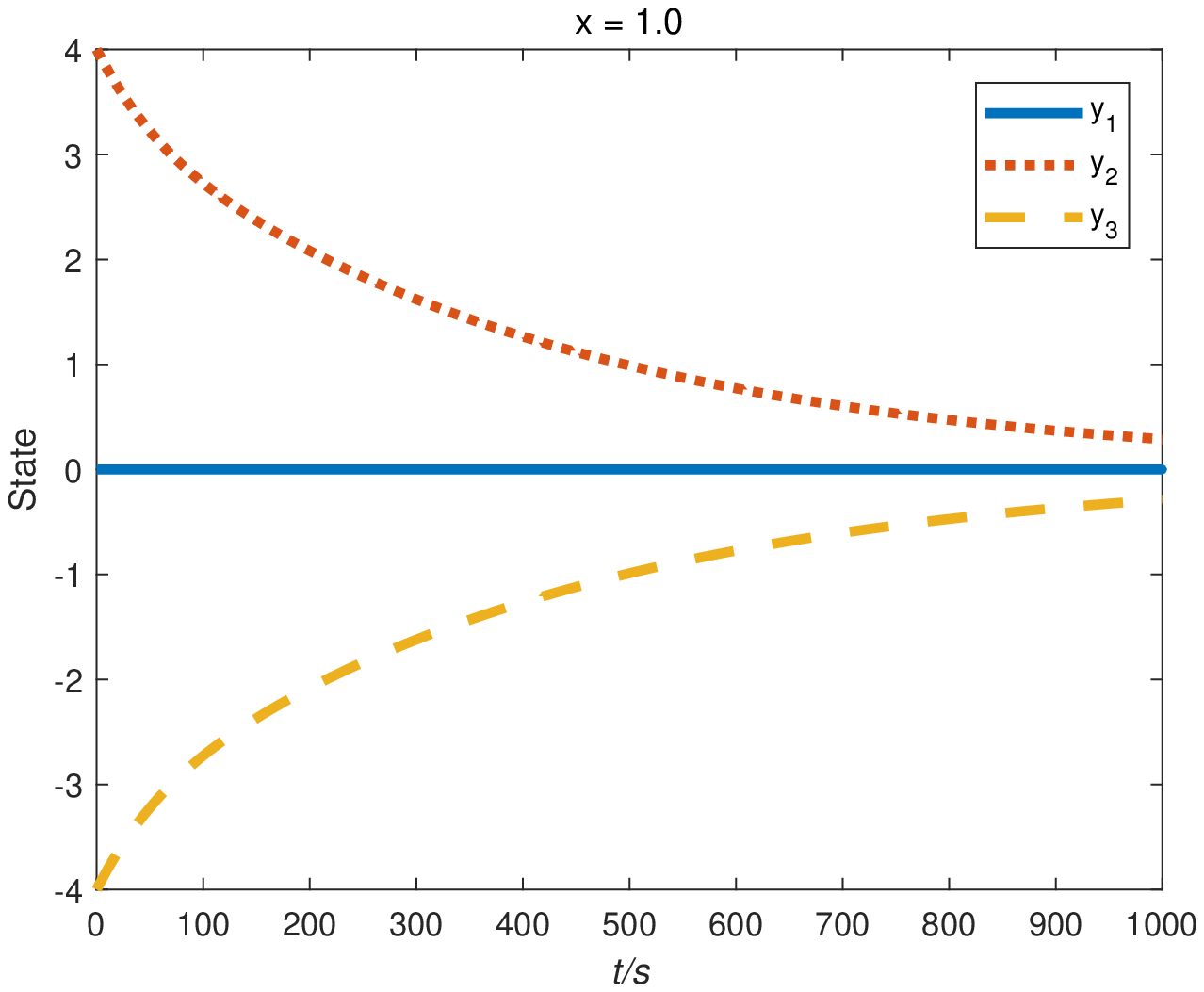}
   \includegraphics[width=4cm]{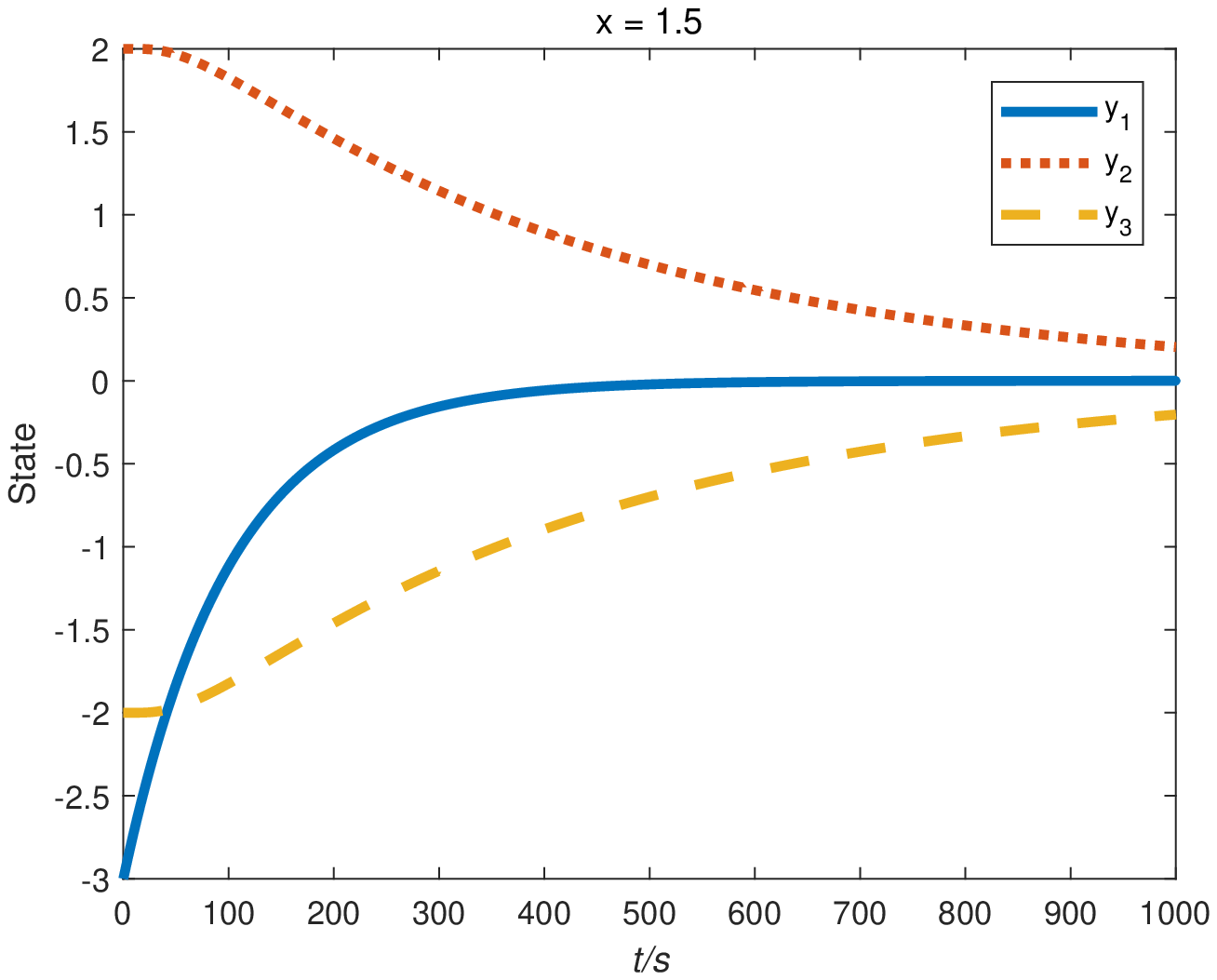}
   \caption{The sections of MASs \eqref{e4.9} with \eqref{e4.13} and without \eqref{e4.10} via the controller \eqref{e4.c2}.}\label{Fig20}
 \end{figure}

The sections of MASs \eqref{e4.9} with \eqref{e4.13} and \eqref{e4.10} via the controllers \eqref{e4.c1} and \eqref{e4.c2} at $x=0.5, 1.0, 1.5$ are shown in Fig.\,\ref{Fig17} and Fig.\,\ref{Fig19}, respectively, which indicate that the controllers \eqref{e4.c1} and \eqref{e4.c2} cannot handle the presence of external interference. When $d_1(x, t)=d_2(x, t)=d_3(x, t)=0$, the sections of MASs \eqref{e4.9} via the controllers \eqref{e4.c1} and \eqref{e4.c2} are shown in Fig.\,\ref{Fig18} and Fig.\,\ref{Fig20}, respectively. By Fig.\,\ref{Fig13}, Fig.\,\ref{Fig15}, Fig.\,\ref{Fig18} and Fig.\,\ref{Fig20}, we can see that the controllers \eqref{e4.11} and \eqref{e4.12} have faster convergence speed than the controllers \eqref{e4.c1} and \eqref{e4.c2}.

\begin{figure}[H]
   \centering
   \includegraphics[width=5cm]{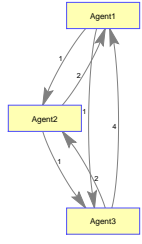}
   \caption{A directed, s-con and d-bal graph with 3 agents.}\label{Fig6}
 \end{figure}

\textbf{Example 2} \quad Consider MASs with $3$ agents, the graph $G$ is directed, s-con and d-bal given by Fig.\,\ref{Fig6}, where the number beside the edge means the adjacency element. Clearly, we can get the adjacency matrix $A=\left(\begin{smallmatrix}
           0 & 1 & 1 \\
           2 & 0 & 1 \\
           4 & 2 & 0
         \end{smallmatrix}\right)$,  and $\omega=[1, \frac{1}{2}, \frac{1}{4}]^T$.
When $\alpha = 0.5$, $\beta = 1.1$, $k=10^{-3} $ and $x \in [0 ,2]$, MASs can be obtained as
\begin{eqnarray}\label{e4.5}
\left\{
\begin{array}{llll}
y_{i,t}=10^{-3}\cdot y_{i,xx}+u_i,\\
y_{i}(0,t)= y_{i}(2,t)=0,\quad i=1, 2, 3,
\end{array}
\right.
\end{eqnarray}
the FTC controller can be formulated as
\begin{eqnarray}\label{e4.6}
u_i=-\sum_{j=1}^{3}a_{ij}\left(\int_0^2\xi_{ij}^{2}\mathrm{d}x\right)^{-0.5}\xi_{ij},
\end{eqnarray}
and the FXC controller can be obtained as follows
\begin{eqnarray}\label{e4.7}
u_i(x,t)=&-\sum_{j=1}^{3}a_{ij}\left(\int_0^2\xi_{ij}^{2}\mathrm{d}x\right)^{-0.5}\xi_{ij}\nonumber\\
&-\sum_{j=1}^{3}a_{ij}\left(\int_0^2\xi_{ij}^{2}\mathrm{d}x\right)^{0.1}\xi_{ij}.
\end{eqnarray}
The initial conditions are given by
\begin{eqnarray}\label{e4.8}
\left\{
\begin{array}{lll}
y_1(x, 0)=y^0_1(x)=3\sin(\pi x),\\
y_2(x, 0)=y^0_2(x)=-2\cos(\pi x)+2,\\
y_3(x, 0)=y^0_3(x)=2\cos(\pi x)-2.
\end{array}
\right.
\end{eqnarray}
It is worth mentioning that the distributed consensus controller considered in \cite{Fu18} was given by
\begin{eqnarray}\label{e41}
u_i=-\sum_{j=1}^{3}a_{ij}\xi_{ij}
\end{eqnarray}
and the consensus boundary controller studied in \cite{Pi16} was designed as follows
\begin{eqnarray}\label{e42}
y_{i,x}(2,t)=u_i(t)=-\sum_{j=1}^{3}a_{ij}\xi_{ij}(2,t).
\end{eqnarray}

\begin{figure}[H]
   \centering
   \includegraphics[width=4cm]{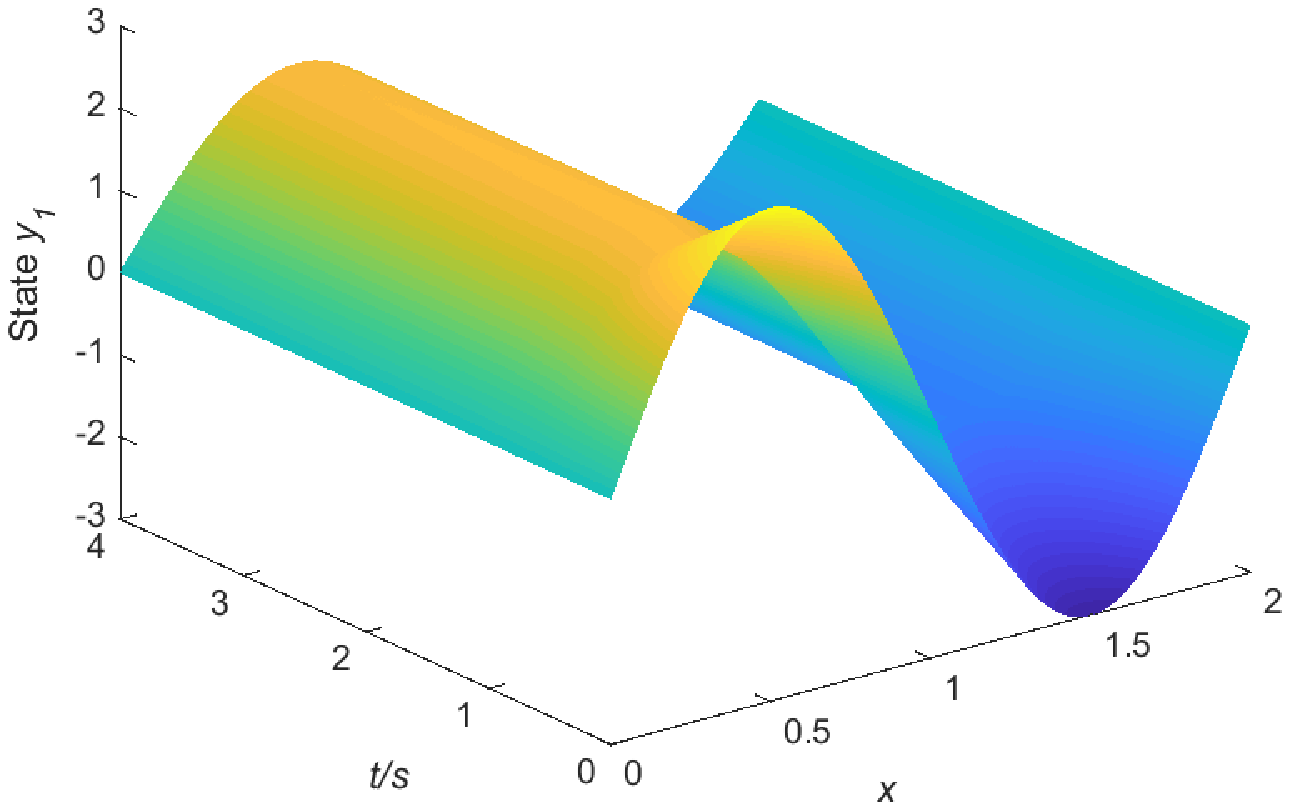}
   \includegraphics[width=4cm]{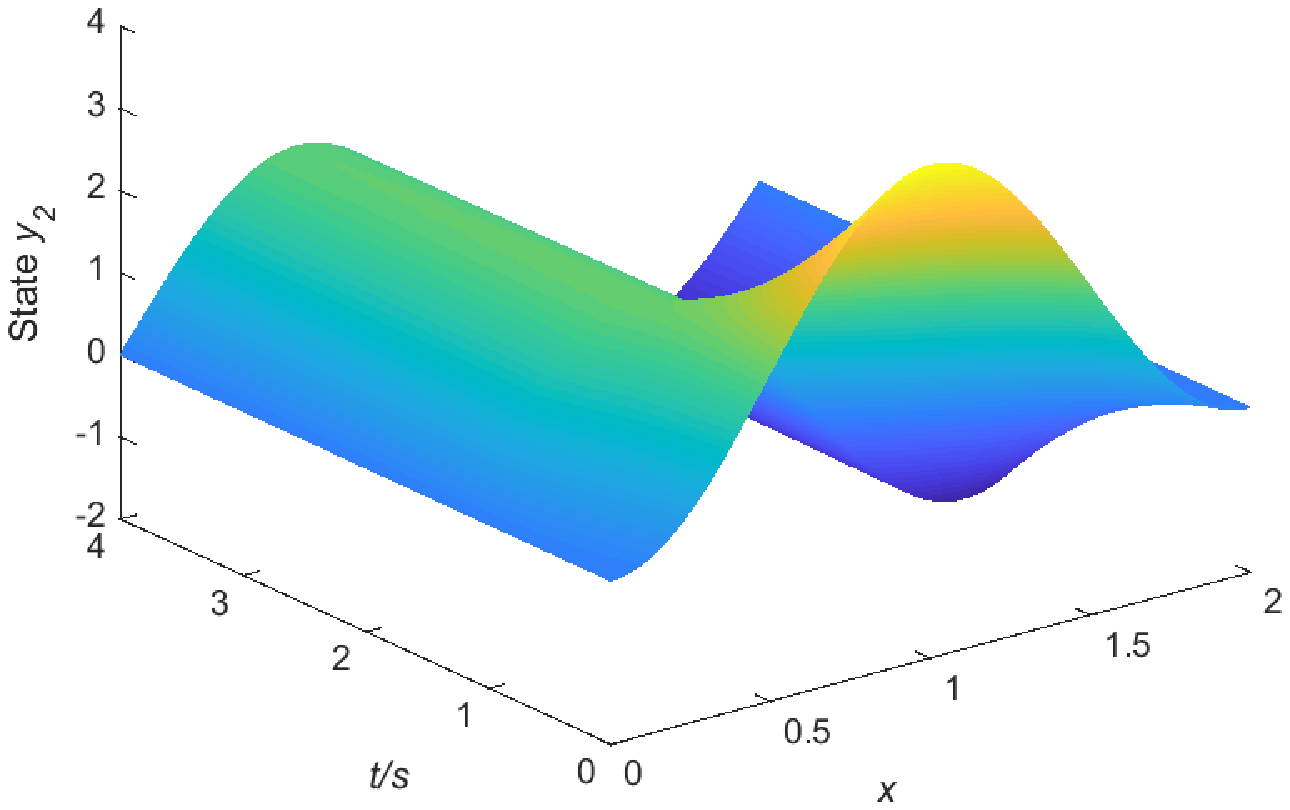}
   \includegraphics[width=4cm]{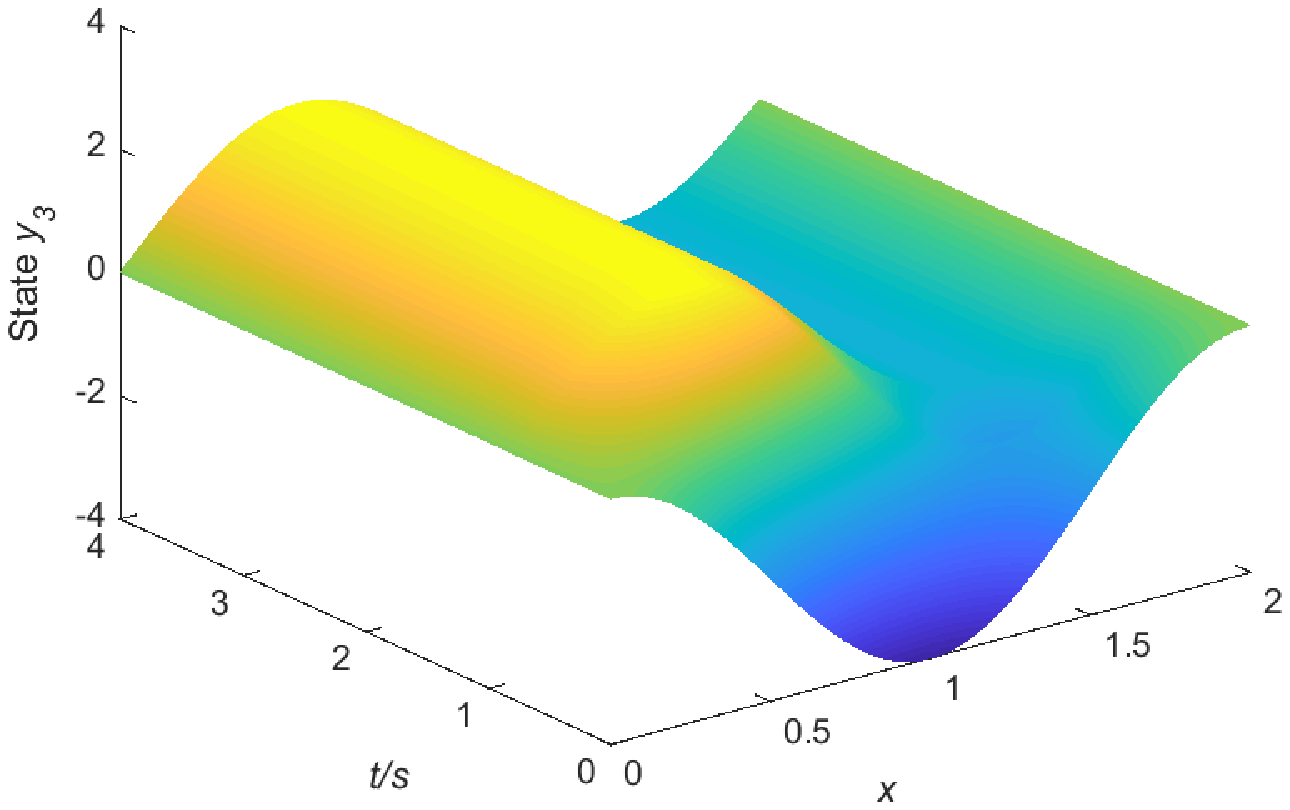}
   \caption{The states of MASs \eqref{e4.5} with \eqref{e4.8} via the controller \eqref{e4.6}. }\label{Fig7}
 \end{figure}

\begin{figure}[H]
   \centering
   \includegraphics[width=4cm]{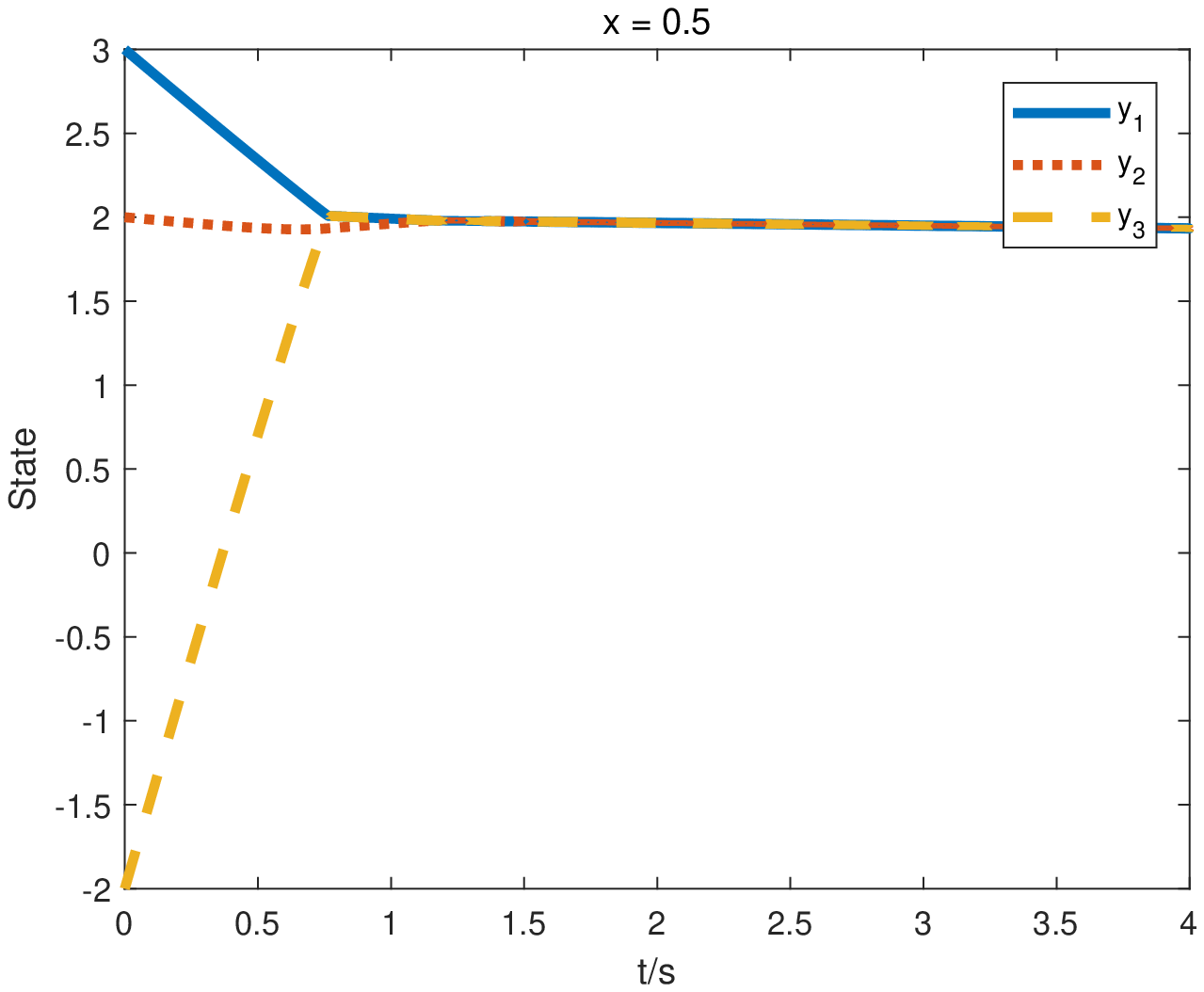}
   \includegraphics[width=4cm]{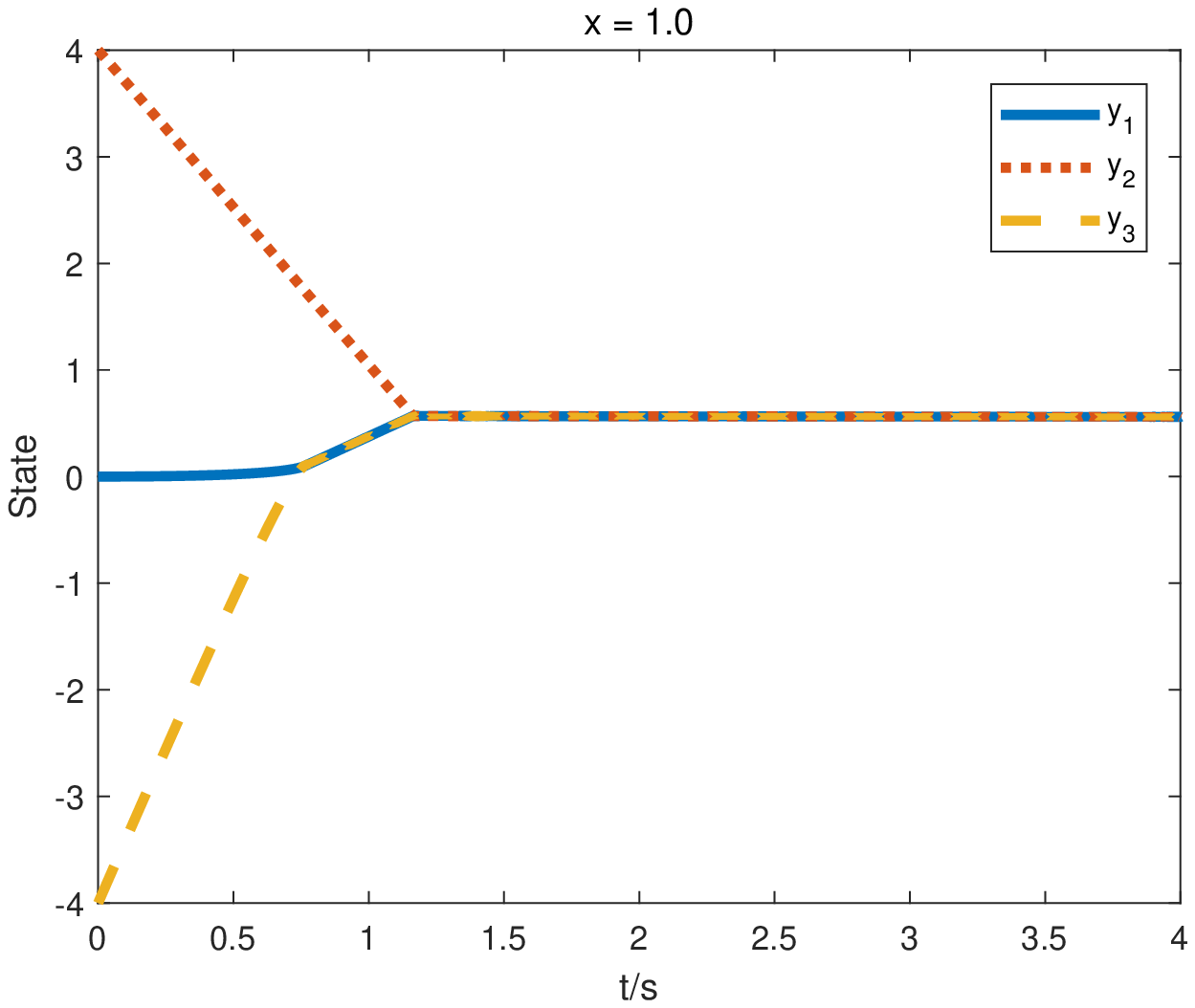}
   \includegraphics[width=4cm]{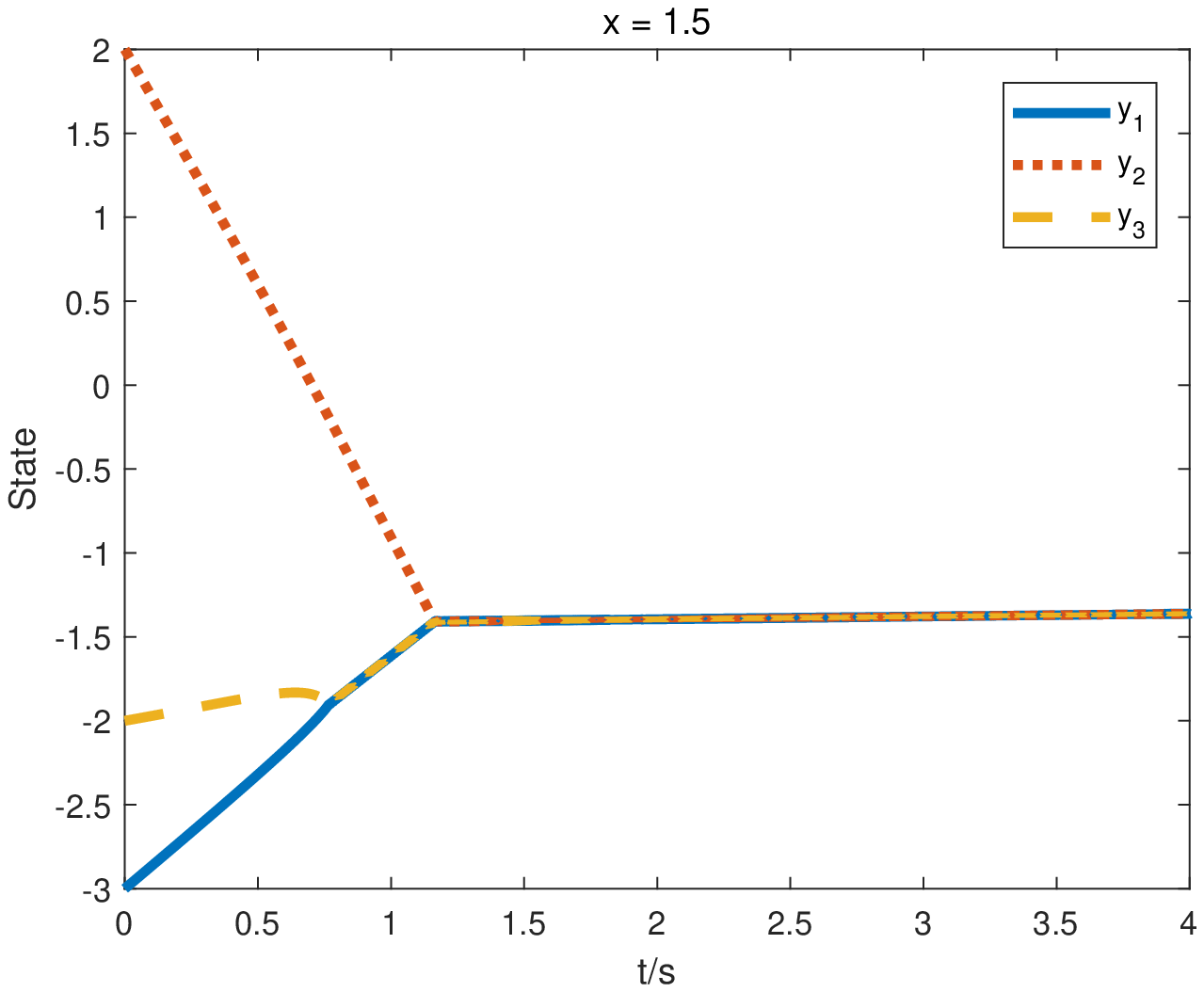}
   \caption{The sections of MASs \eqref{e4.5} with \eqref{e4.8} via the controller \eqref{e4.6}.}\label{Fig8}
 \end{figure}

From Theorem \ref{t3}, it is easy to see that MASs \eqref{e4.5} with \eqref{e4.8} via the controller \eqref{e4.6} can achieve the FTC at $t^*$ and $t^*\leq 4.15$. Fig.\,\ref{Fig7} and Fig.\,\ref{Fig8} show the states and the sections of MASs \eqref{e4.5} with \eqref{e4.8} at $x=0.5, 1.0, 1.5$, respectively.  Fig.\,\ref{Fig21} displays the sections of MASs \eqref{e4.5} without control at $x=0.5, 1.0, 1.5$, respectively. Thus, the controller \eqref{e4.6} can ensure the FTC of MASs \eqref{e4.5} with \eqref{e4.8}.

\begin{figure}[H]
   \centering
   \includegraphics[width=4cm]{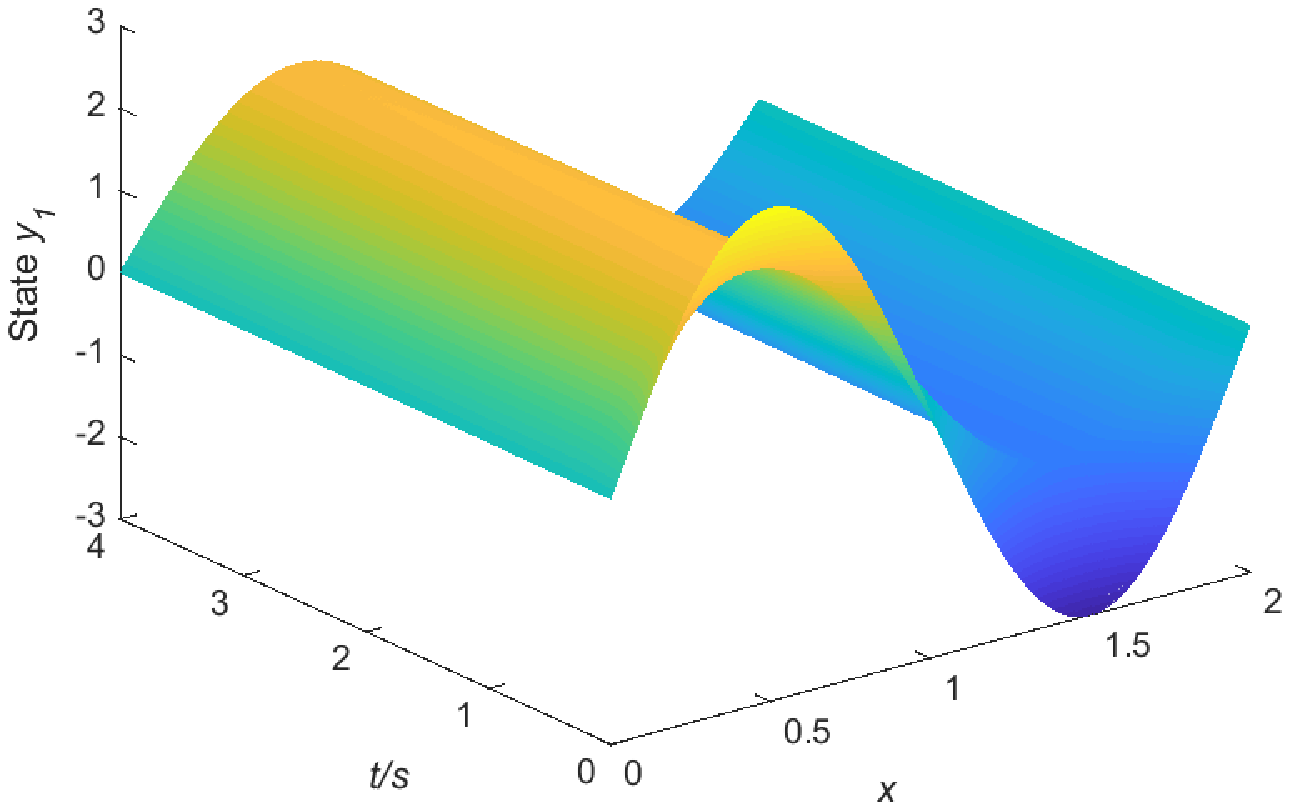}
   \includegraphics[width=4cm]{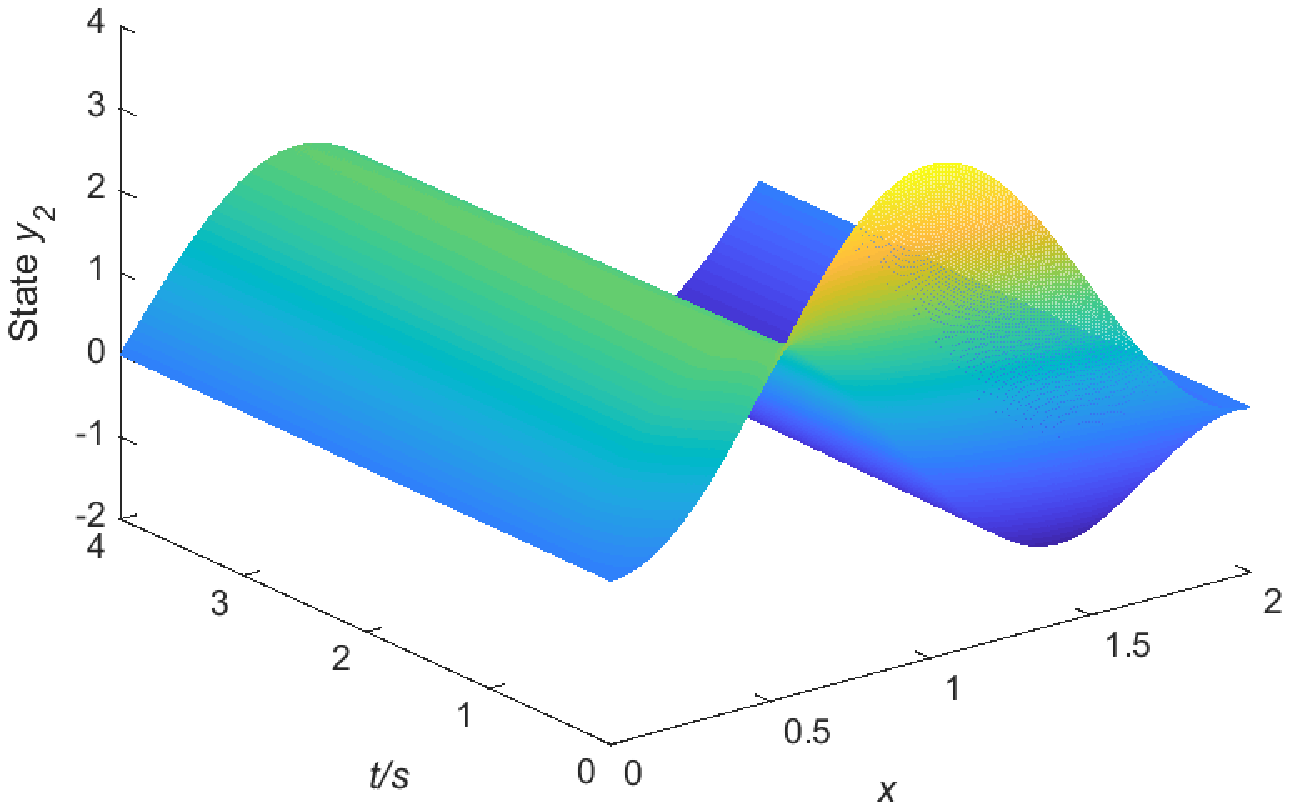}
   \includegraphics[width=4cm]{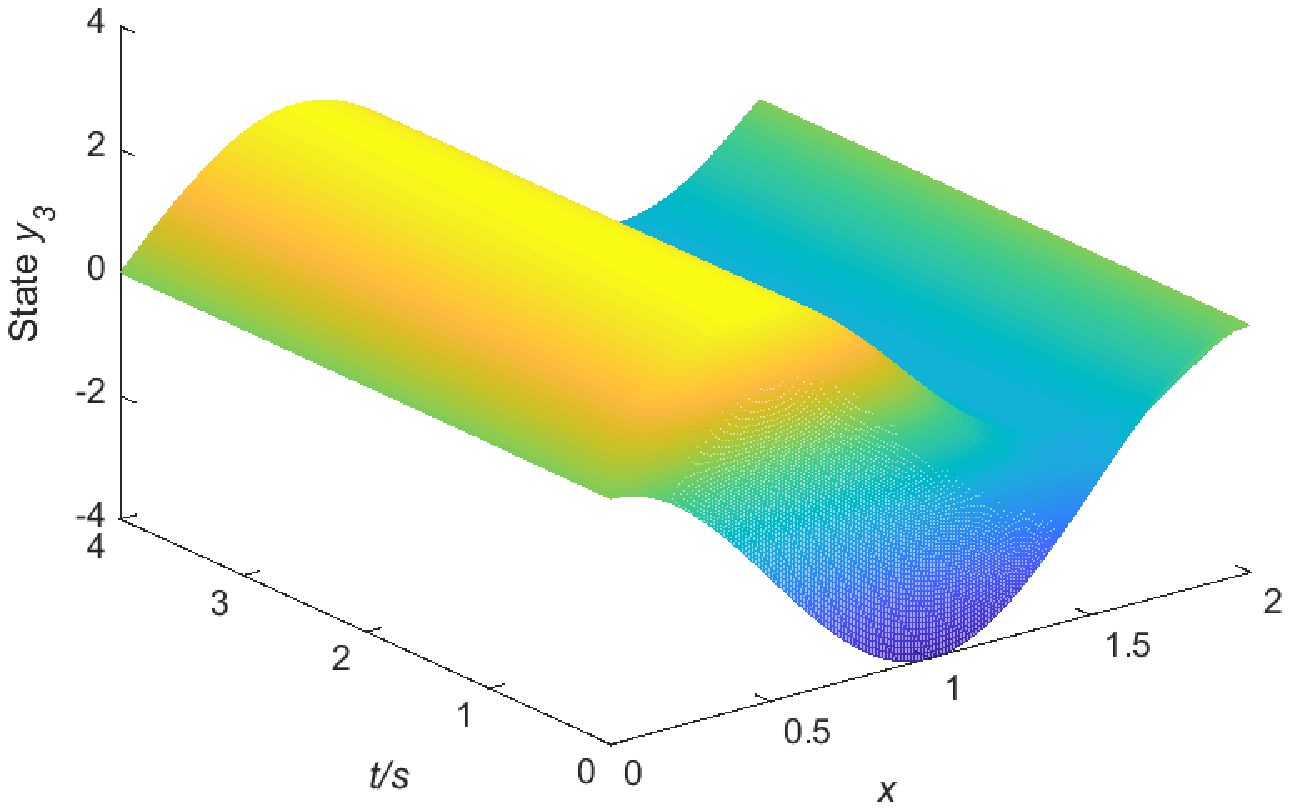}
   \caption{The states of MASs \eqref{e4.5} with \eqref{e4.8} via the controller \eqref{e4.7}.}\label{Fig9}
 \end{figure}

\begin{figure}[H]
   \centering
   \includegraphics[width=4cm]{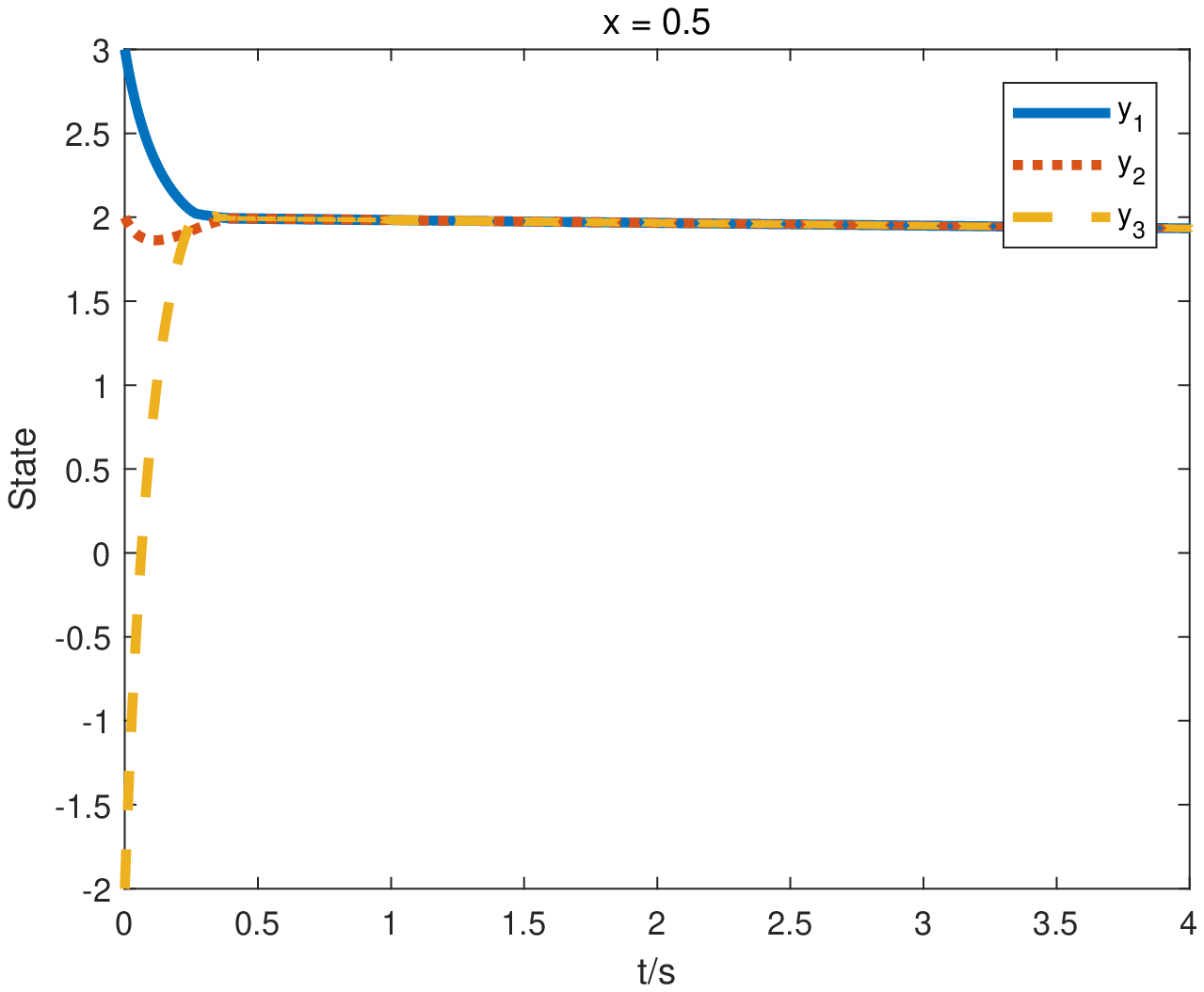}
   \includegraphics[width=4cm]{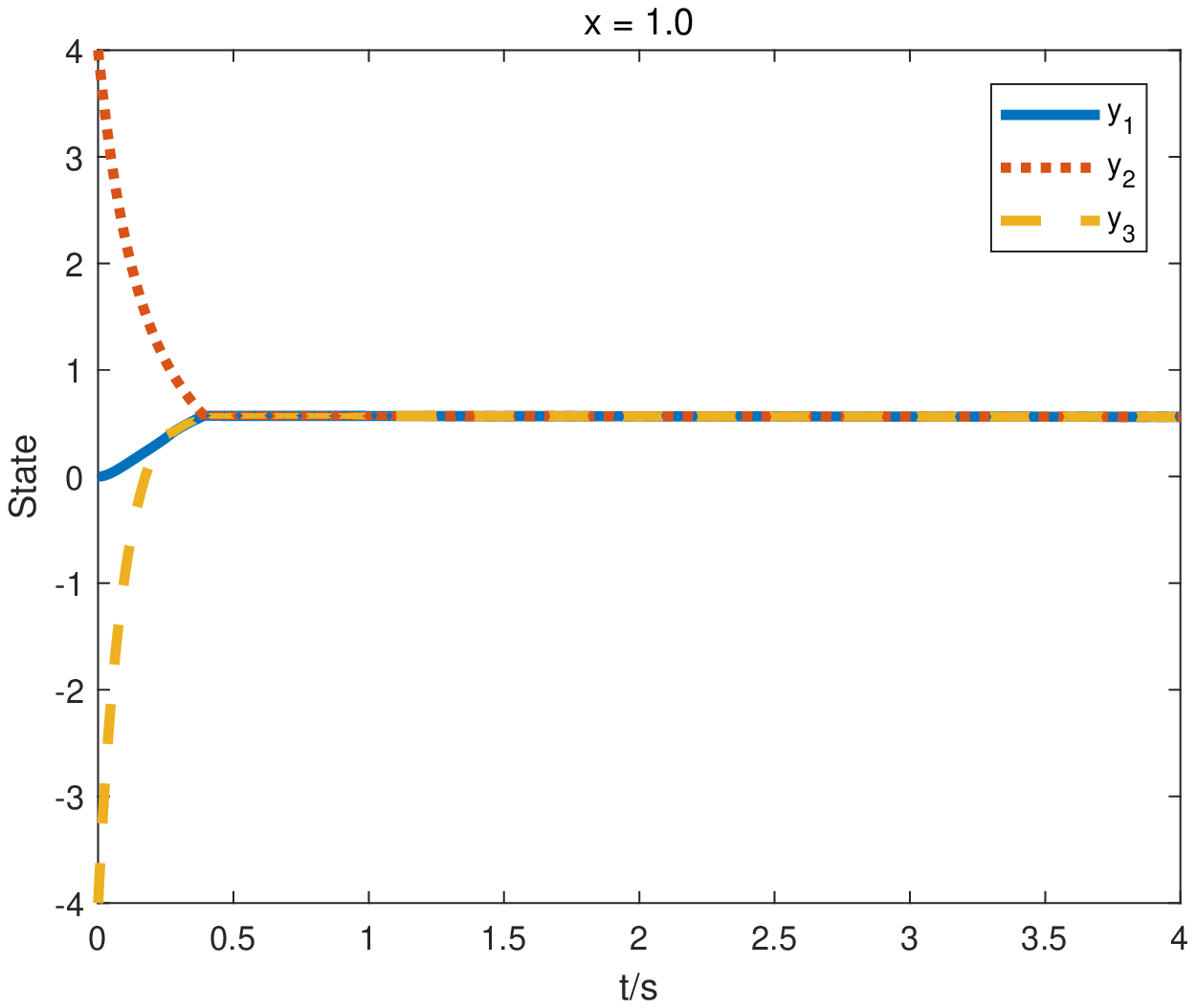}
   \includegraphics[width=4cm]{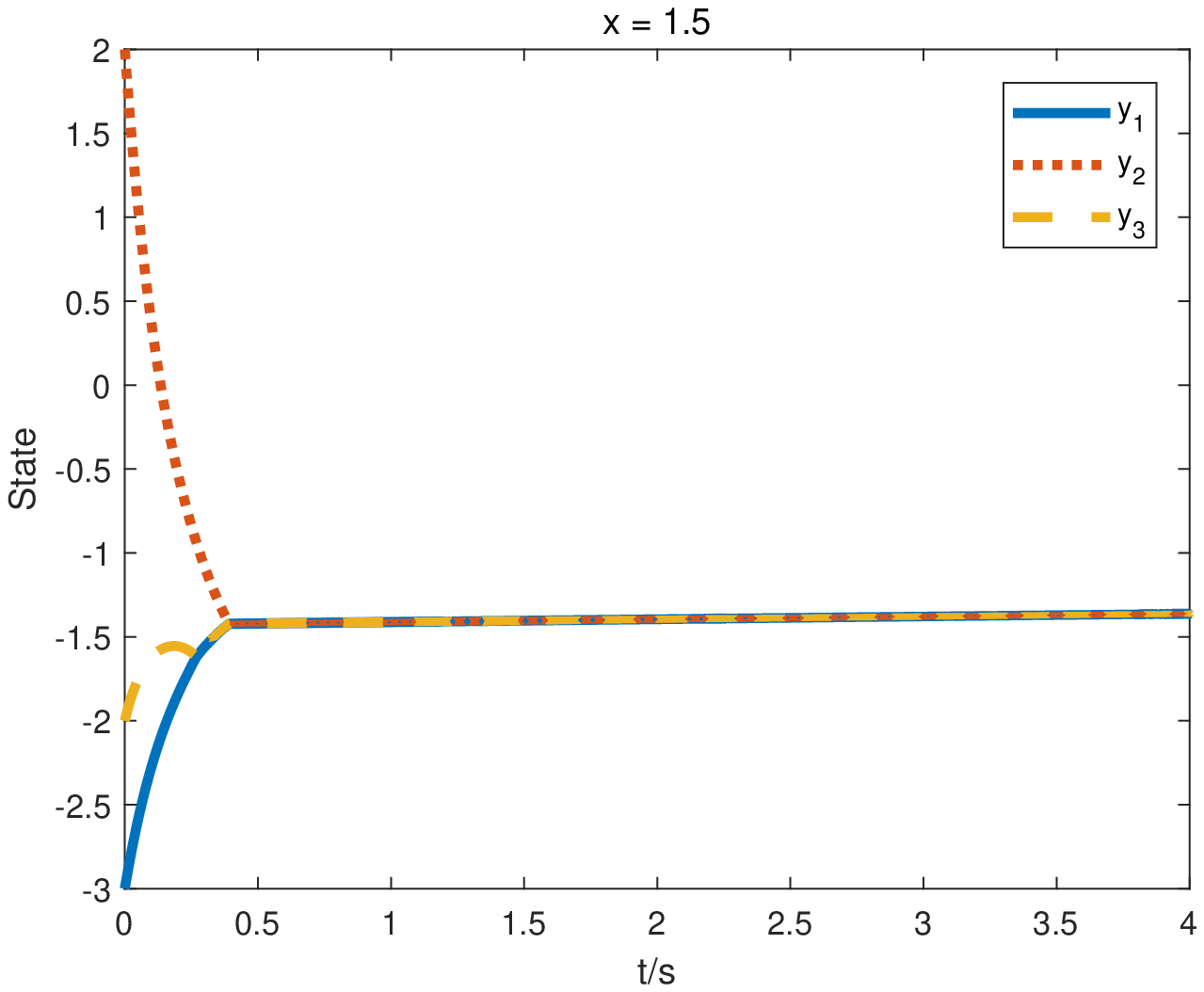}
   \caption{The sections of MASs \eqref{e4.5} with \eqref{e4.8} via the controller \eqref{e4.7}.}\label{Fig10}
 \end{figure}

\begin{figure}[H]
   \centering
   \includegraphics[width=4cm]{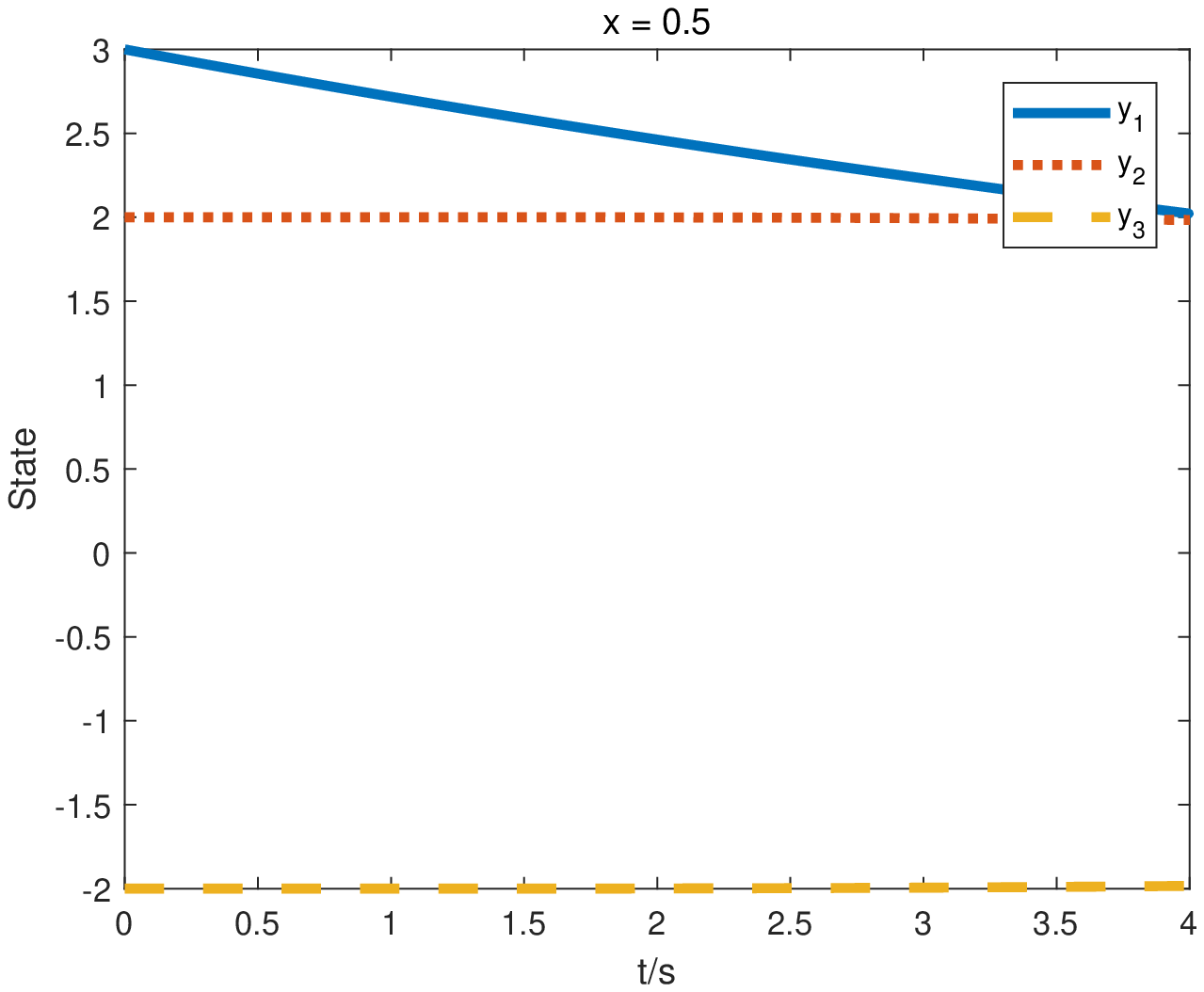}
   \includegraphics[width=4cm]{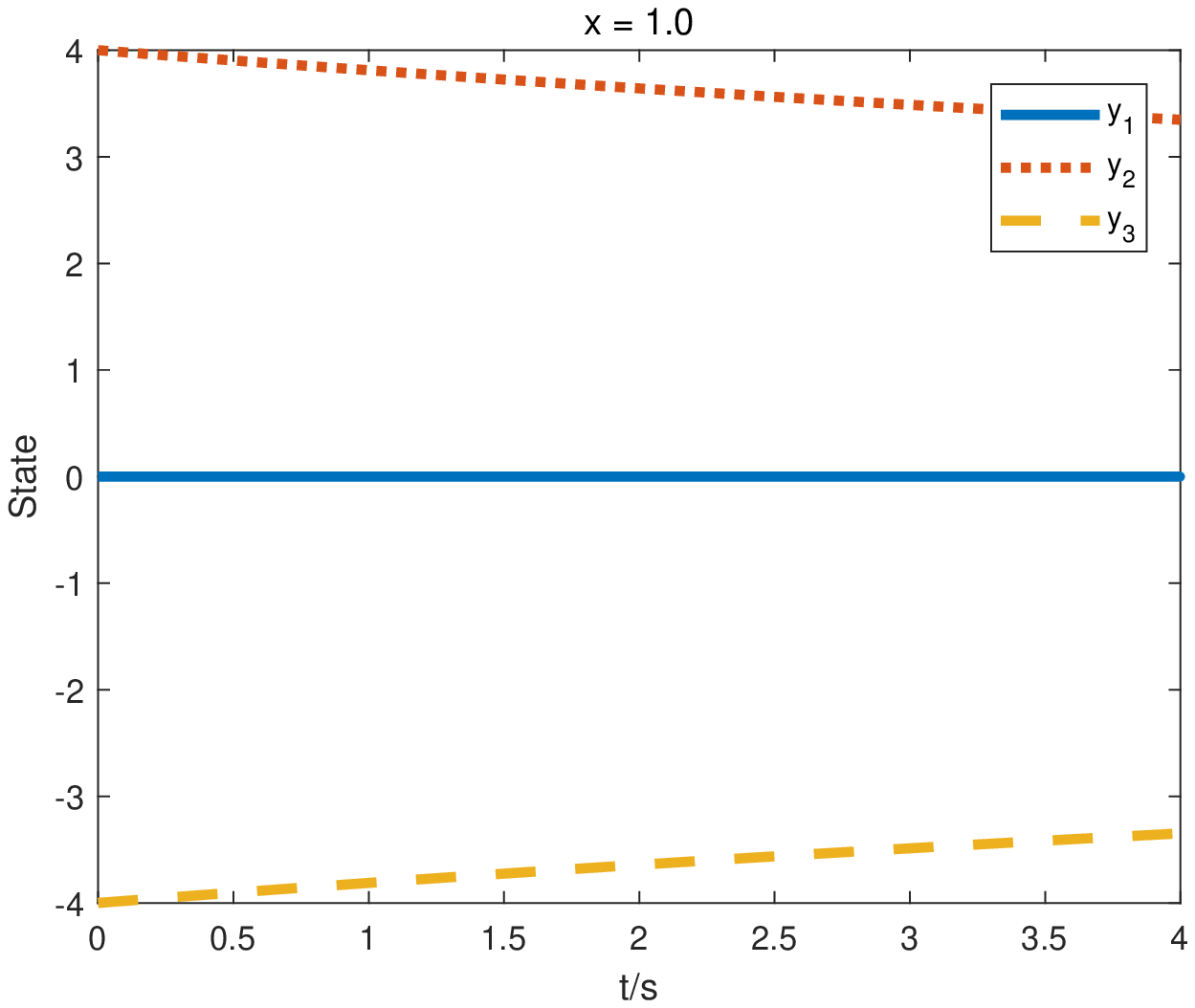}
   \includegraphics[width=4cm]{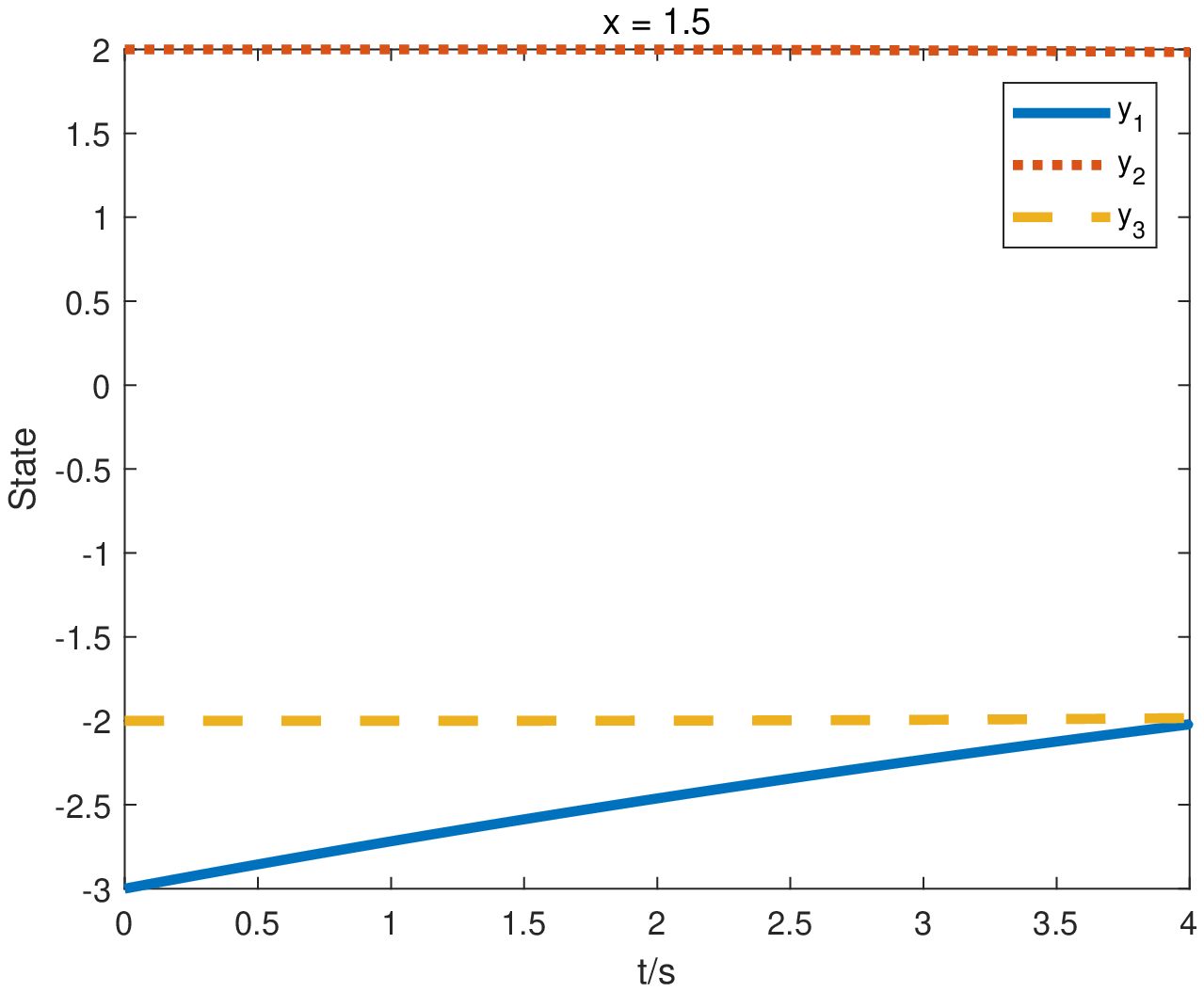}
   \caption{The sections of MASs \eqref{e4.5} with \eqref{e4.8} without control.}\label{Fig21}
 \end{figure}

It follows from Theorem \ref{t4} that MASs \eqref{e4.5} with \eqref{e4.8} via the controller \eqref{e4.7} can reach the FXC at $t^*$ and $t^*\leq T_{\max} = 3.92$. Fig.\,\ref{Fig9} and Fig.\,\ref{Fig10} indicate the states and the sections of MASs \eqref{e4.5} with \eqref{e4.8} at $x=0.5, 1.0, 1.5$, respectively. Therefore, the controller \eqref{e4.7} can guarantee the FXC of MASs \eqref{e4.5} with \eqref{e4.8}.

\begin{figure}[H]
   \centering
   \includegraphics[width=4cm]{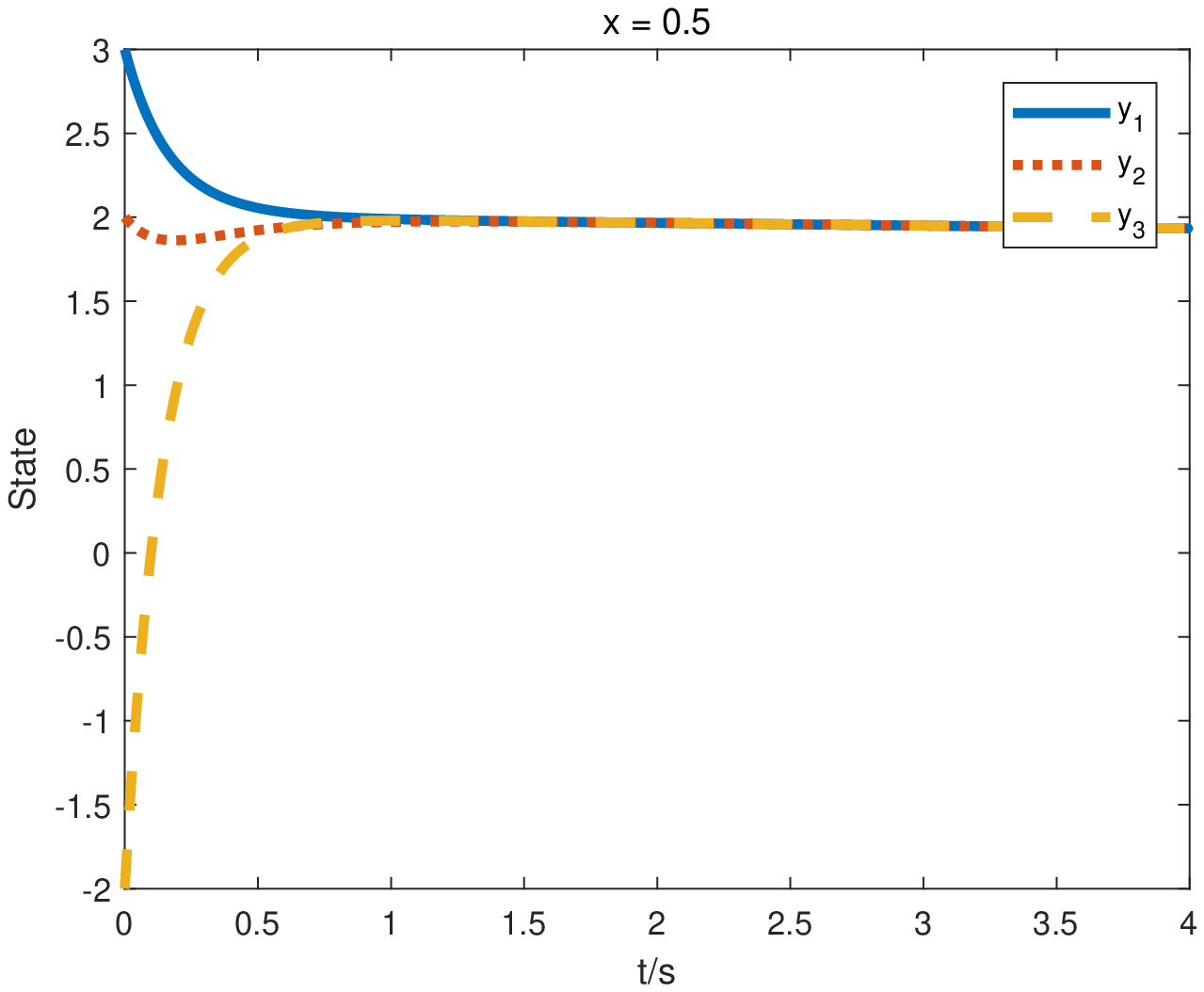}
   \includegraphics[width=4cm]{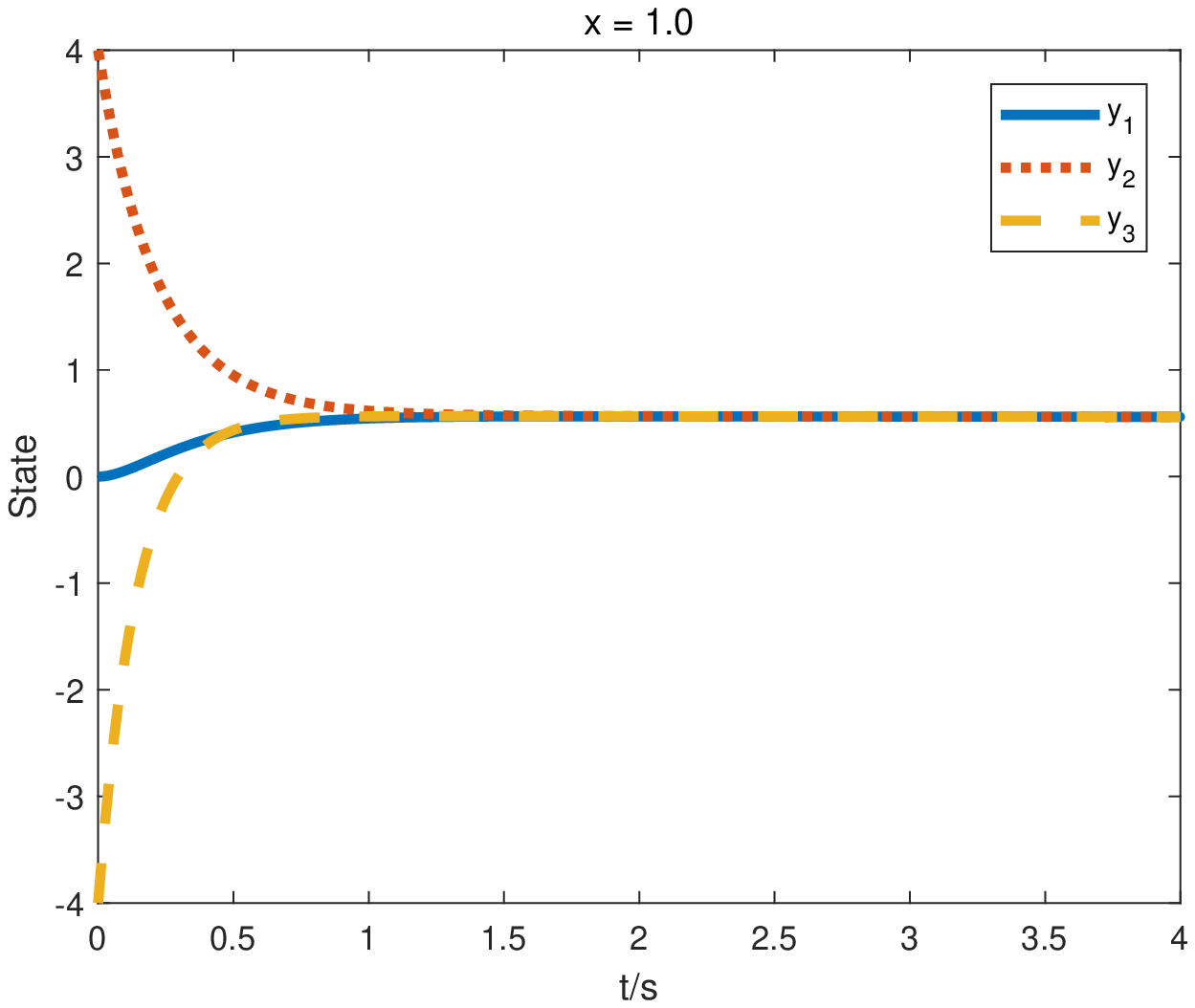}
   \includegraphics[width=4cm]{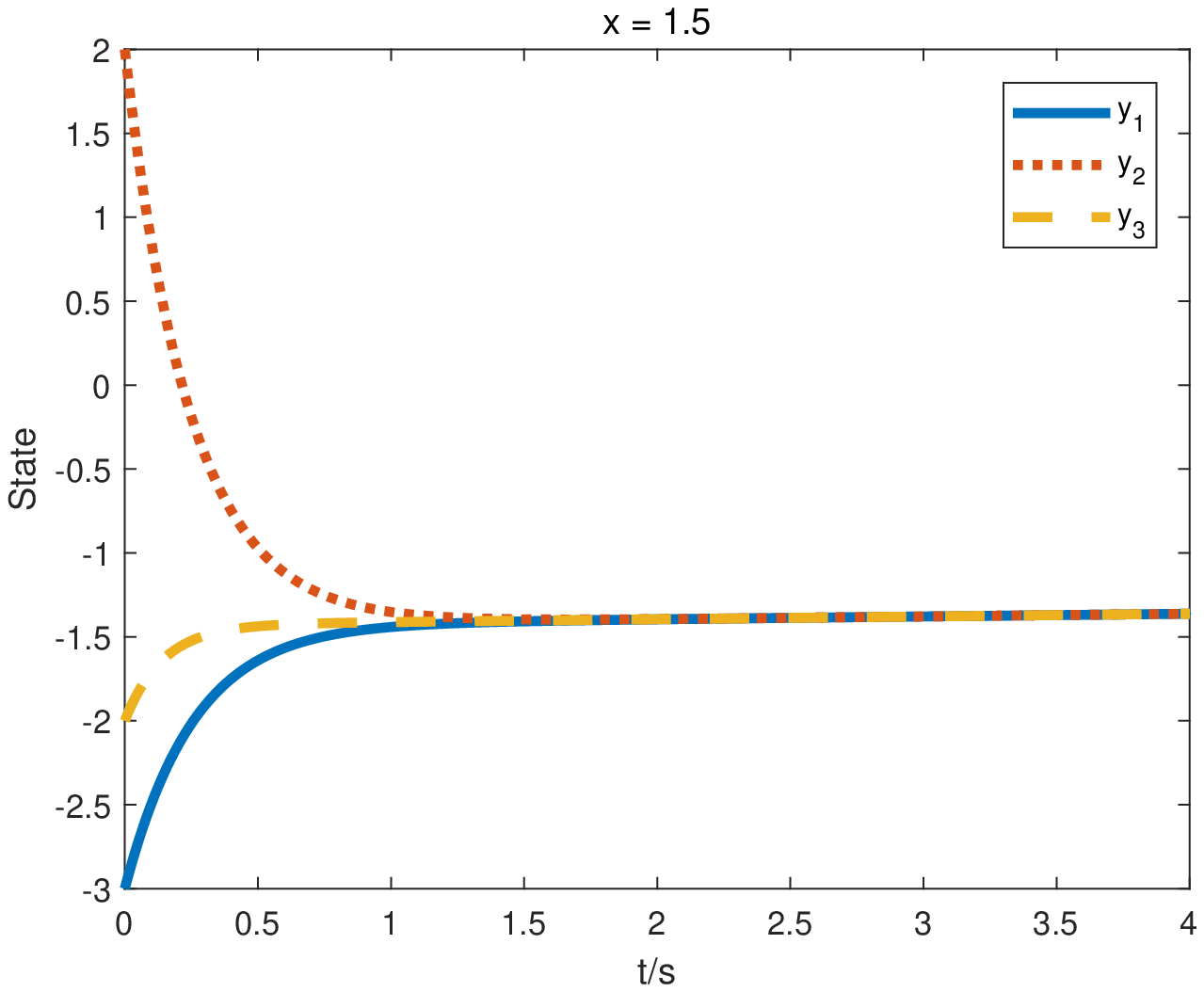}
   \caption{The sections of MASs \eqref{e4.5} with \eqref{e4.8} via the controller \eqref{e41}.}\label{Fig22}
 \end{figure}

\begin{figure}[H]
   \centering
   \includegraphics[width=4cm]{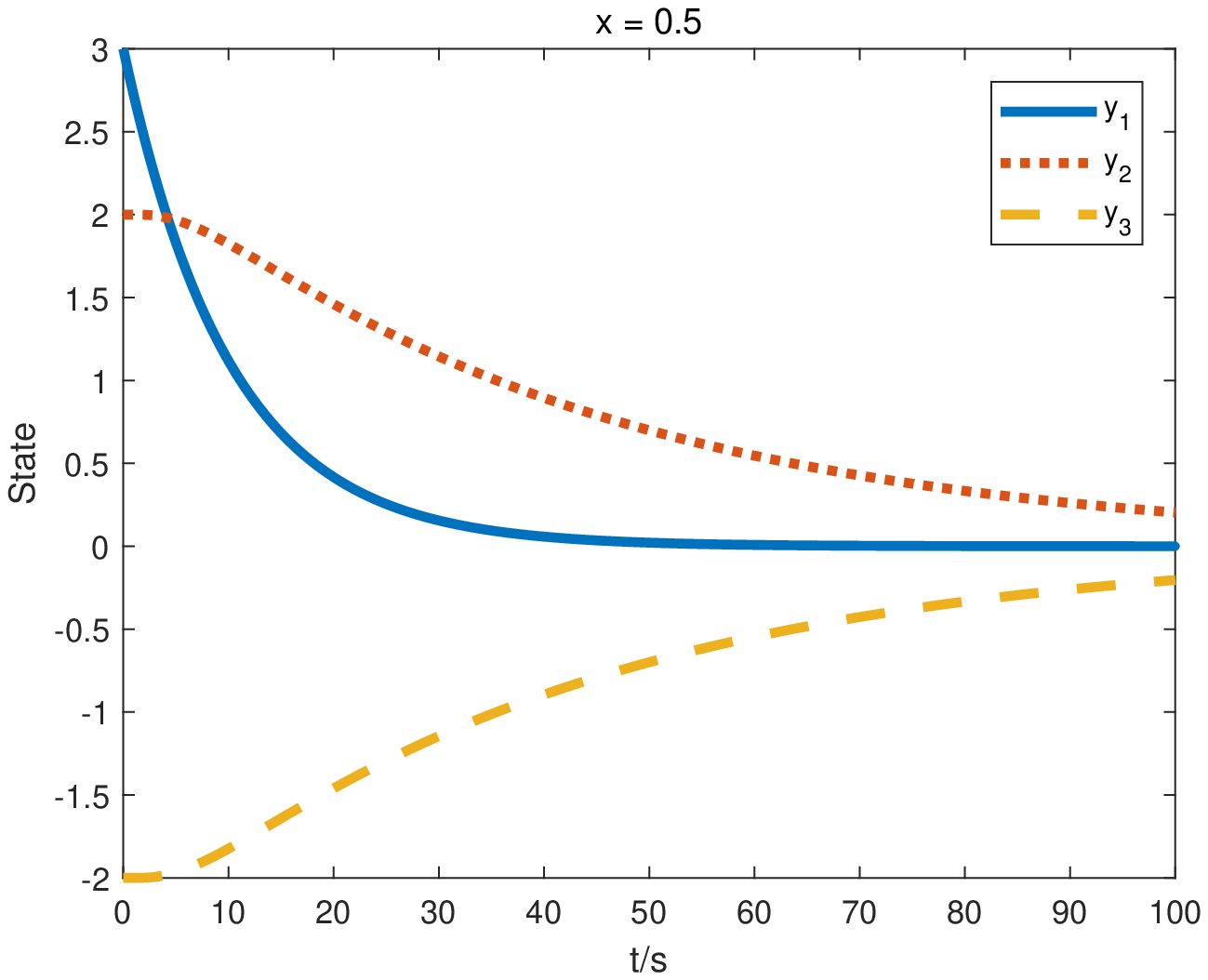}
   \includegraphics[width=4cm]{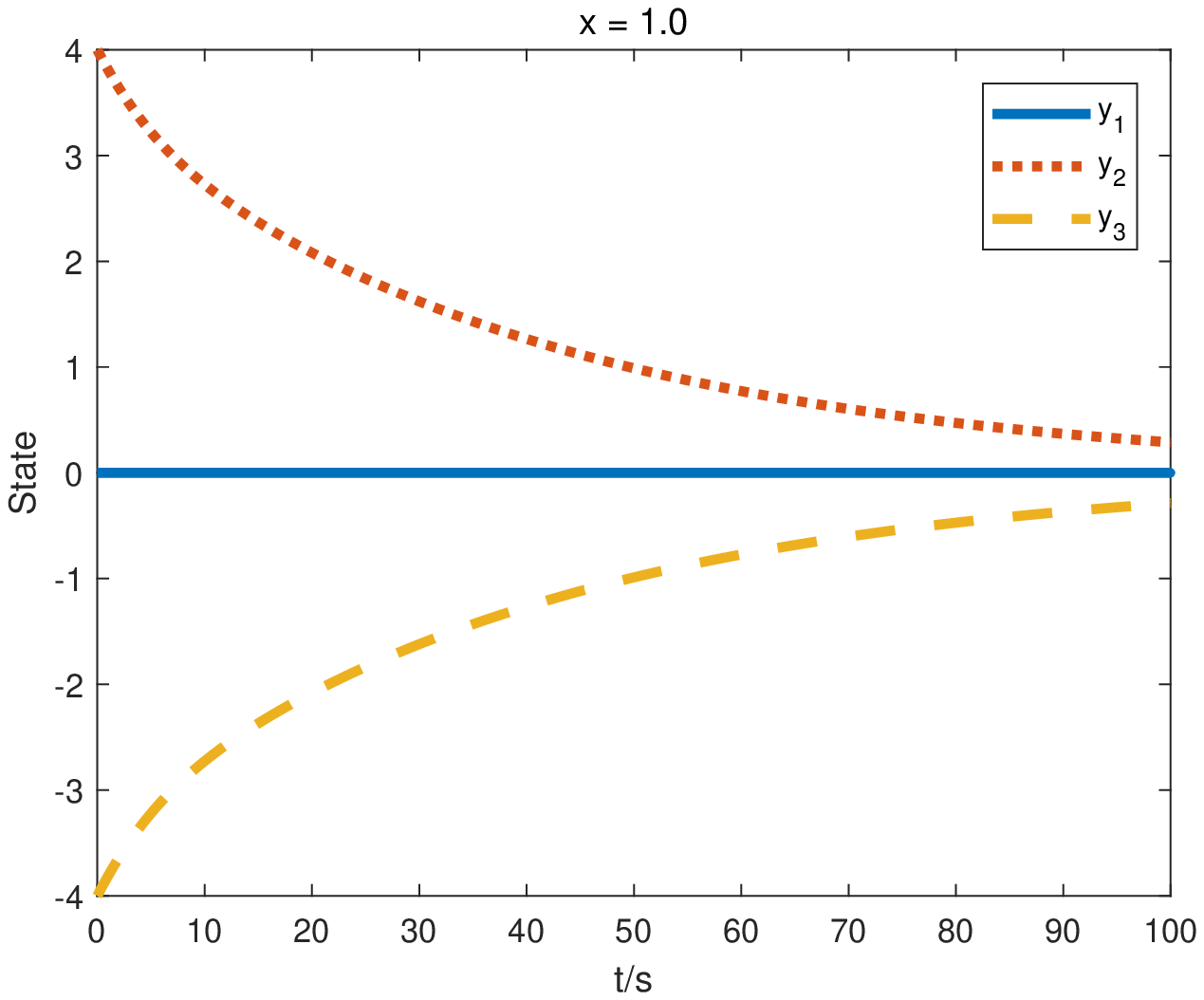}
   \includegraphics[width=4cm]{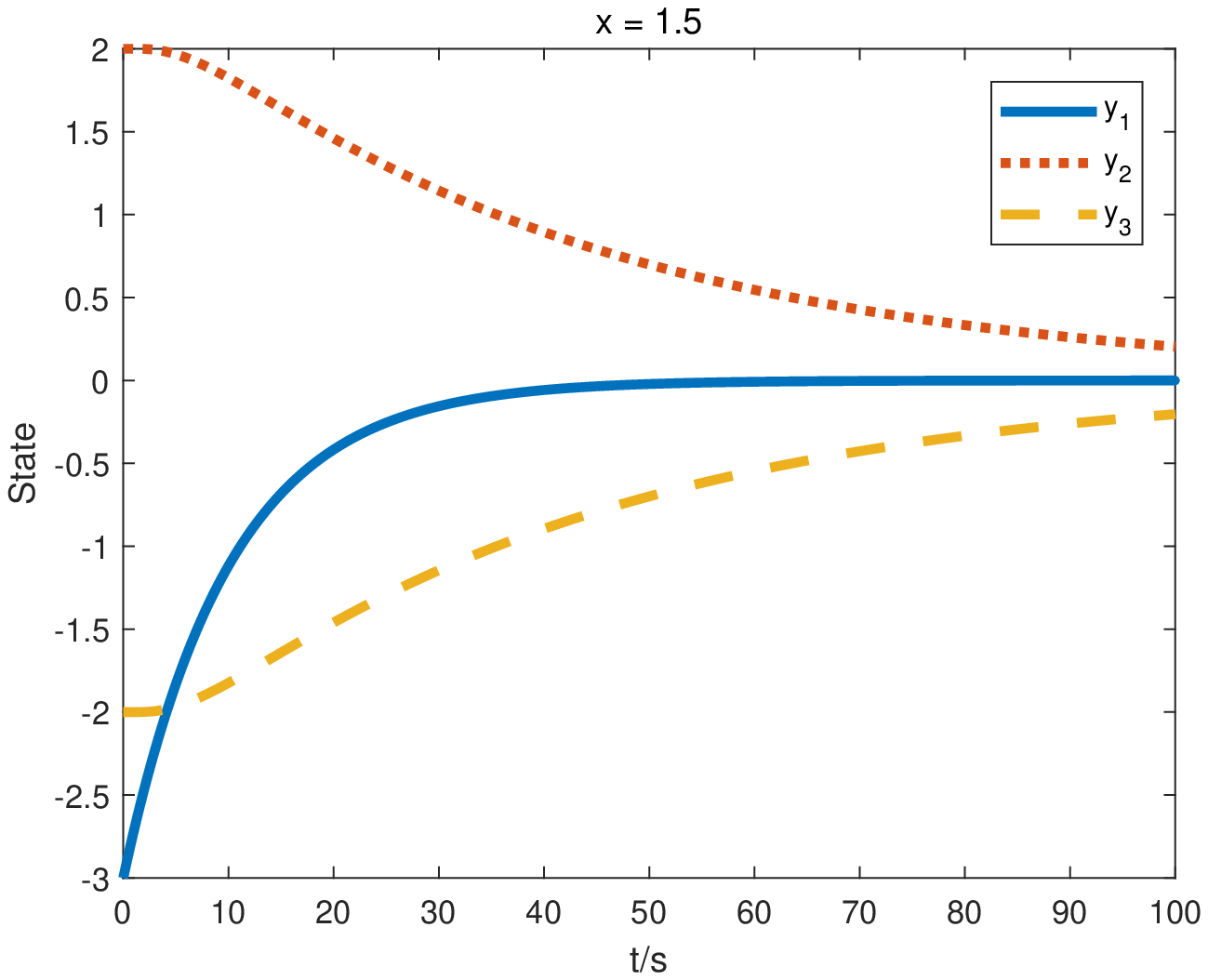}
   \caption{The sections of MASs \eqref{e4.5} with \eqref{e4.8} via the controller \eqref{e42}.}\label{Fig23}
 \end{figure}

The sections of MASs \eqref{e4.5} with \eqref{e4.8} via the controllers \eqref{e41} and \eqref{e42} at $x=0.5, 1.0, 1.5$ are shown in Fig.\,\ref{Fig22} and Fig.\,\ref{Fig23}, respectively. By Fig.\,\ref{Fig8}, Fig.\,\ref{Fig10}, Fig.\,\ref{Fig22} and Fig.\,\ref{Fig23}, we can see that the controllers \eqref{e4.6} and \eqref{e4.7} have faster convergence speed than the controllers \eqref{e41} and \eqref{e42}. Moreover, Fig.\,\ref{Fig8} and Fig.\,\ref{Fig22} display that the convergence speed of the controller \eqref{e41} is faster than the one of the controller \eqref{e4.6} in the first half but it decreases in the second half.

\section{Conclusions}
\noindent

This paper studies the FTC and FXC issues of MASs driven by parabolic PDEs. We design the FTC and FXC distributed controllers to satisfy different control effects in MASs with external disturbance. Theoretical analysis shows that agents can reach the FTC and FXC under an undirected and connected graph or a directed, s-con, and d-bal graph. We also simplify controllers to solve the FTC and FXC issues of MASs without external disturbances. This work provides theoretical support for solving the FTC and FXC control problems of MASs described by parabolic PDEs. Furthermore, it gives some technical guidance towards the design of the controller to satisfy different control effects in MASs with/without external disturbance. However, for MASs driven by parabolic PDEs under event-triggering scheme or cyber attacks, how to ensure the reliability of the controller is still an open problem.

Noting that, in some real situations, it is necessary and important to co-consider other control problems in MASs such as the event-triggered control problem, the fuzzy fast-tracking control problem, the cyber attack problem and so on.  Our future work will consider designing finite-time event-triggered and observer-based control strategies to solve the FTC and FXC problems of MASs driven by more general PDEs.



\section*{Acknowledgements}
\noindent

The original version of this paper was submitted to IEEE Transactions on Systems, Man, and Cybernetics: Systems on June 19 of 2022. The authors are grateful to the editors and reviewers whose helpful comments and suggestions have improved the paper.


\begin{thebibliography}{1}
\bibliographystyle{IEEEtran}

\bibitem{Qin16} J. Qin, W. Fu, H. Gao and W. Zheng, ''Distributed k-means algorithm and fuzzy c-means algorithm for sensor networks based on multiagent consensus theory," \textit{ IEEE Transactions on Cybernetics}, vol. 47, no. 3, pp. 772--783, 2016.


\bibitem{Avram18} R.C. Avram, X. Zhang and J. Muse, ''Nonlinear adaptive fault-tolerant quadrotor altitude and attitude tracking with multiple actuator faults," \textit{ IEEE Transactions on Control Systems Technology}, vol. 26, no. 2, pp. 701--707, 2018.

\bibitem{Ren09} W. Ren, ''Distributed cooperative attitude synchronization and tracking for multiple rigid bodies," \textit{ IEEE Transactions on Control Systems Technology}, vol. 18, no. 2, pp. 383--392, 2010.

\bibitem{Re87} C.W. Reynolds, ''Flocks, herds and schools: A distributed behavioral model," \textit{ Computer Graphics}, vol. 21, no. 4, pp. 25--34, 1987.


\bibitem{Vi95} T. Vicsek, A. Czirook, E. Ben-Jacob, O. Cohen and I. Shochet, ''Novel type of phase transition in a system of self-driven paticles," \textit{ Physical Review Letters}, vol. 75, no. 6, pp. 1226--1229, 1995.

\bibitem{Ol04} R. Olfati-Saber and R.M. Murray, ''Consensus problems in networks of agents with switching topology and time-delays," \textit{ IEEE Transactions on Automatic Control}, vol. 49, no. 9, pp. 1520--1533, 2004.

\bibitem{LiH17} H. Li, Y. Zhu, J. Wang, J. Liu, S. Shen, H. Gao and Y. Sun, ''Consensus of nonlinear second-order multi-agent systems with mixed time-delays and intermittent communications," \textit{ Neurocomputing}, vol. 251, no. 16, pp. 115--126, 2017.

\bibitem{Hamed17} H. Rezaee and F. Abdollahi, ''Discrete-time consensus strategy for a class of high-order linear multiagent systems under stochastic communication topologies," \textit{ Journal of the Franklin Institute}, vol. 354, no. 9, pp. 3690--3705, 2017.

\bibitem{Wang19} X. Wang, X. Li, N. Huang and D. O'Regan, ''Asymptotical consensus of fractional-order multi-agent systems with current and delay states," \textit{ Applied Mathematics and Mechanics (English Edition)}, vol. 40, no. 11, pp. 1677--1694, 2019.

\bibitem{Li05} Z. Lin, B. Francis and M. Maggiore, ''Necessary and sufficient graphical conditional for formation control of unicycles," \textit{ IEEE Transactions on Automatic Control}, vol. 50, no. 1, pp. 121--127, 2005.

\bibitem{Ol07} R. Olfati-Saber, J.A. Fax and R.M. Murray, ''Consensus and cooperation in networked multi-agent systems," \textit{ IEEE Transactions on Automatic Control}, vol. 95, no. 1, pp. 215--233, 2007.

\bibitem{Re07} W. Ren and E. Atkins, ''Distributed multi-vehicle coordinated control via local information exchange," \textit{ International Journal of Robust and Nonlinear Control}, vol. 17, pp. 1002--1033, 2007.

\bibitem{Re08} W. Ren, ''On consensus algorithms for double-integrator dynamics," \textit{ IEEE Transactions on Automatic Control}, vol. 53, no. 6, pp. 1503--1509, 2008.

\bibitem{Yu09} W. Yu, G. Chen, Z. Wang and W. Yang, ''Distributed consensus filtering in sensor networks," \textit{ IEEE Transactions on Systems, Man, and Cybernetics, Part B: Cybernetics}, vol. 39, no. 6, pp. 1568--1577, 2009.

\bibitem{Wu18} Y. Wu, Z. Wang, S. Ding and H. Zhang, ''Leader–follower consensus of multi-agent systems in directed networks with actuator faults," \textit{ Neurocomputing}, vol. 275, pp. 1177--1185, 2018.

\bibitem{Yang13} X. Yang, Z. Wu and J. Cao, ''Finite-time synchronization of complex networks with nonidentical discontinuous nodes," \textit{ Nonlinear Dynamics}, vol. 73, no. 4, pp. 2313--2327, 2013.

\bibitem{Liu15} X. Liu, J. Cao, W. Yu and Q. Song, ''Nonsmooth finite-time synchronization of seitched coupled neural networks," \textit{ IEEE Transactions on Cybernetics}, vol. 46, no. 10, pp. 2360--2371, 2015.

\bibitem{WangL19} L. Wang and T. Chen, ''Finite-time and fixed-time anti-synchronization of neural networks with time-varying delays," \textit{ Neurocomputing}, vol. 329, pp. 165--171, 2019.

\bibitem{Bh00} S.P. Bhat and D.S. Bernstein, ''Finite-time stability of continuous autonomous systems," \textit{ SIAM Journal on Control and Optimization}, vol. 38, no. 3, pp. 751--766, 2000.

\bibitem{Co06}J. Cortes, ''Finite-time convergent gradient flows with applications to network consensus," \textit{ Automatica}, vol. 42, no. 11, pp. 1993--2000, 2006.


\bibitem{Wa10} L. Wang and F. Xiao, ''Finite-time consensus problems for networks of dynamic agents," \textit{ IEEE Transactions on Automatic Control}, vol. 55, no. 4, pp. 950--955, 2010.



\bibitem{Li11} S. Li, H. Du and X. Lin, ''Finite-time consensus algorithm for multi-agent systems with double-integrator dynamics," \textit{ Automatica}, vol. 47, no. 8, pp. 1706--1712, 2011.

\bibitem{Guan12} Z. Guan, F. Sun, Y. Wang and T. Li, ''Finite-time consensus for leader-following second-order multi-agent networks," \textit{ IEEE Transactions on Circuits and Systems I: Regular Papers}, vol. 59, no. 11, pp. 2646--2654, 2012.

\bibitem{Fu16} J. Fu and J. Wang, ''Observer-based finite-time coordinated tracking for general linear multi-agent systems," \textit{ Automatica}, vol. 66, pp. 231--237, 2016.

\bibitem{Duan19} J. Duan, H. Zhang, Y. Liang and Y. Cai, ''Bipartite finite-time output consensus of heterogeneous multi-agent systems by finite-time event-triggered observer," \textit{ Neurocomputing}, vol. 365, pp. 86--93, 2019.

\bibitem{Po12} A. Polyakov, ''Nonlenear feedback design for fixed-time stabilization of linear control systems," \textit{ IEEE Transactions on Automatic Control}, vol. 57, no. 8, pp. 2106--2110, 2012.

\bibitem{Pa13} S.E. Parsegov, A. Polyakov and P.S. Shcherbakov, ''Fixed-time consensus algorithm for multi-agent systems with integrator dynamics," \textit{ IFAC Proceedings Volumes}, vol. 46, no. 27, pp. 110--115, 2013.

\bibitem{De15} M. Defoort, A. Polyakov, G. Demesure, M. Djemai and K. Veluvolu, ''Leader-follower fixed-time consensus for multi-agent systems with unknown non-linear inherent dynamics," \textit{ IET Control Theory and Applications}, vol. 9, no. 14, pp. 2165--2170, 2015.

\bibitem{Zuo15} Z. Zuo, ''Nonsingular fixed-time consensus tracking for second-order multi-agent networks," \textit{ Automatica}, vol. 54, pp. 305--309, 2015.

\bibitem{Zuo16} Z. Zuo and L. Tie, ''Distributed robust finite-time nonlinear consensus protocols for multi-agent systems," \textit{ International Journal of Systems Science}, vol. 47, no. 6, pp. 1366--1375, 2016.

\bibitem{Zuo18} Z. Zuo, Q. Han, B. Ning, X. Ge and X. Zhang, ''An overview of recent advances in fixed-time cooperative control of multiagent systems," \textit{ IEEE Transactions on Industrial Informatics}, vol. 14, no. 6, pp. 2322--2334, 2018.


\bibitem{Zh12} Y. Zheng and L. Wang, ''Finite-time consensus of heterogeneous multi-agent systems with and without velocity measurement," \textit{ Systems Control Letters}, vol. 61, no. 8, pp. 871--878, 2012.

\bibitem{Zh14} Y. Zhu, X. Guan and X. Luo, ''Finite-time consensus of heterogeneous multi-agent systems with linear and nonlinear dynamics," \textit{ Acta Automatica Sinica}, vol. 40, no. 11, pp. 2618--2624, 2014.

\bibitem{Xi14} F. Xiao, L. Wang and T. Chen, ''Finite-time consensus in  networks of integrator-like dynamic agents with directional link failure," \textit{ IEEE Transactions on Automatic Control}, vol. 59, no. 3, pp. 756--762, 2014.

\bibitem{Li17} X. Lin and Y. Zheng, ''Finite-time consensus of switched multiagent systems," \textit{ IEEE Transactions on Systems, Man, and Cybernetics: Systems}, vol. 47, no. 7, pp. 1535--1545, 2017.

\bibitem{Tong18} P. Tong, S. Chen and L. Wang, ''Finite-time consensus of multi-agent systems with continuous time-varying interaction topology," \textit{ Neurocomputing}, vol. 284, pp. 187--193, 2018.

\bibitem{Hu19} B. Hu, Z. Guan and M. Fu, ''Distributed event-driven control for finite-time consensus," \textit{ Automatica}, vol. 103, pp. 88--95, 2019.

\bibitem{Sa19} R. Sakthivel, S. Kanakalakshmi, B. Kaviarasan, Y. Ma and A. Leelamani, ''Finite-time consensus of input delayed multi-agent systems via non-fragile controller subject to switching topology," \textit{ Neurocomputing}, vol. 325, pp. 225--233, 2019.

\bibitem{Liu19} W. Liu, Q. Ma, Q. Wang and H. Feng, ''Finite-time consensus control of heterogeneous nonlinear MASs with uncertainties bounded by positive functions," \textit{ Neurocomputing}, vol. 330, pp. 29--38, 2019.

\bibitem{LiuJ19} J. Liu, C. Wang and X. Cai, ''Global finite-time event-triggered consensus for a class of second-order multi-agent systems with the power of positive odd rational number and quantized control inputs," \textit{ Neurocomputing}, vol. 360, pp. 254--264, 2019.

\bibitem{Li22}P. Li, X. Wu, X. Chen and J. Qiu, ''Distributed adaptive finite-time tracking for multi-agent systems and its application," \textit{ Neurocomputing}, vol. 481, pp. 46--54, 2022.


\bibitem{Ghods12} N. Ghods and M. Krstic, ''Multiagent deployment over a source," \textit{ IEEE Transactions on Control Systems Technology }, vol. 20, no. 1, pp. 277--285, 2012.

\bibitem{Ras19} H. Rastgoftar and E. Atkins, ''Physics-based freely scalable continuum deformation for uas traffic coordination," \textit{ IEEE Transactions on Control of Network Systems}, vol. 7, no. 2 pp. 532--544, 2019.

\bibitem{Cas19}C.G. Cassandras, ''Freeway traffic modelling and control," \textit{ IEEE Control Systems Magazine }, vol. 39, no. 1 pp. 68--71, 2022.


\bibitem{Wa11} J. Wang, H. Wu and H. Li, ''Distributed fuzzy control design of nonlinear hyperbolic pde systems with application to nonisothermal plug-flow reactor," \textit{ IEEE Transactions on Fuzzy Systems}, vol. 19, no. 3, pp. 514--526, 2011.

\bibitem{Wa12} J. Wang, H. Wu and H. Li, ''Distributed proportional-spatial derivative control of nonlinear parabolic systems via fuzzy pde modelling approach," \textit{IEEE Transactions on Systems Man $\&$ Cybernetics Part B}, vol. 42, no. 3, pp. 927--938, 2012.

\bibitem{Ra81} W. Ray, {\it{Advanced Process Control}}, New York, McGraw-Hill, 1981.




\bibitem{Pi16} A. Pilloni, A. Pisano, Y. Orlov and E. Usai, ''Consensus-based conrol for a network of diffusion PDEs with Boundary local interaction," \textit{ IEEE Transactions on Automatic Control}, vol. 61, no. 9, pp. 2708--2713, 2016.



\bibitem{Ya17} C. Yang, H. He, T. Huang, A. Zhang, J. Qiu, J. Cao and X. Li, ''Consnesus for non-linear multi-agent systems modelled by PDEs based on spatial boundary communication," \textit{ IET Control Theory and Applications}, vol.11, no. 17, pp. 3196--3200, 2017.

\bibitem{Ya18} C. Yang, T. Huang, A. Zhang, J. Qiu, J. Cao and F.E. Alsaadi, ''Output consensus of multiagent systems based on PDEs with input constraint: a boundary control approach," \textit{ IEEE Transactions on Systems, Man, and Cybernetics: Systems}, vol. 51, no. 1, pp. 370--377, 2018.

\bibitem{Wang21}X. Wang and N. Huang, ''Finite-time consensus of multi-agent systems driven by hyperbolic partial differential equations via boundary control," \textit{ Applied Mathematics and Mechanics (English Edition)}, vol. 42, no. 12, pp. 1799--1816, 2021.

\bibitem{Wu13}J. Wang and H. Wu, ''Passivity of delayed reaction-diffusion networks with application to a food web model," \textit{ Applied Mathematics and Computation}, vol. 219, pp. 11311--11326, 2013.

\bibitem{Shi20}B. Shi, J. Yuan and Y. Jian, {\it{Fundamentals of Systems Stability and Control (in Chinese)}}, Beijing, Publishing House of Electronics Industry, 2020.

\bibitem{Wu14} H. Wu, J. Wang and H. Li, ''Fuzzy boundary control design for a class of nonlinear parabolic distributed parameter systems," \textit{ IEEE Transactions on Fuzzy Systems}, vol. 22, no. 3, pp. 642--652, 2014.

\bibitem{Ev81} L.C. Evans, {\it{Partial differential equations}}, American Mathematical Society, 2010.

\bibitem{H087} S.L. Hollis, R.H. Martin, Wu and M. Pierre, ''Global existence and boundedness in reaction-diffusion systems," \textit{ SIAM Journal on Mathematical Analysis}, vol. 18, no. 3, pp. 744--761, 1987.


\bibitem{Fu18} Q. Fu, L. Du, G. Xu, J. Wu and P. Yu, ''Consensus control for multi-agent systems with distributed parameter models," \textit{ Neurocomputing}, vol. 365, pp. 86--93, 2019.

\end{thebibliography}
\end{document}